\documentclass[12 pt]{amsart}
\usepackage[legalpaper, margin=1in]{geometry} 

\usepackage{lipsum}
\usepackage{tikz-cd}
\usepackage[unicode]{hyperref}
\usepackage{graphicx,color}
\usepackage{amssymb,amsmath}
\usepackage{stix}
\usepackage{mathtools}
\usepackage{mathrsfs}
\usepackage{enumerate}
\usepackage{ragged2e}
\usepackage{enumitem,pifont}
\usepackage{xypic}
\usepackage{tikz}
\usetikzlibrary{arrows,decorations.pathmorphing,automata,backgrounds}
\usetikzlibrary{backgrounds,positioning}

\usepackage{caption}
\usepackage{subcaption}
\usepackage{float}

\usepackage{egothic}
\usepackage[T1]{fontenc}
\usepackage{yfonts}

\graphicspath{{Figures/}}

\usepackage{calc}

\newcommand{\bE}{\mathbf{E}}
\newcommand{\bQ}{\mathbf{Q}}

\newcommand{\calI}{\mathcal{I}}

\newcommand{\bG}{\mathbf{G}}
\newcommand{\calS}{\mathcal{S}}
\newcommand{\calM}{\mathcal{M}}

\newcommand{\bP}{\mathbf{P}}
\newcommand{\bC}{\mathbf{C}}
\newcommand{\bD}{\mathbf{D}}
\newcommand{\bS}{\mathbf{S}}

\newcommand{\fC}{\mathfrak{C}}

\newcommand{\nb}{\operatorname{nb}}

\newcommand{\colim}{\operatorname{colim}}

\newcommand{\End}{\operatorname{End}}
\newcommand{\Emb}{\operatorname{Emb}}
\newcommand{\map}{\mathbb{R}\operatorname{Map}}
\newcommand{\Gal}{\operatorname{Gal}}

\newcommand{\Aut}{\operatorname{Aut}}

\newcommand{\col}{\operatorname{col}}

\newcounter{dummy} \numberwithin{dummy}{section}
\theoremstyle{definition}
\newtheorem{theorem}[dummy]{Theorem}

\newtheorem{prop}[dummy]{Proposition}

\newtheorem{cor}[dummy]{Corollary}
\newtheorem*{thm*}{Theorem}
\newtheorem*{prop*}{Proposition}
\newtheorem*{definition*}{Definition}

\newtheorem{definition}[dummy]{Definition}

\newtheorem{exercise}[dummy]{Exercise}
\newtheorem{problem}[dummy]{Open Problem}

\newtheorem{example}[dummy]{Example}
\newtheorem{remark}[dummy]{Remark}

\numberwithin{equation}{section}

\usepackage{multicol}

\newcommand{\gt}{\widehat{\mathsf{GT}}}

\newcommand{\ns}{\widehat{\mathsf{NS}}}

\newcommand{\ModOp}{\mathbf{ModOp}}

\newcommand{\bU}{\mathbf{U}}

\title{Lecture notes on modular infinity operads and Grothendieck-Teichm\"uller theory}

\author[O. Borghi]{Olivia Borghi}
\address{School of Mathematics and Statistics \\ The University of Melbourne \\ Melbourne, Victoria, Australia}
\email{oborghi@student.unimelb.edu.au}

\author[M. Robertson]{Marcy Robertson}
\address{School of Mathematics and Statistics \\ The University of Melbourne \\ Melbourne, Victoria, Australia}
\email{marcy.robertson@unimelb.edu.au}

\date{\today}

\begin{document}

\maketitle

\tableofcontents 
\section{Introduction} 
These notes represent the transcript of three, $90$ minute lectures given by the second author at the CRM in Barcelona in 2021 as part of the ``Higher Structures and Operadic Calculus'' workshop. The goal of the series was to introduce and motivate modular $\infty$-operads via their application to what is often called ``Grothendieck-Teichm\"uller'' theory. We therefore start by trying to answer \[\text{What is ``Grothendieck-Teichm\"uller'' theory ?}\]

The absolute Galois group of $\mathbb{Q}$, $\Gal(\mathbb{Q})$, is the 
(topological) group of automorphisms of the separable closure 
$\bar{\mathbb{Q}}$ which fix $\mathbb{Q}$. This is a \emph{profinite group} (Definition~\ref{def: profinite group}), which means that, to write down an element $g$ of the group $\Gal(\mathbb{Q})$, one must describe the image of $g$ in \emph{each} of the finite quotients of $\Gal(\mathbb{Q})$. It turns out, however, we do not know all of the finite quotients of $\Gal(\mathbb{Q})$ (see, for example, \cite{serre_notes_galois}).

The ``Grothendieck'' part of Grothendieck-Teichm\"uller theory is an idea, laid out in \emph{Esquisse d'un Programme} (\cite{grothendieck_esquisse}), to use the actions of  $\Gal(\mathbb{Q})$ on the geometric fundamental groups of the stacks $\calM_{g,n}$ to gain insight into $\Gal(\mathbb{Q})$. This idea is inspired by a theorem of Belyi \cite{belyi} which says that there is a faithful action of $\Gal(\mathbb{Q})$ on $\pi^{et}_1(\calM_{0,4})$ and therefore on general $\pi^{et}_1(\calM_{g,n})$.

Grothendieck suggestion that one could consider the collection of fundamental groups, $\pi^{et}_1(\calM_{g,n})$, together with the natural maps between them, as a single object called the ``Teichm\"uller tower''. He conjectured that the automorphisms of this tower could not only be explicitly described, but that the group of such automorphisms may even be equivalent to $\Gal(\mathbb{Q})$. In addition, Grothendieck suggested that the ideal tower should be constructed from the genus zero and genus one components of the tower (often called the ``two-level principle"). A non-comprehensive list of survey articles about this includes: \cite{Nakamura_survey}, \cite{POP2021107993}, and \cite{LS_braid}.  It seems then a good question to ask is: 
\[ \text{What is the ``Teichm\"uller tower''?}\] 

The geometric fundamental group of $\calM_{g,n}$ is also the profinite completion of the \emph{mapping class group} of a genus $g$ surface with $n$ marked points, $\Gamma_{g, n}$.\footnote{In brief, a point in the moduli space $\mathcal{M}_{g,n}$ is an isomorphism class of Riemann surfaces and a loop (up to homotopy) is an orientation preserving diffeomorphism of the basepoint surface, up to those homotopic to the identity. Full details can be found in \cite{MR1483111}.} Using this fact, Hatcher, Lochak and Schneps \cite{hls} proposed a model for the Teichm\"uller tower as the collection of all profinite mapping class groups $\widehat{\Gamma}_{g,n}$ together with all homomorphisms $\widehat{\Gamma}_{h, m}\hookrightarrow \widehat{\Gamma}_{g,n}$ induced by the inclusion of ``nice'' subsurfaces.  Here, $\Sigma'$ is a nice subsurface of $\Sigma$, if $\Sigma'$ is obtained by cutting along a set of disjoint simple closed curves on $\Sigma$. The natural embedding $\Sigma' \hookrightarrow \Sigma$ induces a map on the mapping class groups via the inclusion of Dehn twists (cf. Section~\ref{sec: mapping class}). 

Any good model for our ideal Teichm\"uller tower would, in particular, have the property that the Galois action on each of the $\widehat{\Gamma}_{g,n}$ necessarily commutes with all the maps in the tower. In practice, this is done by studying the actions of a family of more explicitly defined profinite groups on a proposed model of the tower. The most famous example is the Grothendieck-Teichm\"uller group, $\gt$, introduced by Drinfeld,  we have described (\cite[Section 4]{Drin}, \cite{ihara}). There are several other interesting groups related to $\gt$, but in these lectures we will focus on the group $\ns$, introduced by Nakamura and Schneps, which is known to act on the tower which includes all the higher genus mapping class groups $\widehat{\Gamma}_{g, n}$ (\cite[Theorem 1.3]{ns}).  The advantage in studying actions of $\gt$ and $\ns$ on the Teichm\"uller tower is that they have explicit presentations (Definition~\ref{def: gt} and \ref{def: ns}) and live ``in between'' the absolute Galois group and the mapping class groups in a very explicit way. In particular, Ihara \cite{ihara} showed that the image of the action of $\Gal(\mathbb{Q})$ on the geometric fundamental group $\pi_{1}^{et}(\calM_{0,4})$ lies in $\gt$. Nakamura and Schneps show that their group $\ns$ is a subgroup of $\gt$ which still contains $\Gal(\mathbb{Q})$ (\cite[Theorem 1.2]{ns}).  In other words, \[\Gal(\mathbb{Q})\hookrightarrow\ns\hookrightarrow \gt.\]

\medskip

The main goal of this lecture series is to explain an ``operadic'' interpretation of the Teichm\"uller tower in Hatcher, Lochak and Schneps \cite{hls}. In this approach, we consider the mapping class groups of surfaces of genus $g$ with $n$ boundary components, $\Gamma^{g}_{n}$, rather than marked points, $\Gamma_{g,n}$. With this small change, we can assemble the collection of spaces $B\Gamma^{g}_{n}$ into a \emph{modular operad} (Definition~\ref{def: modular}) with operations defined by gluing along the boundary components. In this case, the nice subsurface inclusions in the tower of Hatcher, Lochak and Schneps are equivalent to modular operad operations as in Figure~\ref{comp and contract}. Topologically, this is only a minor change, as there is a short exact sequence
\[\begin{tikzcd}
0\arrow[r] & \mathbb{Z}^{n}\arrow[r]&\Gamma^{g}_{n}\arrow[r]&\Gamma_{g,n}\arrow[r]& 0
\end{tikzcd}\] where the map $\Gamma^{g}_{n}\rightarrow\Gamma_{g, n}$ collapses boundary components to points. 

\begin{figure}[h!]
\[
\begin{tikzpicture}[x=0.75pt,y=0.75pt,yscale=-1,xscale=1]

\draw    (78.29,147.93) .. controls (87.14,152.99) and (92.11,153.65) .. (102.98,148.15) ;
\draw    (81.05,149.25) .. controls (88.06,145.07) and (92.66,144.19) .. (100.4,149.69) ;

\draw  [color={rgb, 255:red, 4; green, 146; blue, 194 }  ,draw opacity=1 ] (98.76,110.53) .. controls (100.12,108.86) and (107.3,112.45) .. (114.79,118.56) .. controls (122.29,124.66) and (127.26,130.96) .. (125.9,132.63) .. controls (124.54,134.3) and (117.36,130.71) .. (109.87,124.61) .. controls (102.37,118.51) and (97.4,112.2) .. (98.76,110.53) -- cycle ;
\draw   (53.82,125.11) .. controls (55.98,125.11) and (57.73,132.94) .. (57.74,142.61) .. controls (57.74,152.28) and (56,160.11) .. (53.85,160.11) .. controls (51.69,160.12) and (49.94,152.28) .. (49.93,142.61) .. controls (49.93,132.95) and (51.67,125.11) .. (53.82,125.11) -- cycle ;
\draw   (88.04,187.97) .. controls (87.71,185.84) and (95.17,182.89) .. (104.72,181.38) .. controls (114.27,179.87) and (122.28,180.37) .. (122.62,182.5) .. controls (122.95,184.63) and (115.49,187.58) .. (105.94,189.09) .. controls (96.39,190.6) and (88.38,190.1) .. (88.04,187.97) -- cycle ;
\draw    (53.82,125.11) .. controls (63.13,126.75) and (92.04,118.57) .. (98.76,110.53) ;
\draw    (53.85,160.11) .. controls (65.22,162.65) and (84.14,178.87) .. (88.04,187.97) ;
\draw    (122.62,182.5) .. controls (116.66,168.28) and (115.27,143.67) .. (125.9,132.63) ;
\draw   (106.37,62.46) .. controls (104.1,61.55) and (104.23,53.35) .. (106.66,44.14) .. controls (109.08,34.94) and (112.89,28.21) .. (115.16,29.12) .. controls (117.43,30.03) and (117.3,38.23) .. (114.88,47.43) .. controls (112.45,56.64) and (108.64,63.36) .. (106.37,62.46) -- cycle ;
\draw  [color={rgb, 255:red, 4; green, 146; blue, 194 }  ,draw opacity=1 ] (146.98,108.37) .. controls (145.67,110.39) and (138.57,107.06) .. (131.13,100.92) .. controls (123.7,94.78) and (118.73,88.16) .. (120.05,86.14) .. controls (121.36,84.11) and (128.45,87.45) .. (135.89,93.59) .. controls (143.33,99.72) and (148.3,106.34) .. (146.98,108.37) -- cycle ;
\draw   (160.78,12.5) .. controls (161.61,10.15) and (169.08,11.66) .. (177.46,15.87) .. controls (185.84,20.07) and (191.97,25.38) .. (191.14,27.73) .. controls (190.32,30.07) and (182.85,28.56) .. (174.46,24.36) .. controls (166.08,20.15) and (159.95,14.84) .. (160.78,12.5) -- cycle ;
\draw   (202.03,52.58) .. controls (204.39,52.8) and (206.3,60.92) .. (206.29,70.71) .. controls (206.28,80.51) and (204.35,88.27) .. (201.99,88.05) .. controls (199.63,87.83) and (197.72,79.71) .. (197.73,69.91) .. controls (197.74,60.12) and (199.67,52.36) .. (202.03,52.58) -- cycle ;
\draw    (201.99,88.05) .. controls (187.79,81.26) and (153.53,92.23) .. (146.98,108.37) ;
\draw    (160.78,12.49) .. controls (155.82,25.46) and (126.82,33.7) .. (115.16,29.12) ;
\draw    (202.03,52.58) .. controls (192.35,49.39) and (184.35,38.19) .. (191.14,27.73) ;
\draw [color={rgb, 255:red, 4; green, 146; blue, 194 }  ,draw opacity=1 ] [dash pattern={on 0.84pt off 2.51pt}]  (120.05,86.14) .. controls (118.29,99.29) and (108.69,107.3) .. (98.76,110.53) ;
\draw [color={rgb, 255:red, 4; green, 146; blue, 194 }  ,draw opacity=1 ] [dash pattern={on 0.84pt off 2.51pt}]  (146.98,108.37) .. controls (135,114.33) and (130,121.33) .. (125.9,132.63) ;
\draw    (310.29,135.93) .. controls (319.14,140.99) and (324.11,141.65) .. (334.98,136.15) ;
\draw    (313.05,137.25) .. controls (320.06,133.07) and (324.66,132.19) .. (332.4,137.69) ;

\draw   (285.82,113.11) .. controls (287.98,113.11) and (289.73,120.94) .. (289.74,130.61) .. controls (289.74,140.28) and (288,148.11) .. (285.85,148.11) .. controls (283.69,148.12) and (281.94,140.28) .. (281.93,130.61) .. controls (281.93,120.95) and (283.67,113.11) .. (285.82,113.11) -- cycle ;
\draw   (320.04,175.97) .. controls (319.71,173.84) and (327.17,170.89) .. (336.72,169.38) .. controls (346.27,167.87) and (354.28,168.37) .. (354.62,170.5) .. controls (354.95,172.63) and (347.49,175.58) .. (337.94,177.09) .. controls (328.39,178.6) and (320.38,178.1) .. (320.04,175.97) -- cycle ;
\draw    (285.82,113.11) .. controls (295.13,114.75) and (324.04,106.57) .. (330.76,98.53) ;
\draw    (285.85,148.11) .. controls (297.22,150.65) and (316.14,166.87) .. (320.04,175.97) ;
\draw    (354.62,170.5) .. controls (348.66,156.28) and (347.27,131.67) .. (357.9,120.63) ;
\draw   (317.37,74.46) .. controls (315.1,73.55) and (315.23,65.35) .. (317.66,56.14) .. controls (320.08,46.94) and (323.89,40.21) .. (326.16,41.12) .. controls (328.43,42.03) and (328.3,50.23) .. (325.88,59.43) .. controls (323.45,68.64) and (319.64,75.36) .. (317.37,74.46) -- cycle ;
\draw  [color={rgb, 255:red, 4; green, 146; blue, 194 }  ,draw opacity=1 ][dash pattern={on 0.84pt off 2.51pt}] (357.98,120.37) .. controls (356.67,122.39) and (349.57,119.06) .. (342.13,112.92) .. controls (334.7,106.78) and (329.73,100.16) .. (331.05,98.14) .. controls (332.36,96.11) and (339.45,99.45) .. (346.89,105.59) .. controls (354.33,111.72) and (359.3,118.34) .. (357.98,120.37) -- cycle ;
\draw   (371.78,24.5) .. controls (372.61,22.15) and (380.08,23.66) .. (388.46,27.87) .. controls (396.84,32.07) and (402.97,37.38) .. (402.14,39.73) .. controls (401.32,42.07) and (393.85,40.56) .. (385.46,36.36) .. controls (377.08,32.15) and (370.95,26.84) .. (371.78,24.5) -- cycle ;
\draw   (413.03,64.58) .. controls (415.39,64.8) and (417.3,72.92) .. (417.29,82.71) .. controls (417.28,92.51) and (415.35,100.27) .. (412.99,100.05) .. controls (410.63,99.83) and (408.72,91.71) .. (408.73,81.91) .. controls (408.74,72.12) and (410.67,64.36) .. (413.03,64.58) -- cycle ;
\draw    (317.37,74.46) .. controls (332.33,75) and (336,92) .. (331.05,98.14) ;
\draw    (412.99,100.05) .. controls (398.79,93.26) and (364.53,104.23) .. (357.98,120.37) ;
\draw    (371.78,24.49) .. controls (366.82,37.46) and (337.82,45.7) .. (326.16,41.12) ;
\draw    (413.03,64.58) .. controls (403.35,61.39) and (395.35,50.19) .. (402.14,39.73) ;

\draw    (106.37,62.46) .. controls (121.33,63) and (125,80) .. (120.05,86.14) ;
\draw  [color={rgb, 255:red, 4; green, 146; blue, 194 }  ,draw opacity=1 ] (155.14,270.63) .. controls (156.08,268.38) and (164.28,268.63) .. (173.45,271.19) .. controls (182.61,273.76) and (189.28,277.66) .. (188.34,279.92) .. controls (187.39,282.18) and (179.2,281.93) .. (170.03,279.36) .. controls (160.86,276.8) and (154.2,272.89) .. (155.14,270.63) -- cycle ;
\draw  [color={rgb, 255:red, 4; green, 146; blue, 194 }  ,draw opacity=1 ] (134.08,213.28) .. controls (136.1,211.96) and (141.24,217.88) .. (145.56,226.5) .. controls (149.87,235.12) and (151.74,243.18) .. (149.72,244.5) .. controls (147.69,245.82) and (142.56,239.9) .. (138.24,231.28) .. controls (133.92,222.66) and (132.06,214.6) .. (134.08,213.28) -- cycle ;
\draw   (128.57,325.96) .. controls (128.36,328.44) and (120.75,328.83) .. (111.58,326.84) .. controls (102.41,324.85) and (95.16,321.24) .. (95.38,318.76) .. controls (95.6,316.28) and (103.2,315.89) .. (112.37,317.88) .. controls (121.54,319.87) and (128.79,323.49) .. (128.57,325.96) -- cycle ;
\draw   (78.66,297.4) .. controls (76.31,297.77) and (72.45,290.38) .. (70.02,280.89) .. controls (67.6,271.4) and (67.53,263.41) .. (69.87,263.03) .. controls (72.21,262.66) and (76.08,270.05) .. (78.5,279.54) .. controls (80.93,289.03) and (81,297.02) .. (78.66,297.4) -- cycle ;
\draw    (69.87,263.03) .. controls (103.67,253.33) and (102.67,225.33) .. (134.08,213.28) ;
\draw    (128.57,325.96) .. controls (140.67,299.33) and (175.67,310.33) .. (188.34,279.92) ;
\draw    (78.65,297.4) .. controls (88.82,298.08) and (97,306) .. (95.38,318.76) ;
\draw    (155.14,270.63) .. controls (136,282.67) and (122.67,249.33) .. (149.72,244.5) ;
\draw    (102.29,283.93) .. controls (111.14,288.99) and (116.11,289.65) .. (126.98,284.15) ;
\draw    (105.05,285.25) .. controls (112.06,281.07) and (116.66,280.19) .. (124.4,285.69) ;

\draw [color={rgb, 255:red, 4; green, 146; blue, 194 }  ,draw opacity=1 ] [dash pattern={on 0.84pt off 2.51pt}]  (188.34,279.92) .. controls (201,234.33) and (183,207) .. (134.08,213.28) ;
\draw [color={rgb, 255:red, 4; green, 146; blue, 194 }  ,draw opacity=1 ] [dash pattern={on 0.84pt off 2.51pt}]  (155.14,270.63) .. controls (168,264) and (163,247) .. (149.72,244.5) ;

\draw   (349.57,326.96) .. controls (349.36,329.44) and (341.75,329.83) .. (332.58,327.84) .. controls (323.41,325.85) and (316.16,322.24) .. (316.38,319.76) .. controls (316.6,317.28) and (324.2,316.89) .. (333.37,318.88) .. controls (342.54,320.87) and (349.79,324.49) .. (349.57,326.96) -- cycle ;
\draw   (299.66,298.4) .. controls (297.31,298.77) and (293.45,291.38) .. (291.02,281.89) .. controls (288.6,272.4) and (288.53,264.41) .. (290.87,264.03) .. controls (293.21,263.66) and (297.08,271.05) .. (299.5,280.54) .. controls (301.93,290.03) and (302,298.02) .. (299.66,298.4) -- cycle ;
\draw    (290.87,264.03) .. controls (324.67,254.33) and (323.82,222.22) .. (355.08,214.28) ;
\draw    (349.57,326.96) .. controls (353.67,292) and (379,307) .. (400.33,275) .. controls (421.67,243) and (403.67,200) .. (355.08,214.28) ;
\draw    (299.65,298.4) .. controls (309.82,299.08) and (318,307) .. (316.38,319.76) ;
\draw    (323.29,284.93) .. controls (332.14,289.99) and (337.11,290.65) .. (347.98,285.15) ;
\draw    (326.05,286.25) .. controls (333.06,282.07) and (337.66,281.19) .. (345.4,286.69) ;

\draw  [color={rgb, 255:red, 4; green, 146; blue, 194 }  ,draw opacity=1 ][dash pattern={on 0.84pt off 2.51pt}] (372.63,245.02) .. controls (371.36,243.35) and (376.1,236.39) .. (383.21,229.48) .. controls (390.32,222.56) and (397.11,218.31) .. (398.37,219.98) .. controls (399.64,221.65) and (394.9,228.61) .. (387.79,235.52) .. controls (380.68,242.44) and (373.89,246.69) .. (372.63,245.02) -- cycle ;
\draw    (358.29,247.93) .. controls (367.14,252.99) and (372.11,253.65) .. (382.98,248.15) ;
\draw    (361.05,249.25) .. controls (368.06,245.07) and (372.66,244.19) .. (380.4,249.69) ;

\draw    (201,110.5) -- (246,110.5) ;
\draw [shift={(248,110.5)}, rotate = 180] [color={rgb, 255:red, 0; green, 0; blue, 0 }  ][line width=0.75]    (10.93,-3.29) .. controls (6.95,-1.4) and (3.31,-0.3) .. (0,0) .. controls (3.31,0.3) and (6.95,1.4) .. (10.93,3.29)   ;
\draw    (219,270.5) -- (264,270.5) ;
\draw [shift={(266,270.5)}, rotate = 180] [color={rgb, 255:red, 0; green, 0; blue, 0 }  ][line width=0.75]    (10.93,-3.29) .. controls (6.95,-1.4) and (3.31,-0.3) .. (0,0) .. controls (3.31,0.3) and (6.95,1.4) .. (10.93,3.29)   ;

\draw (189,117) node [anchor=north west][inner sep=0.75pt]   [align=left] {{\footnotesize Composition}};
\draw (211,277) node [anchor=north west][inner sep=0.75pt]   [align=left] {{\footnotesize Contraction}};

\end{tikzpicture}\]
    \caption{}
    \label{comp and contract}
\end{figure}

To argue that the modular operad $B\Gamma=\{B\Gamma^{g}_{n}\}$ models the Teichm\"uller tower, we start by replacing the mapping class groups $\Gamma^{g}_{n}$ by homotopy equivalent groupoids $\bS(g, n)$. The collection of spaces $B\bS(g,n)$ assembles into a modular operad, which we call the \emph{modular operad of seamed surfaces}.  The genus zero part of this modular operad contains an operad, $B\calS$, we call the \emph{genus zero surface operad} (\cite[Definition 6.5]{bhr1}). Boavida, Horel and the second author show that there is a faithful action of $\Gal(\mathbb{Q})$ on the profinite completion of the genus zero surface operad, following from the fact that there is an action of $\gt$ on $\widehat{B\calS}$ (\cite[Proposition 10.6]{bhr1}).  

In our third lecture, we demonstrate that there is a relatively straightforward action of $\gt$ on each groupoid $\bS(0, n)$, and this action is faithful when $n$ is at least $3$. The difficulty, which is beyond the scope of these lectures, is to show that this action of $\gt$, and thus the action of $\Gal(\mathbb{Q})$, is compatible with the operad structure. This compatibility is established in \cite{bhr1} by showing:
\begin{thm*}
The profinite Grothendieck-Teichm\"uller group $\gt$ is isomorphic to the group of homotopy automorphisms of the profinite completion of the genus zero surface operad.
\end{thm*} 

In \cite{BR22}, Bonatto and the second author show that the group $\ns$ acts on the profinite completion of the modular operad of seamed surfaces in such a way that, when restricting to genus zero, we recover the action of $\gt$. We can show quite directly that there is an action of $\ns$ on each groupoid $\bS(g,n)$, but the difficulty is showing that these $\ns$-actions are compatible with the modular operad structure maps. This is overcome by showing that $\ns$ is isomorphic to the group of homotopy automorphisms of the profinite completion of modular operad of seamed surfaces. However, unlike in the genus zero case, we make use an operadic version of the ``two-level principle'':
\begin{thm*}
The group $\ns$ acts on the group of homotopy automorphisms of the genus one truncation of the profinite completion of the modular operad of seamed surfaces.
\end{thm*} 

In this context, a \emph{genus one truncation} of $B\bS$ is a modular operad $B\bS_{\leq 1}$ in which $B\bS(g,n)=\emptyset$ if $g\geq2$. This translates into an operadic version of the two-level principle because one can show that there is a homotopy equivalence \[\mathbb{R}\End(\widehat{B\bS})\simeq \mathbb{R}\End(\widehat{B\bS}_{\leq 1}),\] where $\mathbb{R}\End(-)$ denotes the derived endomorphism space of modular dendroidal spaces (Definition~\ref{def: modular spaces}). 

\medskip
What our brief discussion on truncation illustrates is one of the many homotopical difficulties which arise due to profinite completion.  In particular, the profinite completion of a (modular) operad in spaces is no longer a 
(modular) operad.  In very special cases, however, the profinite completion of a (modular) operad can be considered as a (modular) operad whose operations hold ``up to homotopy''.  In order to describe these ``up to homotopy'', or $\infty$, modular operads the second author, together with Hackney and Yau, developed a Segal model for modular $\infty$-operads which provides a good setting for working with the profinite completion of modular operads. 

In the first two lectures of this series we introduce this homotopical background. In the first lecture describes a category $\bU$, whose objects are undirected, connected graphs with loose ends (Definition~\ref{def: graphs}). Morphisms are given by ‘blowing up’ vertices of the source into ``subgraphs'' of the target in a way that reflects iterated operations in a modular operad (Definition~\ref{def: graphical maps}). This graphical category models (discrete, coloured) modular operads in a very explicit way.
\begin{thm*}(Theorem~\ref{thm: nerve})
There exists an equivalence of categories \[\begin{tikzcd} \ModOp\arrow[r, shift left =.1 cm, "N"] & \arrow[l,  shift left =.1 cm] (\mathbf{Set}^{\bU^{op}})_{Segal}. \end{tikzcd}\]
\end{thm*}
The category on the right-hand side is a category of modular dendroidal sets (Definition~\ref{def:modular_dendroidal_sets}) satisfying a strict Segal condition (Definition~\ref{def: Segal map}). 

\medskip

The second lecture builds on the first, weakening the Segal condition of Theorem~\ref{thm: nerve} to provide a model for modular $\infty$-operads (Definition~\ref{def: weak Segal 2}). To do this we consider space-valued $\bU$-presheaves $X:\bU^{op}\rightarrow \mathbf{sSet}$ and say that (Definition~\ref{def: weak Segal 1}):
\begin{definition*}
A modular dendroidal  space $X\in \textbf{sSet}^{\bU^{op}}$ is \emph{Segal} if for all $G\in\bU$, the Segal map \[\begin{tikzcd} \map(\mathbf{U}[G], X)  \arrow[r] & \map(\mathbf{Sc}[G], X)\end{tikzcd}\] is a weak equivalence. 
\end{definition*}
Informally, what this definition says is that a Segal modular operad is a space-valued $\bU$-presheaf $X$ in which the value of $X$ at a graph $G$ is determined, up to homotopy, by the value of $X$ at each of the vertices of $G$. Finally, we describe profinite completion of modular operads and demonstrate why this operation results in a Segal modular operad.

\subsection{Structure and intentions of these notes}
These notes represent the transcript of three, $90$ minute lectures and, as such, fall unfortunately short of being a complete survey and introduction to either Grothendieck-Teichm\"uller theory or (modular) $\infty$-operads. For a more comprehensive overview of Grothendieck-Teichm\"uller theory there is an excellent collection \cite{MR1483106} edited by Lochak and Schneps and a rather recent survey article by Pop \cite{POP2021107993}.  The theory of modular $\infty$-operads is somewhat new, initiated in \cite{hry1, hry2}, and still in development, but there are several more comprehensive resources for an introduction to $\infty$-operads. We recommend the lecture series of Moerdijk \cite{Moerdijk_notes} as well as the book \cite{HM}.  

Throughout this series we take for granted that the reader will be familiar with the theory of operads. We give a definition of (coloured) modular operads in the first lecture but assume throughout that the reader, while maybe not familiar with all the details, is aware of the fact that there are adjunctions 
\[\begin{tikzcd} \textbf{Operad}  \arrow[r, shift left = .1 cm] & \textbf{Cyc} \arrow[l, shift left=.1cm]  \arrow[r, shift left=.1cm] & \arrow[l, shift left =.1 cm] \ModOp.\end{tikzcd}\] where the right adjoint functors are ``forgetful'' functors. In addition, the open problems we give in Section~\ref{sec:further directions 1} assume the reader is familiar with the other ``operad-like'' objects in the literature such as properads, dioperads and wheeled properads. These results are scattered throughout the literature, but a good first introduction is the survey article \cite{markl_survey}.  

We avoid delving too deep into the homotopical properties of Segal modular operads and suppress many arguments which require the use of Quillen model categories and/or $\infty$-categories. Our goal for these lectures is to give the reader some basic understanding of modular $\infty$-operads and their applications without having to understand the full theory. The only real homotopical prerequisite is that the reader have some understanding of \emph{derived mapping spaces} or \emph{homotopy function complexes} which we denote by $\map(X,Y)$, throughout. A good standard reference is \cite[Chapter 9]{hirschhorn}.  
\medskip

In the original lecture series we gave several exercises and a series of open problems to which we are reasonably confident there is an answer. Most of these open problems are in Section~\ref{sec:further directions 1} and Section~\ref{sec:open problems 2}, though some are scattered throughout the notes.

\subsection{Acknowledgments} 
The work presented in these lectures covers joint work of the second author with Luciana Basualdo Bonatto, Pedro Boavida de Brito, Philip Hackney, Geoffroy Horel, and Donald Yau. We would like to thank all of them for their work, comments, and suggestions and take full responsibility for any typos, errors or bad writing contained in these notes. 

Over the years this work has also benefited immensely from feedback and comments from Ezra Getzler, Joachim Kock, Ieke Moerdijk, Sophie Raynor, Marco Robalo, Leila Schneps, Michelle Strumila and many others.  We would like to specifically acknowledge Luciana Basualdo Bonatto, Philip Hackney and Sophie Raynor who provided feedback on these lectures as they were being prepared and the graduate students at the University of Melbourne who were subjected to the first version(s) of these lectures. A special thank-you goes to Santiago Nahuel Martinez who made many of the beautiful figures in these notes. 

Lastly, we would like to thank the organizers and participants of the workshop on ``Higher Structures and Operadic Calculus''. We are particularly grateful for the efforts made to make this event possible during a (hopefully!) once-in-a-lifetime pandemic and for making it possible for the authors to participate remotely from the USA and Australia, respectively.

\section{Lecture 1: Graphs and Modular Operads} 
Modular operads are a generalization of operads which allow one to encode algebraic structures that come equipped with a ``bilinear form'' or ``contraction operation''. They were introduced by Getzler and Kapronov \cite{gk_modular} where they gave the canonical example of a modular operad built from the (compactified) moduli spaces of genus $g$ curves, $\bar{\mathcal{M}}_{g,n}$. 

This first lecture introduces a very general definition of coloured modular operads in which we allow for an involutive set of colours (Definition~\ref{def: modular}). This differs slightly from the usual definition of a coloured modular operad in the literature (\cite{hvz10}, \cite{DM16}, \cite{kw10}, etc) where the colour set has trivial involution. Interesting examples of operads with involutive colour sets include \cite{DCH_Dwyer_Kan} where the authors study coloured \emph{cyclic operads} with involutive set of colours and the examples of modular operads coloured by involutive groupoids in \cite{pet}. In the case of one-coloured, or monochrome, modular operads the involution on colour sets is trivial and thus all the reader's favorite modular operads are still examples of our definition. 
As we will discuss in the open problems section (Section~\ref{sec:further directions 1}), an advantage of considering involutive colour sets is that coloured cyclic operads, coloured operads, and coloured dioperads can all be considered as special types of coloured modular operads. Similarly, we can consider wheeled properads as a special case of modular operads with involutive colour sets.

The main theorem in this lecture is a so-called \emph{nerve theorem} (Theorem~\ref{thm: nerve}). This is an extension of the classical theorem which says that the inclusion of the simplex category $\Delta$ into the category of small categories induces a fully faithful functor from $\mathbf{Cat}$ into the category of simplicial sets $\mathbf{sSet}^{\Delta^{op}}$ whose essential image consists of the Segal objects, i.e. those $X\in\mathbf{sSet}^{\Delta^{op}}$ with
\[X_n\cong X_1\times_{X_0}X_1\times_{X_0}\ldots \times_{X_0}X_1\] for all $n\geq 2.$

\subsection{Cyclic Operads}
Modular operads are cyclic operads with contraction operations. As such, we first introduce the notion of a coloured cyclic operad: 

\begin{definition} Let $\mathfrak{C}$ be a non-empty set.  
A $\mathfrak{C}$-coloured \emph{cyclic operad} $\bP$ consists of:   
\begin{enumerate} 
\item An involutive set $\mathfrak{C}$ which we call the set of colours and denote by $\col(\bP)$. We will write $c\mapsto c^{\dagger}$ for the action of the involution on an element $c\in\mathfrak{C}$.
\item A collection of sets $\bP=\{\bP(c_1,\ldots, c_n)\}$, in which, for each $c_1,\ldots, c_n\in\mathfrak{C}$, the set $\bP(c_1,\ldots, c_n)$ in equipped with a right $\Sigma_n$-action \[\begin{tikzcd}
\bP(c_1,\ldots, c_n) \arrow[r, "\sigma^*"] &\bP(c_{\sigma(1)},\ldots c_{\sigma(n)}) \end{tikzcd}.\]
\item A set of distinguished unit elements $id_c\in\bP(c^{\dagger},c)$, one for each $c\in\mathfrak{C}$. 
\item A family of associative, unital and equivariant composition operations \[\begin{tikzcd}\bP(c_1,\ldots, c_n) \times \bP(d_1,\ldots, d_m) \arrow[rr, "\circ_{ij}"]&& \bP(c_1,\ldots, \hat{c}_i,\ldots, d_1,\ldots \hat{d}_j,\ldots, d_m), \end{tikzcd}\] whenever $c_i=d^{\dagger}_j$,  $(i, j) \in [1, n] \times [1, m]$.\\
\end{enumerate} 
\end{definition} 

A morphism of coloured cyclic operads $f: \bP\rightarrow \bQ$ consists of an involutive function $f: \col(\bP)\rightarrow \col(\bQ)$ together with a family of $\Sigma_n$-equivariant maps \[\begin{tikzcd} \bP(c_1,\ldots, c_n) \arrow[r, "f_{\underline{c}}"] & \bQ(f(c_1),\ldots,f(c_n))\end{tikzcd}\] for every list $c_1,\ldots, c_n\in\col(\bP)$ which commute with composition and identities. We will denote the category of all coloured cyclic operads by $\mathsf{Cyc}$. 

\begin{remark}
\begin{enumerate} 
\item In this first lecture we are only discussing discrete cyclic and modular operads, but one may extend the above to define coloured cyclic operads enriched in any closed, symmetric monoidal category $\bE=(\bE,\otimes,1)$. In this case, we denote the category of cyclic operads by $\mathsf{Cyc}(\bE)$. 

\item In order to avoid introducing too much notation, we have opted not to include the full list of axioms for coloured cyclic operads here. A good reference for these axioms is Definition 2.6 of \cite{DCH_Dwyer_Kan}. 
\end{enumerate} 
\end{remark} 

\begin{figure}[h!]
\[
\begin{tikzpicture}[x=0.75pt,y=0.75pt,yscale=-.75,xscale=.75]

\draw [color={rgb, 255:red, 208; green, 2; blue, 27 }  ,draw opacity=1 ]   (76,92) -- (110,133) ;
\draw  [fill={rgb, 255:red, 255; green, 255; blue, 255 }  ,fill opacity=1 ] (120.97,130.07) .. controls (129.25,130) and (136.02,136.66) .. (136.09,144.94) .. controls (136.16,153.22) and (129.5,160) .. (121.22,160.07) .. controls (112.93,160.14) and (106.16,153.48) .. (106.09,145.19) .. controls (106.02,136.91) and (112.68,130.14) .. (120.97,130.07) -- cycle ;
\draw [color={rgb, 255:red, 65; green, 117; blue, 5 }  ,draw opacity=1 ]   (111,158) -- (73,203) ;
\draw [color={rgb, 255:red, 144; green, 19; blue, 254 }  ,draw opacity=1 ]   (130,156) -- (170,204) ;
\draw    (133,136) -- (170,91) ;
\draw  [fill={rgb, 255:red, 255; green, 255; blue, 255 }  ,fill opacity=1 ] (416.97,88.07) .. controls (425.25,88) and (432.02,94.66) .. (432.09,102.94) .. controls (432.16,111.22) and (425.5,118) .. (417.22,118.07) .. controls (408.93,118.14) and (402.16,111.48) .. (402.09,103.19) .. controls (402.02,94.91) and (408.68,88.14) .. (416.97,88.07) -- cycle ;
\draw [color={rgb, 255:red, 208; green, 2; blue, 27 }  ,draw opacity=1 ]   (354,49) -- (403,96) ;
\draw    (430,94) -- (473,50) ;
\draw    (423,118) -- (465,164) ;
\draw  [fill={rgb, 255:red, 255; green, 255; blue, 255 }  ,fill opacity=1 ] (475.97,161.07) .. controls (484.25,161) and (491.02,167.66) .. (491.09,175.94) .. controls (491.16,184.22) and (484.5,191) .. (476.22,191.07) .. controls (467.93,191.14) and (461.16,184.48) .. (461.09,176.19) .. controls (461.02,167.91) and (467.68,161.14) .. (475.97,161.07) -- cycle ;
\draw [color={rgb, 255:red, 65; green, 117; blue, 5 }  ,draw opacity=1 ]   (476.22,191.07) -- (476,249) ;
\draw [color={rgb, 255:red, 144; green, 19; blue, 254 }  ,draw opacity=1 ]   (492,171) -- (565,170) ;

\draw (113,135) node [anchor=north west][inner sep=0.75pt]    {$p_v$};
\draw (409,94.4) node [anchor=north west][inner sep=0.75pt]    {$q_{v}$};
\draw (466,166) node [anchor=north west][inner sep=0.75pt]    {$q_{w}$};

\end{tikzpicture}\]
\caption{Operations of a cyclic operad pictured as decorated graphs.}\label{cyclic operations}
\end{figure}
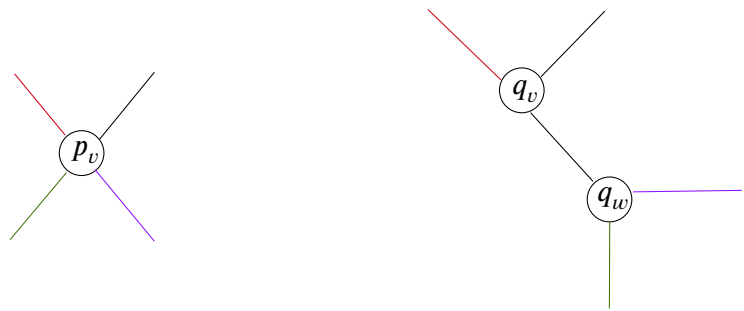

Informally, operations of a cyclic operad can be pictured as \emph{simply connected} graphs with $n$-free edges in the \emph{boundary} whose vertices are \emph{decorated} by elements of the underlying collection $\bP=\{\bP(c_1,\ldots,c_n)\}$.\footnote{A precise definition of decoration can be found in \cite[Definition 2.7]{hry2}.} Two examples of such decorated graphs are depicted in Figure~\ref{cyclic operations}.  The operation on the left is a decoration of a \emph{star}--a simply connected graph with a single vertex $v$. The element $p_v\in\bP(\textcolor{red}{c_1},c_2, \textcolor{violet}{c_3},\textcolor{teal}{c_4})$.  Similarly, the operation on the right is depicting an operation obtained by the composition of two stars. In a moment we will introduce a specific way to label the loose ends of our graphs and make precise the graphical interpretation of when two ends can be composed, but for this informal discussion we just use colours. This informal depiction will be made precise via Theorem~\ref{thm: nerve}.

If we restrict to one-coloured, or monochrome, cyclic operads the involution is trivial, and as such the familiar examples of one-coloured cyclic operads are still objects in $\mathsf{Cyc}(\bE)$. 

\begin{example}
An $\ast$-\emph{autonomous category} $\mathcal{C}$ is a closed, symmetric monoidal category in which every object $x$ has a dual $x^{\dagger}$, satisfying the property that $x\cong (x^{\dagger})^{\dagger}$ (\cite{Barr79}). When the double dual relation is strict, $x =(x^{\dagger})^{\dagger}$, $*$-autonomous categories are examples of coloured cyclic operads (Example 2.12 of \cite{DCH_Dwyer_Kan}). To our knowledge, it is an open problem to show that all $*$-autonomous categories are coloured cyclic operads.
\end{example} 

\subsection{Modular Operads}

\begin{definition}\label{def: modular}
A $\mathfrak{C}$-coloured \emph{modular operad} is a $\mathfrak{C}$-coloured cyclic operad which also has a family of contraction operations \[\begin{tikzcd}\bP(c_1,\ldots, c_n) \arrow[r, "\xi_{ij}"] & \bP(c_1,\ldots, \hat{c}_i, \ldots,  \hat{c}_j,\ldots  c_n),\end{tikzcd}\] whenever $c_i=c^{\dagger}_j$, $0\leq i < j \leq n$.  
\end{definition}

\begin{remark}
Modular operads, as we have defined them in this first lecture, are called ``compact symmetric multicategories'' in \cite{jk11}, \cite{Sophie_phdthesis} and \cite{Sophie_21}.
\end{remark}

We require the composition and contraction operations to satisfy a series of commutativity, associativity, unitality and equivariance constraints. We do not list the full list of axioms here, as doing so requires one to be extremely careful with re-indexing and requires the introduction of a significant amount of notion. We do not think this comes at much of a cost because, as we will show, the applications in these notes will often use the identification in Theorem~\ref{thm: nerve}.  For the curious reader, many of the axioms needed are illustrated beautifully in Definition 1.24 of \cite{Sophie_21} and we also recommend \cite[Section 13.4]{yau_bv}, where a complete set of axioms is given for coloured modular operads with non-involutive colour sets. To give one example, Figure~\ref{cyclic and contraction operations} illustrates compatibility of composition and contraction operations.  The top of the diagram depicts first composing the two stars along the red edges and then contracting the two purple edges. The bottom of the diagram depicts first composing along the purple edges and then contracting the red edges. 

\begin{figure}[h!]

\[\begin{tikzpicture}[x=0.75pt,y=0.75pt,yscale=-1,xscale=1]

\draw [fill={rgb, 255:red, 255; green, 255; blue, 255 }  ,fill opacity=1 ] (587.97,97.07) .. controls (596.25,97) and (603.02,103.66) .. (603.09,111.94) .. controls (603.16,120.22) and (596.5,127) .. (588.22,127.07) .. controls (579.93,127.14) and (573.16,120.48) .. (573.09,112.19) .. controls (573.02,103.91) and (579.68,97.14) .. (587.97,97.07) -- cycle ;
\draw [color={rgb, 255:red, 0; green, 0; blue, 0 }  ,draw opacity=1 ]   (536,62) -- (574,105) ;
\draw    (601,103) -- (644,59) ;
\draw  [fill={rgb, 255:red, 255; green, 255; blue, 255 }  ,fill opacity=1 ] (589.97,180.07) .. controls (598.25,180) and (605.02,186.66) .. (605.09,194.94) .. controls (605.16,203.22) and (598.5,210) .. (590.22,210.07) .. controls (581.93,210.14) and (575.16,203.48) .. (575.09,195.19) .. controls (575.02,186.91) and (581.68,180.14) .. (589.97,180.07) -- cycle ;
\draw [color={rgb, 255:red, 208; green, 2; blue, 27 }  ,draw opacity=1 ]   (575,190) .. controls (567,181) and (566,135) .. (576,123) ;
\draw [color={rgb, 255:red, 144; green, 19; blue, 254 }  ,draw opacity=1 ]   (599,123) .. controls (607,136) and (609,176) .. (603,187) ;
\draw    (590.22,210.07) -- (591,255) ;
\draw [color={rgb, 255:red, 0; green, 0; blue, 0 }  ,draw opacity=1 ]   (37.1,34.81) -- (60.94,68.19) ;
\draw  [fill={rgb, 255:red, 255; green, 255; blue, 255 }  ,fill opacity=1 ] (68.63,65.81) .. controls (74.43,65.75) and (79.18,71.17) .. (79.23,77.92) .. controls (79.28,84.66) and (74.61,90.17) .. (68.8,90.23) .. controls (62.99,90.29) and (58.25,84.87) .. (58.2,78.12) .. controls (58.15,71.38) and (62.82,65.86) .. (68.63,65.81) -- cycle ;
\draw [color={rgb, 255:red, 208; green, 2; blue, 27 }  ,draw opacity=1 ]   (61.64,88.55) -- (35,125.19) ;
\draw [color={rgb, 255:red, 144; green, 19; blue, 254 }  ,draw opacity=1 ]   (74.96,86.92) -- (103,126) ;
\draw    (77.06,70.64) -- (103,34) ;
\draw [color={rgb, 255:red, 208; green, 2; blue, 27 }  ,draw opacity=1 ]   (37.1,174.81) -- (60.94,208.19) ;
\draw  [fill={rgb, 255:red, 255; green, 255; blue, 255 }  ,fill opacity=1 ] (68.63,205.81) .. controls (74.43,205.75) and (79.18,211.17) .. (79.23,217.92) .. controls (79.28,224.66) and (74.61,230.17) .. (68.8,230.23) .. controls (62.99,230.29) and (58.25,224.87) .. (58.2,218.12) .. controls (58.15,211.38) and (62.82,205.86) .. (68.63,205.81) -- cycle ;
\draw [color={rgb, 255:red, 0; green, 0; blue, 0 }  ,draw opacity=1 ]   (68.8,230.23) -- (69.16,271.68) ;
\draw [color={rgb, 255:red, 144; green, 19; blue, 254 }  ,draw opacity=1 ]   (77.06,209.64) -- (103,173) ;
\draw  [fill={rgb, 255:red, 255; green, 255; blue, 255 }  ,fill opacity=1 ] (301.34,27.28) .. controls (308.09,27.23) and (313.61,31.48) .. (313.67,36.76) .. controls (313.72,42.04) and (308.3,46.36) .. (301.55,46.41) .. controls (294.8,46.45) and (289.28,42.21) .. (289.22,36.92) .. controls (289.17,31.64) and (294.59,27.32) .. (301.34,27.28) -- cycle ;
\draw [color={rgb, 255:red, 0; green, 0; blue, 0 }  ,draw opacity=1 ]   (259,4.91) -- (289.96,32.34) ;
\draw    (311.96,31.06) -- (347,3) ;
\draw  [fill={rgb, 255:red, 255; green, 255; blue, 255 }  ,fill opacity=1 ] (302.97,80.21) .. controls (309.72,80.17) and (315.24,84.41) .. (315.3,89.7) .. controls (315.35,94.98) and (309.93,99.3) .. (303.18,99.34) .. controls (296.43,99.39) and (290.91,95.14) .. (290.85,89.86) .. controls (290.8,84.58) and (296.22,80.26) .. (302.97,80.21) -- cycle ;
\draw [color={rgb, 255:red, 208; green, 2; blue, 27 }  ,draw opacity=1 ]   (290.78,86.55) .. controls (284.26,80.81) and (283.44,51.47) .. (291.59,43.82) ;
\draw    (303.18,99.34) -- (303.81,128) ;
\draw [color={rgb, 255:red, 144; green, 19; blue, 254 }  ,draw opacity=1 ]   (309,44.91) -- (326,59) ;
\draw [color={rgb, 255:red, 144; green, 19; blue, 254 }  ,draw opacity=1 ]   (313.96,84.06) -- (330,69) ;
\draw  [fill={rgb, 255:red, 255; green, 255; blue, 255 }  ,fill opacity=1 ] (304.38,206.53) .. controls (310.98,206.49) and (316.37,210.43) .. (316.42,215.33) .. controls (316.48,220.23) and (311.18,224.24) .. (304.58,224.28) .. controls (297.98,224.32) and (292.59,220.38) .. (292.54,215.48) .. controls (292.48,210.58) and (297.78,206.57) .. (304.38,206.53) -- cycle ;
\draw [color={rgb, 255:red, 0; green, 0; blue, 0 }  ,draw opacity=1 ]   (263,185.78) -- (293.26,211.22) ;
\draw    (314.76,210.04) -- (349,184) ;
\draw [fill={rgb, 255:red, 255; green, 255; blue, 255 }  ,fill opacity=1 ] (305.97,255.65) .. controls (312.57,255.61) and (317.96,259.55) .. (318.02,264.45) .. controls (318.07,269.36) and (312.77,273.36) .. (306.17,273.41) .. controls (299.58,273.45) and (294.18,269.51) .. (294.13,264.6) .. controls (294.07,259.7) and (299.38,255.69) .. (305.97,255.65) -- cycle ;
\draw [color={rgb, 255:red, 144; green, 19; blue, 254 }  ,draw opacity=1 ]   (313.17,221.88) .. controls (319.54,229.57) and (321.13,253.24) .. (316.35,259.76) ;
\draw    (306.17,273.41) -- (306.8,300) ;
\draw [color={rgb, 255:red, 208; green, 2; blue, 27 }  ,draw opacity=1 ]   (294.64,220.55) -- (277,240) ;
\draw [color={rgb, 255:red, 208; green, 2; blue, 27 }  ,draw opacity=1 ]   (275,256) -- (294.13,264.6) ;
\draw [->]   (117,79) -- (224,58) ;
\draw [->]   (392,63) -- (497,101) ;
\draw [->]   (120,216) -- (227,241) ;
\draw[->]    (390,240) -- (497,192) ;

\end{tikzpicture}\]
\caption{}\label{cyclic and contraction operations}
\end{figure}

A morphism of coloured modular operads $f:\bP\rightarrow\bQ$ consists of an involutive function $f: \col(\bP)\rightarrow \col(\bQ)$ together with a family of $\Sigma_n$-equivariant maps \[\begin{tikzcd} \bP(c_1,\ldots, c_n) \arrow[r, "f_{\underline{c}}"] & \bQ(f(c_1),\ldots,f(c_n))\end{tikzcd}\] which commute with composition, contractions, and identities. We denote the category of all modular operads 
by $\ModOp (\bE)$, suppressing the $\bE$ when $\bE=\textbf{Set}$. 

\medskip

\begin{figure}[h!]
\[
\begin{tikzpicture}[x=0.75pt,y=0.75pt,yscale=-1,xscale=1]

\draw [color={rgb, 255:red, 0; green, 0; blue, 0 }  ,draw opacity=1 ]   (53,90) -- (87,131) ;
\draw  [fill={rgb, 255:red, 255; green, 255; blue, 255 }  ,fill opacity=1 ] (97.97,128.07) .. controls (106.25,128) and (113.02,134.66) .. (113.09,142.94) .. controls (113.16,151.22) and (106.5,158) .. (98.22,158.07) .. controls (89.93,158.14) and (83.16,151.48) .. (83.09,143.19) .. controls (83.02,134.91) and (89.68,128.14) .. (97.97,128.07) -- cycle ;
\draw [color={rgb, 255:red, 0; green, 0; blue, 0 }  ,draw opacity=1 ]   (88,156) -- (50,201) ;
\draw [color={rgb, 255:red, 0; green, 0; blue, 0 }  ,draw opacity=1 ]   (107,154) -- (147,202) ;
\draw    (110,134) -- (147,89) ;
\draw  [fill={rgb, 255:red, 255; green, 255; blue, 255 }  ,fill opacity=1 ] (286.97,92.07) .. controls (295.25,92) and (302.02,98.66) .. (302.09,106.94) .. controls (302.16,115.22) and (295.5,122) .. (287.22,122.07) .. controls (278.93,122.14) and (272.16,115.48) .. (272.09,107.19) .. controls (272.02,98.91) and (278.68,92.14) .. (286.97,92.07) -- cycle ;
\draw [color={rgb, 255:red, 0; green, 0; blue, 0 }  ,draw opacity=1 ]   (235,57) -- (273,100) ;
\draw    (300,98) -- (343,54) ;
\draw  [fill={rgb, 255:red, 255; green, 255; blue, 255 }  ,fill opacity=1 ] (288.97,175.07) .. controls (297.25,175) and (304.02,181.66) .. (304.09,189.94) .. controls (304.16,198.22) and (297.5,205) .. (289.22,205.07) .. controls (280.93,205.14) and (274.16,198.48) .. (274.09,190.19) .. controls (274.02,181.91) and (280.68,175.14) .. (288.97,175.07) -- cycle ;
\draw    (274,185) .. controls (266,176) and (265,130) .. (275,118) ;
\draw    (298,118) .. controls (306,131) and (308,171) .. (302,182) ;
\draw    (289.22,205.07) -- (290,250) ;
\draw  [fill={rgb, 255:red, 255; green, 255; blue, 255 }  ,fill opacity=1 ] (489.75,102.94) .. controls (497.09,102.88) and (503.08,109) .. (503.14,116.63) .. controls (503.2,124.25) and (497.31,130.48) .. (489.97,130.54) .. controls (482.64,130.6) and (476.64,124.48) .. (476.58,116.86) .. controls (476.52,109.24) and (482.42,103.01) .. (489.75,102.94) -- cycle ;
\draw    (434,67) -- (477.39,110.24) ;
\draw    (501.29,108.4) -- (539.36,67.92) ;
\draw    (495.09,130.48) -- (532.28,172.8) ;
\draw  [fill={rgb, 255:red, 255; green, 255; blue, 255 }  ,fill opacity=1 ] (541.99,170.1) .. controls (549.33,170.04) and (555.32,176.16) .. (555.38,183.79) .. controls (555.44,191.41) and (549.55,197.64) .. (542.21,197.7) .. controls (534.88,197.76) and (528.88,191.64) .. (528.82,184.02) .. controls (528.76,176.4) and (534.66,170.17) .. (541.99,170.1) -- cycle ;
\draw    (542.21,197.7) -- (542.02,251) ;
\draw    (548.22,170.04) .. controls (573.9,99.2) and (643.84,219.72) .. (556.19,193.04) ;

\draw (93,136) node [anchor=north west][inner sep=0.75pt]    {$p$};
\draw (279,99) node [anchor=north west][inner sep=0.75pt]    {$p_{v}$};
\draw (281,183) node [anchor=north west][inner sep=0.75pt]    {$p_{w}$};
\draw (482,108) node [anchor=north west][inner sep=0.75pt]    {$p_{v}$};
\draw (533,177) node [anchor=north west][inner sep=0.75pt]    {$p_{w}$};

\end{tikzpicture}\]

\caption{A graphical interpretation of operations in a modular operad}\label{operations modular}
\end{figure}
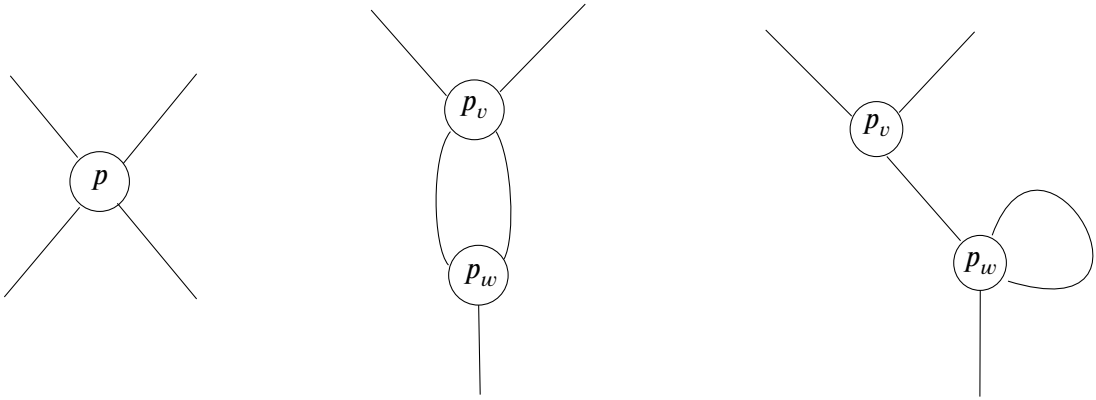

Returning to our informal description, operations of modular operads can be depicted as decorated connected graphs. In a modular operad operations might have more than one \emph{internal edge}, such as those depicted in Figure~\ref{operations modular}. In this interpretation, simply connected operations, such as those depicted on the right in Figure~\ref{cyclic operations} are obtained by cyclic operad composition. Contraction, depicted graphically, is obtained by identifying two boundary edges to form a new internal edge. These depictions will be given a precise meaning in terms of \emph{graphical maps} shortly. 

 As we remarked earlier, when we restrict to one-coloured, or monochrome, modular operads the involution on colour sets is trivial and, as such, the known examples of one-coloured modular operads meet the requirements of our definition. 
 
 \begin{example} 
 A \emph{compact closed category} is a $*$-autonomous category in which the dualizing functor is monoidal (\cite{Kelly_CCC}). Raynor mentions in \cite[Example 1.27]{Sophie_21} that \emph{involutive} compact closed categories are examples of modular operads. To our knowledge a full proof of this fact does not yet exist in the literature. 
 \end{example}

\subsection{Graphs} 

We are now ready to introduce a category $\bU$ whose objects are connected, \underline{u}ndirected graphs.\footnote{The category is called $\bU$ because "u" stands for undirected.} The category $\bU$ models modular operad in the sense that 
there is an equivalence of categories \[\ModOp \cong (\textbf{Set}^{\bU^{op}})_{Segal},\] where the right-hand side is a subcategory of $\bU$-presheaves which satisfy a \emph{strict Segal condition}. This equivalence of categories justifies our depiction of cyclic and modular operad operations with decorated graphs. 

In these lectures graphs are allowed to have ``loose ends''—meaning that it is not necessary for both ends (or either end) of an edge to touch a vertex. A typical example is in Figure~\ref{general graph}. To make this precise will use a combinatorial definition of graphs, \emph{Feynman graphs}, due to Joyal and Kock \cite{jk11}.  This model for graphs has the advantage that it is extremely easy to write down and the drawback in that it does not fully capture all the graphs we need for defining modular operads (see Remark~\ref{rmk: loop} below and Remark 1.1 in \cite{hry1}). We note that there are several other combinatorial definitions of graphs in the literature, all of which can be shown to be equivalent to the definition we use here by combining Proposition 15.2, Proposition 15.6, and Proposition 15.8 of \cite{BB}. 

\begin{figure}[h!]
\begin{centering} 
\begin{tikzpicture}[x=0.75pt,y=0.75pt,yscale=-.5,xscale=.5]

\draw    (284,29) .. controls (275,46) and (287,54) .. (231.2,73.97) ;
\draw  [fill={rgb, 255:red, 255; green, 255; blue, 255 }  ,fill opacity=1 ] (210.99,58.96) .. controls (218.78,56.14) and (227.38,60.18) .. (230.2,67.97) .. controls (233.01,75.76) and (228.98,84.36) .. (221.19,87.17) .. controls (213.4,89.99) and (204.8,85.96) .. (201.98,78.16) .. controls (199.17,70.37) and (203.2,61.78) .. (210.99,58.96) -- cycle ;
\draw    (226,84) .. controls (267.81,111.42) and (292.82,196.77) .. (222,214.99) ;
\draw    (164,4) .. controls (145,44) and (164.32,32.72) .. (206.13,60.14) ;
\draw    (226.45,234.13) .. controls (268.26,261.56) and (262,266) .. (266,283) ;
\draw  [fill={rgb, 255:red, 255; green, 255; blue, 255 }  ,fill opacity=1 ] (207.79,208.98) .. controls (215.58,206.16) and (224.18,210.2) .. (227,217.99) .. controls (229.81,225.78) and (225.78,234.38) .. (217.99,237.19) .. controls (210.2,240.01) and (201.6,235.98) .. (198.78,228.18) .. controls (195.97,220.39) and (200,211.79) .. (207.79,208.98) -- cycle ;
\draw    (201.47,209.14) .. controls (143.88,170.41) and (184.9,292.75) .. (208.58,240.59) ;
\draw    (82,162) .. controls (122,132) and (161.98,108.16) .. (201.98,78.16) ;
\draw    (391,72) -- (391,217) ;
\draw  [fill={rgb, 255:red, 255; green, 255; blue, 255 }  ,fill opacity=1 ] (63.99,158.96) .. controls (71.78,156.14) and (80.38,160.18) .. (83.2,167.97) .. controls (86.01,175.76) and (81.98,184.36) .. (74.19,187.17) .. controls (66.4,189.99) and (57.8,185.96) .. (54.98,178.16) .. controls (52.17,170.37) and (56.2,161.78) .. (63.99,158.96) -- cycle ;
\draw  [fill={rgb, 255:red, 255; green, 255; blue, 255 }  ,fill opacity=1 ] (331.4,136.07) .. controls (339.68,136.24) and (346.26,143.09) .. (346.09,151.37) .. controls (345.92,159.66) and (339.07,166.23) .. (330.78,166.06) .. controls (322.5,165.89) and (315.93,159.04) .. (316.09,150.76) .. controls (316.26,142.48) and (323.11,135.9) .. (331.4,136.07) -- cycle ;
\draw    (85,178) .. controls (136,173) and (210,163) .. (257,159) ;
\draw    (277,156) -- (316.09,150.76) ;
\draw    (330.78,166.06) -- (307,195) ;
\draw    (331.4,136.07) -- (303,112) ;
\draw   (423,239.5) .. controls (423,222.66) and (436.66,209) .. (453.5,209) .. controls (470.34,209) and (484,222.66) .. (484,239.5) .. controls (484,256.34) and (470.34,270) .. (453.5,270) .. controls (436.66,270) and (423,256.34) .. (423,239.5) -- cycle ;
\draw  [fill={rgb, 255:red, 255; green, 255; blue, 255 }  ,fill opacity=1 ] (452.61,200.01) .. controls (460.9,200.18) and (467.47,207.03) .. (467.3,215.31) .. controls (467.13,223.59) and (460.28,230.17) .. (452,230) .. controls (443.72,229.83) and (437.14,222.98) .. (437.31,214.7) .. controls (437.48,206.41) and (444.33,199.84) .. (452.61,200.01) -- cycle ;

\end{tikzpicture}\end{centering}
\caption{}\label{general graph}
\end{figure} 

Edges of our graphs are all comprised of two \emph{distinct} ``half edges''  which you can picture as a copy of the interval $(0,1)$ equipped with a chosen orientation. We denote the set of half edges by $A$.\footnote{Half edges are also often called \emph{arcs} in the literature, hence the use of the letter $A$ for this set.} This set is equipped with a free involution $i$ which identifies a pair of half edges with the opposite orientation; edges are the orbits of this involution. Half edges are attached to vertices via a partial function $t: A \rightarrow V$, where $V$ denotes the set of vertices.  Not all half edges will be attached to a vertex, so to distinguish the set of those edges which are adjacent to vertices, we write $D \subseteq A$ for those half edges which are in the \emph{domain} of the function $t$. 

\begin{definition}\label{def: graphs}  A \emph{graph} $G$ consists of: 
\begin{itemize} 
\item  a diagram of finite sets 
\[\begin{tikzcd} 
\arrow[loop left]{l}{i}A & \arrow[l, swap, "s"] D \arrow[r, "t"] & V
\end{tikzcd} \] where
\begin{itemize} 
\item $i$ is a fixed point-free involution and
\item $s$ is a monomorphism. 
\end{itemize} 
\end{itemize} 
\end{definition}

In Figure~\ref{example: first graph} the graph depicted has the set of arcs $A=\{1,2,\ldots,9,10\}$ and vertices $V=\{v_1,v_2\}$.  The arcs adjacent to the vertices are $D=\{2,3,9,4,5,6,7\}$ and the involution on $A$ interchanges $2n$ and $2n-1$ for $n=1,2,3,4,5$. 
\begin{figure}[h!]

\begin{tikzpicture}[x=0.75pt,y=0.75pt,yscale=-.75,xscale=.75]

\draw  [fill={rgb, 255:red, 255; green, 255; blue, 255 }  ,fill opacity=1 ] (241.97,125.07) .. controls (250.25,125) and (257.02,131.66) .. (257.09,139.94) .. controls (257.16,148.22) and (250.5,155) .. (242.22,155.07) .. controls (233.93,155.14) and (227.16,148.48) .. (227.09,140.19) .. controls (227.02,131.91) and (233.68,125.14) .. (241.97,125.07) -- cycle ;
\draw    (179,86) -- (228,133) ;
\draw    (255,131) -- (298,87) ;
\draw    (248,155) -- (290,201) ;
\draw  [fill={rgb, 255:red, 255; green, 255; blue, 255 }  ,fill opacity=1 ] (300.97,198.07) .. controls (309.25,198) and (316.02,204.66) .. (316.09,212.94) .. controls (316.16,221.22) and (309.5,228) .. (301.22,228.07) .. controls (292.93,228.14) and (286.16,221.48) .. (286.09,213.19) .. controls (286.02,204.91) and (292.68,198.14) .. (300.97,198.07) -- cycle ;
\draw    (301.22,228.07) -- (301,286) ;
\draw    (308,198) .. controls (337,121) and (416,252) .. (317,223) ;

\draw (233,132) node [anchor=north west][inner sep=0.75pt]    {$v_{1}$};
\draw (291,205) node [anchor=north west][inner sep=0.75pt]    {$v_{2}$};
\draw (173,91.4) node [anchor=north west][inner sep=0.75pt] [font=\footnotesize] {$1$};
\draw (208,126.4) node [anchor=north west][inner sep=0.75pt]   [font=\footnotesize] {$2$};
\draw (239,160.4) node [anchor=north west][inner sep=0.75pt] [font=\footnotesize]   {$3$};
\draw (271,195.4) node [anchor=north west][inner sep=0.75pt]  [font=\footnotesize] {$4$};
\draw (321,182) node [anchor=north west][inner sep=0.75pt] [font=\footnotesize]{$5$};
\draw (321,207) node [anchor=north west][inner sep=0.75pt] [font=\footnotesize]{$6$};
\draw (285,237.4) node [anchor=north west][inner sep=0.75pt]  [font=\footnotesize]{$7$};
\draw (285,267.4) node [anchor=north west][inner sep=0.75pt]  [font=\footnotesize]{$8$};
\draw (290,94.4) node [anchor=north west][inner sep=0.75pt]  [font=\footnotesize]{$10$};
\draw (264,126.4) node [anchor=north west][inner sep=0.75pt]  [font=\footnotesize]{$9$};

\end{tikzpicture}
\caption{}\label{example: first graph} 
\end{figure} 

The involution on half-edges determines the \emph{edges} of a graph as follows.  If we write the action of the involution as $i: a\mapsto a^{\dagger}$, an \emph{edge} of the graph $G$ is an $i$-orbit $[a, a^{\dagger}]$. We write $E(G)=A/i$ for the set of edges of the graph $G$. An \emph{internal edge} is an edge of the form $[b, b^{\dagger}]$ where both $b$ and $b^{\dagger}$ are in $D$. In other words an internal edge is an edge $[b, b^{\dagger}]$ where both $b$ and $b^{\dagger}$ are adjacent to vertices. In Figure~\ref{example: first graph} the two internal edges are given by the orbits $e_1=[3,4]$ and $e_2=[5,6]$.  The \emph{boundary} of a graph is the set of half edges not adjacent to a vertex, $\partial(G) =A \setminus D.$ The \emph{neighbourhood} of a vertex $v\in V(G)$ consists of the half-edges which are adjacent to the vertex, $\nb(v)=t^{-1}(v)\subseteq D.$

\begin{example} There are several graphs which warrant special names. 

\begin{enumerate} 
\item The exceptional \emph{edge}, $G=\updownarrow$ has exactly two arcs $A=\{a, a^{\dagger}\}$ and $D=V=\emptyset$. The boundary of the edge $\partial(\updownarrow)=A$. 

\item For $n\geq 0$, we write $G=\medwhitestar_n$ for the \emph{$n$-star}. This graph has $A= \{1, 1^{\dagger}, \ldots, n, n^{\dagger}\}$, $D = \{1,\ldots, n\}$ and $V=\{v\}$. The boundary of  $G=\medwhitestar_n$ is the set $\{1^{\dagger},\ldots, n^{\dagger}\}$ and the neighbourhood of the vertex $v$ is $\{1,\ldots, n\}$. The $4$-star is depicted in Figure~\ref{edge and star}

\end{enumerate} 
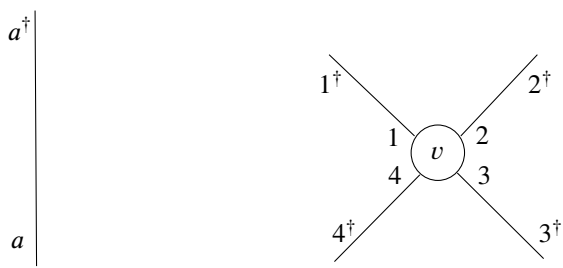
\begin{figure}[h!]

\begin{tikzpicture}[x=0.75pt,y=0.75pt,yscale=-.75,xscale=.75]

\draw    (139,80) -- (140,251) ;
\draw  [fill={rgb, 255:red, 255; green, 255; blue, 255 }  ,fill opacity=1 ] (409.69,156.59) .. controls (419.6,156.5) and (427.71,164.74) .. (427.8,174.98) .. controls (427.88,185.22) and (419.91,193.59) .. (409.99,193.68) .. controls (400.07,193.76) and (391.97,185.53) .. (391.88,175.29) .. controls (391.8,165.05) and (399.77,156.68) .. (409.69,156.59) -- cycle ;
\draw    (336.71,109.52) -- (394.17,163.92) ;
\draw    (425.29,163.92) -- (476.77,109.52) ;
\draw    (422.9,188.65) -- (481,246) ;
\draw    (397.76,189.89) -- (340.3,248) ;

\draw (121,229.4) node [anchor=north west][inner sep=0.75pt] [font=\footnotesize]{$a$};
\draw (119,82.4) node [anchor=north west][inner sep=0.75pt]  [font=\footnotesize]{$a^{\dagger }$};
\draw (403,170) node [anchor=north west][inner sep=0.75pt] [font=\footnotesize]{$v$};
\draw (328,118) node [anchor=north west][inner sep=0.75pt] [font=\footnotesize]{$1^{\dagger }$};
\draw (373.8,156.32) node [anchor=north west][inner sep=0.75pt]  [font=\footnotesize]{$1$};
\draw (433.66,155.09) node [anchor=north west][inner sep=0.75pt]  [font=\footnotesize] {$2$};
\draw (469,118) node [anchor=north west][inner sep=0.75pt]  [font=\footnotesize]{$2^{\dagger }$};
\draw (435,182) node [anchor=north west][inner sep=0.75pt] [font=\footnotesize]{$3$};
\draw (476,218) node [anchor=north west][inner sep=0.75pt]  [font=\footnotesize] {$3^{\dagger }$};
\draw (375,182) node [anchor=north west][inner sep=0.75pt]  [font=\footnotesize]{$4$};
\draw (337.28,218) node [anchor=north west][inner sep=0.75pt]  [font=\footnotesize]{$4^{\dagger }$};

\end{tikzpicture}
\caption{The exceptional edge $\updownarrow$ and the $4$-star $\medwhitestar_4$. }\label{edge and star}
\end{figure}
\end{example}

\begin{figure}[h!]

\begin{tikzpicture}[x=0.75pt,y=0.75pt,yscale=-.5,xscale=.5]

\draw   (419,169) .. controls (419,137.52) and (444.52,112) .. (476,112) .. controls (507.48,112) and (533,137.52) .. (533,169) .. controls (533,200.48) and (507.48,226) .. (476,226) .. controls (444.52,226) and (419,200.48) .. (419,169) -- cycle ;
\draw   (127,169.8) .. controls (127,137.96) and (153.86,112.15) .. (187,112.15) .. controls (220.14,112.15) and (247,137.96) .. (247,169.8) .. controls (247,201.64) and (220.14,227.45) .. (187,227.45) .. controls (153.86,227.45) and (127,201.64) .. (127,169.8) -- cycle ;
\draw  [fill={rgb, 255:red, 255; green, 255; blue, 255 }  ,fill opacity=1 ] (187.43,91.74) .. controls (198.66,91.96) and (207.56,101.27) .. (207.33,112.55) .. controls (207.09,123.82) and (197.79,132.78) .. (186.57,132.56) .. controls (175.34,132.34) and (166.44,123.02) .. (166.67,111.75) .. controls (166.91,100.47) and (176.21,91.51) .. (187.43,91.74) -- cycle ;
\draw  [fill={rgb, 255:red, 255; green, 255; blue, 255 }  ,fill opacity=1 ] (185.31,202.17) .. controls (196.53,202.4) and (205.44,211.71) .. (205.2,222.99) .. controls (204.96,234.26) and (195.67,243.22) .. (184.44,243) .. controls (173.22,242.78) and (164.31,233.46) .. (164.55,222.18) .. controls (164.79,210.91) and (174.08,201.95) .. (185.31,202.17) -- cycle ;

\end{tikzpicture}
\caption{A loop with $2$ vertices and a nodeless loop.}\label{loops}
\end{figure}
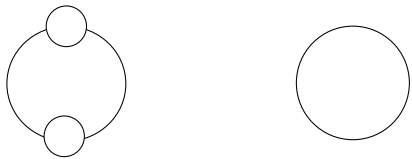

Many of our graphs have empty boundary such as the loop with $n$-vertices depicted in Figure~\ref{loops}. 

\begin{remark}\label{rmk: loop}
The one graph that we \emph{cannot} describe using our chosen formalism of graphs is the \emph{nodeless loop}, depicted on the right of Figure~\ref{loops}. As we explain in Remark 1.1 in \cite{hry1}, if one attempts to describe the nodeless loop in our chosen formalism, we end up with a graph which cannot be distinguished from the \emph{edge}. The technical point will not play a further role in this lecture series, but it plays an essential role in defining the monad for modular operads in \cite{hry2}. 
\end{remark}

Every graph $G$ has an associated star, $\medwhitestar_G$, determined by its boundary. The graph $\medwhitestar_G$ is the one-vertex graph with $A =\partial(G)\sqcup \partial(G)^{\dagger}$ and $D=\partial(G)^{\dagger}$. By definition, $\partial(\medwhitestar_G) = A \setminus D=\partial(G)$ and the neighbourhood of the unique vertex is $D=\partial(G)^{\dagger}$. Similarly, every vertex of a graph $G$ has an associated star, $\medwhitestar_v$,  with $V(G)=\{v\}$, $D=\nb(v)$, and $A =\nb(v)\sqcup \nb(v)^{\dagger}$. The boundary of $\partial(\medwhitestar_v) = \nb(v)^{\dagger}$. 

\begin{exercise}\
Draw a graph $G$ with $A=\{1,2,3,4,5,6,7,8\}$, $D=\{1,2,3,4,5,6,7,8\}$, $V=\{v_1, v_2, v_3, v_4\}$ and $i(n)=n-1$ for $n=2,4,6,8$.  For the graph $G$ you have drawn, write down $\medwhitestar_G$ and $\medwhitestar_v$ for each $v\in V(G)$.
\end{exercise}

\subsubsection{Morphisms of Graphs} 
Our definition of graphical morphisms is designed to capture the units, composition, and contraction operations of modular operads.  The first definition of graphical map we give in this lecture series is not, necessarily, the best or most practical definition of graphical map. Instead, we have chosen to present the material in such a way as to motivate how one might arrive at this definition: start with the most ``obvious'' description of a map between graphs and then modify morphisms until we get to our ideal definition. 
 
A graph $G$ (Definition~\ref{def: graphs}) is a diagram in the category of finite sets in the shape of 
\[\begin{tikzcd} \calI : = & \arrow[loop left]{l}{i}\bullet & \arrow[l, swap, "s"] \bullet \arrow[r, "t"] & \bullet \end{tikzcd} \] where the arrow $s$ is sent to a monomorphism and the generating endomorphism $i$ is a free involution.  A morphism between two graphs should preserve some structure of the graphs, e.g., a vertex with four adjacent edges should not map to a vertex with three adjacent edges.  This leads to the first guess for a definition of graphical map: a graphical map should be maps in the functor category $\mathbf{FinSet}^{\calI}$ which preserve the local structure of graphs. 

\begin{definition}\label{def: embedding}
Let $G$ and $G'$ be two connected graphs.  A natural transformation $\phi: G \rightarrow G'$ is called an \emph{embedding} if the right-hand square of: 
 \[\begin{tikzcd} 
\arrow[loop left]{l}{i}A \arrow[d, "\phi_A"]& \arrow[l, swap, "s"] D \arrow[d, "\phi_D"] \arrow[r, "t"] & V \arrow[d, "\phi_V"]  \\
\arrow[loop left]{l}{i'}A' & \arrow[l, swap, "s'"] D' \arrow[r, "t' "] & V'
\end{tikzcd} \] is a pullback \emph{and} the map $V\rightarrow V' $  a monomorphism.
\end{definition} 

The pullback condition of Definition~\ref{def: embedding} makes sure that the local information i.e., the neighborhoods of vertices, is preserved.\footnote{A natural transformation $\phi:G\rightarrow G'$ which only satisfies the pullback condition of Definition~\ref{def: embedding} is called \`{e}tale in \cite{jk11}. \`{E}tale maps play a key role in the description of ``graphical species'' in \cite{jk11} and \cite{Sophie_21}, but are not quite what we want. See, for example, Remark 2.4~\cite{phil_graphs} for more details on this fine point.} The requirement that $V\rightarrow V' $ be a monomorphism \emph{almost} implies that every graphical map is an inclusion of a \emph{subgraph}. In particular, for every vertex $v$ of a graph $G$ there is a canonical embedding $\medwhitestar_v\longrightarrow G$: 
\begin{equation}\label{star_inclusion}\tag{$\bigstar$}\begin{tikzcd} 
 \nb(v)\sqcup \nb(v)^{\dagger} \arrow[d]& \arrow[l, swap, "s"] \nb(v) \arrow[d] \arrow[r, "t"] & \{v\} \arrow[d]  \\
 A & \arrow[l, swap, "s'"] D \arrow[r, "t' "] & V.
\end{tikzcd}
\end{equation}The left-hand map in this diagram is just the inclusion $\nb(v) \rightarrow D\rightarrow A$ on the first component, while the second component (which is forced by compatibility with the involution) sends $a^{\dagger}$ to $ia$. Similarly, every edge of a graph $G$ corresponds to an embedding $\updownarrow\longrightarrow G$.

\begin{exercise}
Write down the explicit natural transformation for an edge inclusion $\updownarrow\longrightarrow G$. 
\end{exercise}

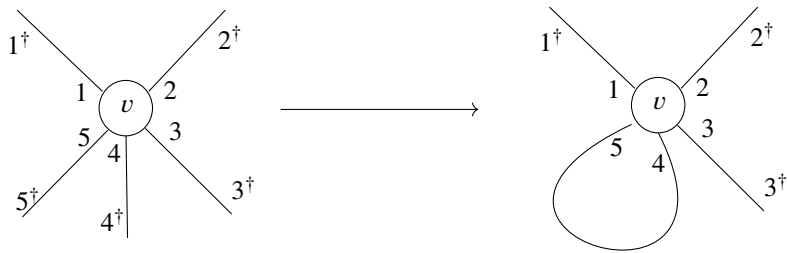
\begin{figure}[h]
\[\begin{tikzpicture}[x=0.75pt,y=0.75pt,yscale=-.75,xscale=.75]

\draw  [fill={rgb, 255:red, 255; green, 255; blue, 255 }  ,fill opacity=1 ] (120.69,110.59) .. controls (130.6,110.5) and (138.71,118.74) .. (138.8,128.98) .. controls (138.88,139.22) and (130.91,147.59) .. (120.99,147.68) .. controls (111.07,147.76) and (102.97,139.53) .. (102.88,129.29) .. controls (102.8,119.05) and (110.77,110.68) .. (120.69,110.59) -- cycle ;
\draw    (47.71,63.52) -- (105.17,117.92) ;
\draw    (136.29,117.92) -- (187.77,63.52) ;
\draw    (133.9,142.65) -- (192,200) ;
\draw    (108.76,143.89) -- (51.3,202) ;
\draw    (120.99,147.68) -- (122,216) ;
\draw  [fill={rgb, 255:red, 255; green, 255; blue, 255 }  ,fill opacity=1 ] (478.69,108.59) .. controls (488.6,108.5) and (496.71,116.74) .. (496.8,126.98) .. controls (496.88,137.22) and (488.91,145.59) .. (478.99,145.68) .. controls (469.07,145.76) and (460.97,137.53) .. (460.88,127.29) .. controls (460.8,117.05) and (468.77,108.68) .. (478.69,108.59) -- cycle ;
\draw    (405.71,61.52) -- (463.17,115.92) ;
\draw    (494.29,115.92) -- (545.77,61.52) ;
\draw    (491.9,140.65) -- (550,198) ;
\draw    (461,140) .. controls (306,213) and (549,283) .. (478.99,145.68) ;
\draw[->]    (225,130) -- (358,130) ;

\draw (115.32,122) node [anchor=north west][inner sep=0.75pt] [font=\footnotesize]   {$v$};
\draw (39.11,75.82) node [anchor=north west][inner sep=0.75pt]  [font=\footnotesize]  {$1^{\dagger }$};
\draw (84.8,110.32) node [anchor=north west][inner sep=0.75pt]  [font=\footnotesize]  {$1$};
\draw (144.66,109.09) node [anchor=north west][inner sep=0.75pt] [font=\footnotesize]   {$2$};
\draw (181.16,72.11) node [anchor=north west][inner sep=0.75pt]   [font=\footnotesize] {$2^{\dagger }$};
\draw (148.25,136.29) node [anchor=north west][inner sep=0.75pt]  [font=\footnotesize]  {$3$};
\draw (190.34,172.84) node [anchor=north west][inner sep=0.75pt]  [font=\footnotesize]  {$3^{\dagger }$};
\draw (107,151) node [anchor=north west][inner sep=0.75pt] [font=\footnotesize]   {$4$};
\draw (102.28,194.26) node [anchor=north west][inner sep=0.75pt][font=\footnotesize]    {$4^{\dagger }$};
\draw (86.8,140.32) node [anchor=north west][inner sep=0.75pt]  [font=\footnotesize]  {$5$};
\draw (45.28,180.26) node [anchor=north west][inner sep=0.75pt]  [font=\footnotesize]  {$5^{\dagger }$};
\draw (473.32,120) node [anchor=north west][inner sep=0.75pt]  [font=\footnotesize]  {$v$};
\draw (397.11,73.82) node [anchor=north west][inner sep=0.75pt] [font=\footnotesize]   {$1^{\dagger }$};
\draw (442.8,108.32) node [anchor=north west][inner sep=0.75pt]  [font=\footnotesize]  {$1$};
\draw (502.66,107.09) node [anchor=north west][inner sep=0.75pt]  [font=\footnotesize]  {$2$};
\draw (539.16,70.11) node [anchor=north west][inner sep=0.75pt]  [font=\footnotesize]  {$2^{\dagger }$};
\draw (506.25,134.29) node [anchor=north west][inner sep=0.75pt] [font=\footnotesize]   {$3$};
\draw (548.34,170.84) node [anchor=north west][inner sep=0.75pt]  [font=\footnotesize]  {$3^{\dagger }$};
\draw (473,158) node [anchor=north west][inner sep=0.75pt] [font=\footnotesize]   {$4$};
\draw (444.8,147.32) node [anchor=north west][inner sep=0.75pt]  [font=\footnotesize]  {$5$};

\end{tikzpicture}\]
\caption{Embedding a $5$-star into a contracted $5$-star.}\label{fig:contracted cor embed}
\end{figure}

Our definition of embedding, however, is more general than a subgraph inclusion, because an embedding is not necessarily injective on the set of half edges.  For example, Figure~\ref{fig:contracted cor embed} depicts an embedding $\phi:\medwhitestar_5\longrightarrow G$, 
\[\begin{tikzcd} 
\{1, 2, 3, 4, 5, 1^{\dagger}, 2^{\dagger}, 3^{\dagger}, 4^{\dagger},5^{\dagger}\} \arrow[d]& \arrow[l, swap, "s"] \{1, 2, 3, 4, 5\} \arrow[d] \arrow[r, "t"] & \{v\} \arrow[d]  \\
 \{1, 2, 3, 4, 5, 1^{\dagger}, 2^{\dagger}, 3^{\dagger}\} & \arrow[l, swap, "s'"] \{1, 2, 3, 4, 5\} \arrow[r, "t' "] & \{v\}.
\end{tikzcd} \] The graph $G$ in Figure~\ref{fig:contracted cor embed} is called a \emph{contracted star}. Explicitly, $G$ is the graph with one vertex $v$, set of half edges $A = \{1, 2, 3, 4, 5, 1^{\dagger}, 2^{\dagger}, 3^{\dagger}\}$ and $D=\{1, 2, 3, 4, 5\}$. The involution $i(n)=n^{\dagger}$, $n=1,2,3$ and $i(4)=5$ (and $i(5)=4$).

\medskip 

Let $\Emb(G)$ denote the set of isomorphism classes of embeddings into $G$. We are now ready to define graphical maps: 

\begin{definition}\label{def: graphical maps}
A graphical map $\varphi: G\rightarrow G'$ consists of: 
\begin{itemize} 
\item a map of involutive sets $\varphi_0: A_G \rightarrow A_{G'}$;
\item a function $\varphi_1: V_{G}\rightarrow \Emb(G')$ satisfying the following conditions:
\begin{itemize} 
\item The embeddings $\varphi_1(v)$ do not \emph{overlap} at vertices -- no vertex in $G'$ is contained in two graphs $\varphi_1(v)$ and $\varphi_1(v')$;
\item For each v, we have a (necessarily unique) bijection making the diagram: 
\[\begin{tikzcd} \nb(v)\arrow[d, "\cong", swap] \arrow[r, "i"] & A\arrow[d, "\varphi_0"] \\
\partial(\varphi_1(v)) \arrow[r]& A'
\end{tikzcd}\] commute, where the top map $i$ is the restriction of the involution on A. 

\item If $\partial(G) =\emptyset$, then there exists a $v$ in $V$ so that $\varphi_1(v) \not= \ \updownarrow$.
\end{itemize}
\end{itemize} 
\end{definition}

The first two conditions on $\varphi_1: V\rightarrow \Emb(G')$ imply that a map $\varphi: G\rightarrow G'$ is obtained by ``blowing up'', or replacing, the vertices of $G$ with another graph $H_v$.\footnote{The notion of ``blowing up'' a vertex can be made precise using the language of \emph{graph substitution} which is described for Feynman graphs in Construction 1.18 \cite{hry1}.} The requirements on the embedding \[\begin{tikzcd}H_v\arrow[r, "\varphi_1(v)", hook] &G'\end{tikzcd}\] guarantee that we replace the vertex $v$ by a graph $H_v$ in such a way that we have an isomorphism $\medwhitestar_{H_v} \cong \medwhitestar_v.$ The third condition is about avoiding the collapse into a \emph{nodeless loop}. In Figure~\ref{map}, we have circled the subgraph $H_v=\varphi_1(v)$ in red. 

\begin{exercise} 
Write down an explicit description of the graphical map depicted in Figure~\ref{map}. 
\end{exercise} 

\begin{figure} 
\[
\begin{tikzpicture}[x=0.75pt,y=0.75pt,yscale=-1,xscale=1]

\draw    (114.84,123.97) -- (143.88,188.12) ;
\draw    (73.89,176.06) -- (107.15,133.2) ;
\draw    (51.03,124.07) -- (102.31,123.97) ;
\draw    (143.88,188.12) -- (115.59,245.03) ;
\draw  [color={rgb, 255:red, 4; green, 146; blue, 194 }  ,draw opacity=1 ][fill={rgb, 255:red, 255; green, 255; blue, 255 }  ,fill opacity=1 ][dash pattern={on 4.5pt off 4.5pt}] (102.31,123.97) .. controls (102.31,117.05) and (107.92,111.44) .. (114.84,111.44) .. controls (121.76,111.44) and (127.37,117.05) .. (127.37,123.97) .. controls (127.37,130.89) and (121.76,136.5) .. (114.84,136.5) .. controls (107.92,136.5) and (102.31,130.89) .. (102.31,123.97) -- cycle ;
\draw    (122.04,114.52) .. controls (166.89,52.51) and (199.98,153.31) .. (126.66,129.14) ;
\draw    (170.22,268.12) -- (115.59,245.03) ;
\draw  [color={rgb, 255:red, 237; green, 130; blue, 14 }  ,draw opacity=1 ][fill={rgb, 255:red, 255; green, 255; blue, 255 }  ,fill opacity=1 ][dash pattern={on 4.5pt off 4.5pt}] (156,191.3) .. controls (154.25,197.99) and (147.41,202) .. (140.71,200.25) .. controls (134.01,198.5) and (130.01,191.65) .. (131.76,184.95) .. controls (133.51,178.26) and (140.36,174.25) .. (147.05,176) .. controls (153.75,177.75) and (157.76,184.6) .. (156,191.3) -- cycle ;
\draw    (121.93,256.52) .. controls (159.11,323.41) and (53.78,310.69) .. (106.73,254.51) ;
\draw  [color={rgb, 255:red, 140; green, 219; blue, 54 }  ,draw opacity=1 ][fill={rgb, 255:red, 255; green, 255; blue, 255 }  ,fill opacity=1 ][dash pattern={on 4.5pt off 4.5pt}] (127.71,248.21) .. controls (125.96,254.9) and (119.12,258.91) .. (112.42,257.16) .. controls (105.72,255.41) and (101.72,248.56) .. (103.47,241.86) .. controls (105.22,235.17) and (112.07,231.16) .. (118.76,232.91) .. controls (125.46,234.66) and (129.47,241.51) .. (127.71,248.21) -- cycle ;
\draw  [color={rgb, 255:red, 4; green, 146; blue, 194 }  ,draw opacity=1 ][dash pattern={on 4.5pt off 4.5pt}] (42.04,70.45) .. controls (42.04,50.91) and (57.88,35.07) .. (77.42,35.07) .. controls (96.96,35.07) and (112.8,50.91) .. (112.8,70.45) .. controls (112.8,89.99) and (96.96,105.83) .. (77.42,105.83) .. controls (57.88,105.83) and (42.04,89.99) .. (42.04,70.45) -- cycle ;
\draw  [color={rgb, 255:red, 237; green, 130; blue, 14 }  ,draw opacity=1 ][dash pattern={on 4.5pt off 4.5pt}] (172.85,175.51) .. controls (172.85,160.37) and (185.12,148.1) .. (200.26,148.1) .. controls (215.4,148.1) and (227.68,160.37) .. (227.68,175.51) .. controls (227.68,190.65) and (215.4,202.93) .. (200.26,202.93) .. controls (185.12,202.93) and (172.85,190.65) .. (172.85,175.51) -- cycle ;
\draw    (453.36,104.93) -- (480.57,165.38) ;
\draw    (349.18,152.84) -- (379.77,113.42) ;
\draw    (328.16,105.03) -- (375.32,104.93) ;
\draw    (398.37,104.93) -- (441.84,104.93) ;
\draw    (421.52,51.15) -- (386.85,104.93) ;
\draw    (421.52,51.15) -- (453.36,104.93) ;
\draw  [fill={rgb, 255:red, 255; green, 255; blue, 255 }  ,fill opacity=1 ] (375.32,104.93) .. controls (375.32,98.57) and (380.48,93.41) .. (386.85,93.41) .. controls (393.21,93.41) and (398.37,98.57) .. (398.37,104.93) .. controls (398.37,111.3) and (393.21,116.45) .. (386.85,116.45) .. controls (380.48,116.45) and (375.32,111.3) .. (375.32,104.93) -- cycle ;
\draw  [fill={rgb, 255:red, 255; green, 255; blue, 255 }  ,fill opacity=1 ] (410,51.15) .. controls (410,44.79) and (415.15,39.63) .. (421.52,39.63) .. controls (427.88,39.63) and (433.04,44.79) .. (433.04,51.15) .. controls (433.04,57.52) and (427.88,62.68) .. (421.52,62.68) .. controls (415.15,62.68) and (410,57.52) .. (410,51.15) -- cycle ;
\draw  [fill={rgb, 255:red, 255; green, 255; blue, 255 }  ,fill opacity=1 ] (441.84,104.93) .. controls (441.84,98.57) and (447,93.41) .. (453.36,93.41) .. controls (459.73,93.41) and (464.89,98.57) .. (464.89,104.93) .. controls (464.89,111.3) and (459.73,116.45) .. (453.36,116.45) .. controls (447,116.45) and (441.84,111.3) .. (441.84,104.93) -- cycle ;
\draw    (91.99,70.14) -- (103.35,95.25) ;
\draw    (47.65,90.52) -- (60.67,73.75) ;
\draw    (42.04,70.45) -- (58.78,70.14) ;
\draw    (91.99,70.14) -- (112,78.43) ;
\draw    (91.99,70.14) -- (109.19,53.12) ;
\draw    (68.59,70.14) -- (87.08,70.14) ;
\draw    (78.44,47.25) -- (63.68,70.14) ;
\draw    (78.44,47.25) -- (91.99,70.14) ;
\draw  [fill={rgb, 255:red, 255; green, 255; blue, 255 }  ,fill opacity=1 ] (58.78,70.14) .. controls (58.78,67.43) and (60.97,65.23) .. (63.68,65.23) .. controls (66.39,65.23) and (68.59,67.43) .. (68.59,70.14) .. controls (68.59,72.85) and (66.39,75.04) .. (63.68,75.04) .. controls (60.97,75.04) and (58.78,72.85) .. (58.78,70.14) -- cycle ;
\draw  [fill={rgb, 255:red, 255; green, 255; blue, 255 }  ,fill opacity=1 ] (73.53,47.25) .. controls (73.53,44.55) and (75.73,42.35) .. (78.44,42.35) .. controls (81.14,42.35) and (83.34,44.55) .. (83.34,47.25) .. controls (83.34,49.96) and (81.14,52.16) .. (78.44,52.16) .. controls (75.73,52.16) and (73.53,49.96) .. (73.53,47.25) -- cycle ;
\draw  [fill={rgb, 255:red, 255; green, 255; blue, 255 }  ,fill opacity=1 ] (87.08,70.14) .. controls (87.08,67.43) and (89.28,65.23) .. (91.99,65.23) .. controls (94.69,65.23) and (96.89,67.43) .. (96.89,70.14) .. controls (96.89,72.85) and (94.69,75.04) .. (91.99,75.04) .. controls (89.28,75.04) and (87.08,72.85) .. (87.08,70.14) -- cycle ;

\draw    (187.17,152.19) -- (198.74,177.74) ;
\draw    (198.74,177.74) -- (187.47,200.41) ;
\draw    (201.73,173.18) .. controls (219.59,148.48) and (232.77,188.63) .. (203.56,179) ;
\draw  [fill={rgb, 255:red, 255; green, 255; blue, 255 }  ,fill opacity=1 ] (193.75,177.74) .. controls (193.75,174.98) and (195.98,172.75) .. (198.74,172.75) .. controls (201.49,172.75) and (203.73,174.98) .. (203.73,177.74) .. controls (203.73,180.5) and (201.49,182.73) .. (198.74,182.73) .. controls (195.98,182.73) and (193.75,180.5) .. (193.75,177.74) -- cycle ;
\draw    (480.57,165.38) -- (441.22,250.03) ;
\draw    (487.6,154.65) .. controls (529.62,96.56) and (560.61,190.99) .. (491.93,168.35) ;
\draw  [fill={rgb, 255:red, 255; green, 255; blue, 255 }  ,fill opacity=1 ] (468.83,165.38) .. controls (468.83,158.89) and (474.08,153.64) .. (480.57,153.64) .. controls (487.05,153.64) and (492.31,158.89) .. (492.31,165.38) .. controls (492.31,171.86) and (487.05,177.12) .. (480.57,177.12) .. controls (474.08,177.12) and (468.83,171.86) .. (468.83,165.38) -- cycle ;
\draw    (491.23,271.46) -- (440.85,250.17) ;
\draw    (441.22,250.03) -- (390.84,228.74) ;
\draw  [fill={rgb, 255:red, 255; green, 255; blue, 255 }  ,fill opacity=1 ] (379.26,228.74) .. controls (379.26,222.34) and (384.45,217.16) .. (390.84,217.16) .. controls (397.24,217.16) and (402.43,222.34) .. (402.43,228.74) .. controls (402.43,235.14) and (397.24,240.33) .. (390.84,240.33) .. controls (384.45,240.33) and (379.26,235.14) .. (379.26,228.74) -- cycle ;
\draw  [fill={rgb, 255:red, 255; green, 255; blue, 255 }  ,fill opacity=1 ] (429.27,250.17) .. controls (429.27,243.78) and (434.45,238.59) .. (440.85,238.59) .. controls (447.25,238.59) and (452.44,243.78) .. (452.44,250.17) .. controls (452.44,256.57) and (447.25,261.76) .. (440.85,261.76) .. controls (434.45,261.76) and (429.27,256.57) .. (429.27,250.17) -- cycle ;
\draw    (447.15,260.8) .. controls (481.33,322.31) and (384.47,310.61) .. (433.17,258.95) ;
\draw    (459.98,95.54) .. controls (501.23,38.51) and (531.65,131.21) .. (464.23,108.98) ;
\draw  [color={rgb, 255:red, 4; green, 146; blue, 194 }  ,draw opacity=1 ][dash pattern={on 4.5pt off 4.5pt}] (362.3,89.16) .. controls (362.3,57.92) and (387.62,32.59) .. (418.87,32.59) .. controls (450.11,32.59) and (475.44,57.92) .. (475.44,89.16) .. controls (475.44,120.41) and (450.11,145.73) .. (418.87,145.73) .. controls (387.62,145.73) and (362.3,120.41) .. (362.3,89.16) -- cycle ;
\draw  [color={rgb, 255:red, 237; green, 130; blue, 14 }  ,draw opacity=1 ][dash pattern={on 4.5pt off 4.5pt}] (455.36,165.38) .. controls (455.36,151.45) and (466.64,140.17) .. (480.57,140.17) .. controls (494.49,140.17) and (505.78,151.45) .. (505.78,165.38) .. controls (505.78,179.3) and (494.49,190.59) .. (480.57,190.59) .. controls (466.64,190.59) and (455.36,179.3) .. (455.36,165.38) -- cycle ;
\draw  [color={rgb, 255:red, 140; green, 219; blue, 54 }  ,draw opacity=1 ][dash pattern={on 4.5pt off 4.5pt}] (368.91,232.75) .. controls (368.91,204.39) and (391.9,181.39) .. (420.27,181.39) .. controls (448.64,181.39) and (471.63,204.39) .. (471.63,232.75) .. controls (471.63,261.12) and (448.64,284.12) .. (420.27,284.12) .. controls (391.9,284.12) and (368.91,261.12) .. (368.91,232.75) -- cycle ;
\draw    (79.07,214.04) -- (69.3,235.44) ;
\draw    (79.72,240.37) -- (69.13,235.51) ;
\draw    (69.3,235.44) -- (46.16,225.67) ;
\draw  [fill={rgb, 255:red, 255; green, 255; blue, 255 }  ,fill opacity=1 ] (40.84,225.67) .. controls (40.84,222.73) and (43.22,220.35) .. (46.16,220.35) .. controls (49.1,220.35) and (51.48,222.73) .. (51.48,225.67) .. controls (51.48,228.61) and (49.1,230.99) .. (46.16,230.99) .. controls (43.22,230.99) and (40.84,228.61) .. (40.84,225.67) -- cycle ;
\draw  [color={rgb, 255:red, 140; green, 219; blue, 54 }  ,draw opacity=1 ][dash pattern={on 4.5pt off 4.5pt}] (36.09,227.51) .. controls (36.09,214.48) and (46.65,203.92) .. (59.67,203.92) .. controls (72.7,203.92) and (83.26,214.48) .. (83.26,227.51) .. controls (83.26,240.54) and (72.7,251.1) .. (59.67,251.1) .. controls (46.65,251.1) and (36.09,240.54) .. (36.09,227.51) -- cycle ;
\draw    (69.13,235.51) -- (55.35,250.77) ;
\draw    (69.3,235.44) -- (74.85,245.89) ;
\draw  [fill={rgb, 255:red, 255; green, 255; blue, 255 }  ,fill opacity=1 ] (63.81,235.51) .. controls (63.81,232.57) and (66.19,230.19) .. (69.13,230.19) .. controls (72.07,230.19) and (74.45,232.57) .. (74.45,235.51) .. controls (74.45,238.45) and (72.07,240.83) .. (69.13,240.83) .. controls (66.19,240.83) and (63.81,238.45) .. (63.81,235.51) -- cycle ;

\draw    (259.5,179.74) -- (318.5,179.74) ;
\draw [shift={(320.5,179.74)}, rotate = 180] [color={rgb, 255:red, 0; green, 0; blue, 0 }  ][line width=0.75]    (10.93,-3.29) .. controls (6.95,-1.4) and (3.31,-0.3) .. (0,0) .. controls (3.31,0.3) and (6.95,1.4) .. (10.93,3.29)   ;

\draw (109.45,118.88) node [anchor=north west][inner sep=0.75pt]   [align=left] {$\displaystyle v$};
\draw (137.57,182.75) node [anchor=north west][inner sep=0.75pt]   [align=left] {$\displaystyle w$};
\draw (110.1,239.69) node [anchor=north west][inner sep=0.75pt]   [align=left] {$\displaystyle u$};
\draw (104.99,24.86) node [anchor=north west][inner sep=0.75pt]   [align=left] {$\displaystyle H_{v}$};
\draw (219.29,134.89) node [anchor=north west][inner sep=0.75pt]   [align=left] {$\displaystyle H_{w}$};
\draw (22,184.83) node [anchor=north west][inner sep=0.75pt]   [align=left] {$\displaystyle H_{u}$};

\end{tikzpicture}
\]
\caption{}\label{map}
\end{figure}

 \begin{definition}\label{def: graphical category} 
The graphical category $\mathbf{U}$ is the category whose objects are connected Feynman graphs (Definition~\ref{def: graphs}). The morphisms are the graphical maps from Definition~\ref{def: graphical maps}. 
\end{definition} 
 
As we mentioned at the start of this section, the definition of graphical map we have given is a bit cumbersome. Luckily, one can show that all graphical maps can be described (up to isomorphism) as the composite of three elementary classes of graphical maps: inner coface maps, outer coface maps and codegeneracies (Theorem 2.7 \cite{hry1}). We have already seen examples of outer coface maps: outer coface maps are embeddings.

\begin{definition}\label{def: outer coface} 
An \emph{outer coface map} is either an \emph{embedding}  $d_e: G\rightarrow G'$ in which $G'$ has precisely \emph{one} more internal edge than $G$ or an \emph{edge inclusion} $\updownarrow \rightarrow \medwhitestar_n$. 
\end{definition} 

\begin{example} 
Below we have depicted the outer coface map $d_{e}: G\rightarrow G'$ where $G'$ has the additional inner edge $e=[a,b]$ (highlighted in red). On half edges the graphical map is $(d_e)_0(a^{\dagger})=b$, $(d_e)_0(b^{\dagger})=a$ and is the identity elsewhere. An outer coface map will not change any of the vertices of $G$, which is explicitly written as $(d_e)_1(v)=\medwhitestar_v$ and $(d_e)_1(w)=\medwhitestar_w$. 

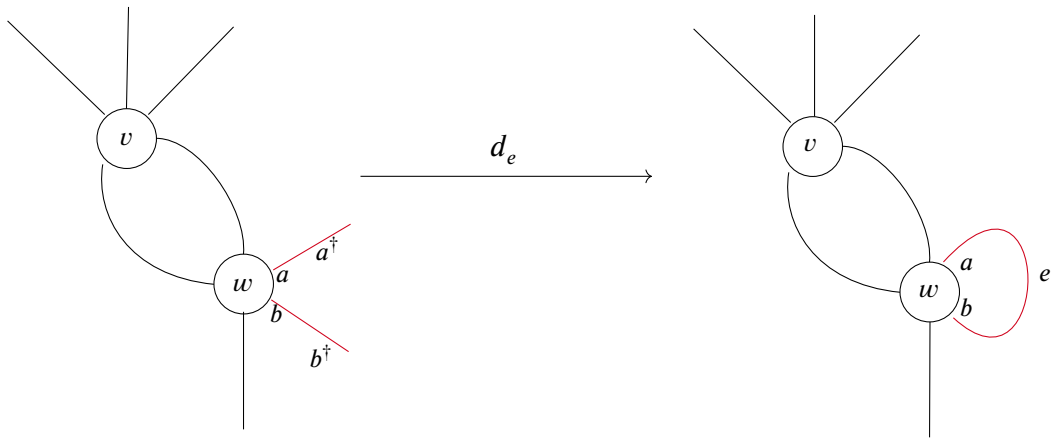
\[\begin{tikzpicture}[x=0.75pt,y=0.75pt,yscale=-1,xscale=1]

\draw  [fill={rgb, 255:red, 255; green, 255; blue, 255 }  ,fill opacity=1 ] (498.97,110.07) .. controls (507.25,110) and (514.02,116.66) .. (514.09,124.94) .. controls (514.16,133.22) and (507.5,140) .. (499.22,140.07) .. controls (490.93,140.14) and (484.16,133.48) .. (484.09,125.19) .. controls (484.02,116.91) and (490.68,110.14) .. (498.97,110.07) -- cycle ;
\draw    (439,66) -- (469.15,94.92) -- (488,113) ;
\draw    (510,113) -- (553,69) ;
\draw  [fill={rgb, 255:red, 255; green, 255; blue, 255 }  ,fill opacity=1 ] (557.97,183.07) .. controls (566.25,183) and (573.02,189.66) .. (573.09,197.94) .. controls (573.16,206.22) and (566.5,213) .. (558.22,213.07) .. controls (549.93,213.14) and (543.16,206.48) .. (543.09,198.19) .. controls (543.02,189.91) and (549.68,183.14) .. (557.97,183.07) -- cycle ;
\draw    (558.22,213.07) -- (558,271) ;
\draw [color={rgb, 255:red, 208; green, 2; blue, 27 }  ,draw opacity=1 ]   (565,183) .. controls (622,121) and (620,257) .. (570,211) ;
\draw    (487,138) .. controls (483,164) and (500,194) .. (543.09,198.19) ;
\draw    (514.09,124.94) .. controls (535,125) and (560,160) .. (557.97,183.07) ;
\draw    (500,59) -- (500,110) ;
\draw  [fill={rgb, 255:red, 255; green, 255; blue, 255 }  ,fill opacity=1 ] (152.97,106.07) .. controls (161.25,106) and (168.02,112.66) .. (168.09,120.94) .. controls (168.16,129.22) and (161.5,136) .. (153.22,136.07) .. controls (144.93,136.14) and (138.16,129.48) .. (138.09,121.19) .. controls (138.02,112.91) and (144.68,106.14) .. (152.97,106.07) -- cycle ;
\draw    (93,62) -- (123.15,90.92) -- (142,109) ;
\draw    (164,109) -- (207,65) ;
\draw  [fill={rgb, 255:red, 255; green, 255; blue, 255 }  ,fill opacity=1 ] (211.97,179.07) .. controls (220.25,179) and (227.02,185.66) .. (227.09,193.94) .. controls (227.16,202.22) and (220.5,209) .. (212.22,209.07) .. controls (203.93,209.14) and (197.16,202.48) .. (197.09,194.19) .. controls (197.02,185.91) and (203.68,179.14) .. (211.97,179.07) -- cycle ;
\draw   (212,208) -- (212,267) ;
\draw    (141,134) .. controls (137,160) and (154,190) .. (197.09,194.19) ;
\draw    (168.09,120.94) .. controls (189,121) and (214,156) .. (211.97,179.07) ;
\draw    (154,55) -- (152.97,106.07) ;
\draw [color={rgb, 255:red, 208; green, 2; blue, 27 }  ,draw opacity=1 ]   (226,202) -- (265,228) ;
\draw [color={rgb, 255:red, 208; green, 2; blue, 27 }  ,draw opacity=1 ]   (227,187) -- (266,164) ;
\draw [->]   (271,140) -- (418,140) ;

\draw (335,118) node [anchor=north west][inner sep=0.75pt]    {$d_{e}$};
\draw (612,185.4) node [anchor=north west][inner sep=0.75pt] [font=\footnotesize]    {$e$};
\draw (148,117) node [anchor=north west][inner sep=0.75pt]  [font=\footnotesize]   {$v$};
\draw (493,120) node [anchor=north west][inner sep=0.75pt]  [font=\footnotesize]   {$v$};
\draw (205,190) node [anchor=north west][inner sep=0.75pt] [font=\footnotesize]    {$w$};
\draw (551,193) node [anchor=north west][inner sep=0.75pt]  [font=\footnotesize]   {$w$};
\draw (247,169) node [anchor=north west][inner sep=0.75pt] [font=\footnotesize]    {$a^{\dagger }$};
\draw (227,186) node [anchor=north west][inner sep=0.75pt]  [font=\footnotesize]   {$a$};
\draw (224,203) node [anchor=north west][inner sep=0.75pt]  [font=\footnotesize]   {$b$};
\draw (244,223) node [anchor=north west][inner sep=0.75pt][font=\footnotesize]     {$b^{\dagger }$};
\draw (572,180) node [anchor=north west][inner sep=0.75pt]  [font=\footnotesize]   {$a$};
\draw (572,200) node [anchor=north west][inner sep=0.75pt] [font=\footnotesize]    {$b$};

\end{tikzpicture}\]
\end{example} 

\begin{definition}\label{def: inner coface}
An \emph{inner coface map}  $d_v: G\rightarrow G'$ is a graphical map defined by blowing-up a single vertex $v$ in $G$ by a graph which has precisely \emph{one} internal edge. 
\end{definition} 

\begin{example} 
The picture below depicts two possible inner coface maps defined at the vertex $v$. 
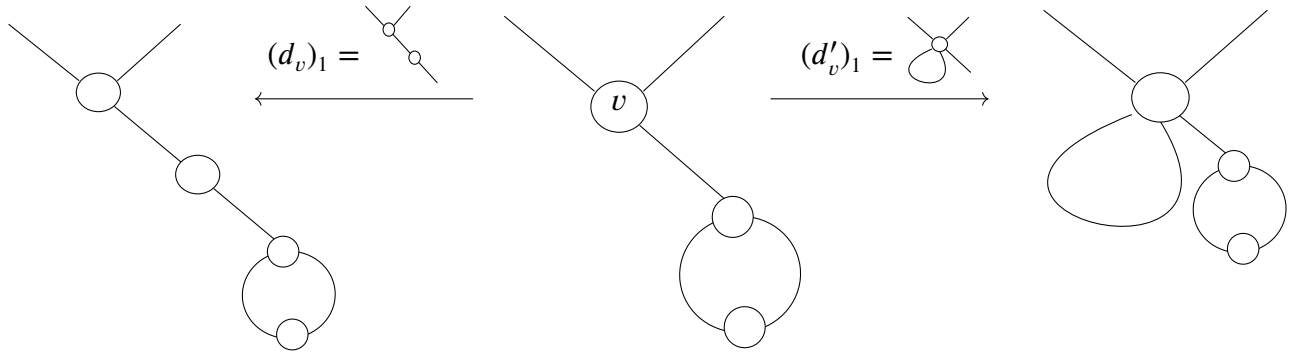
\[
\begin{tikzpicture}[x=0.75pt,y=0.75pt,yscale=-1,xscale=1]

\draw  [fill={rgb, 255:red, 255; green, 255; blue, 255 }  ,fill opacity=1 ] (586.52,116.92) .. controls (594.52,116.86) and (601.05,122.36) .. (601.12,129.19) .. controls (601.19,136.02) and (594.76,141.6) .. (586.77,141.66) .. controls (578.78,141.72) and (572.24,136.23) .. (572.17,129.4) .. controls (572.11,122.57) and (578.53,116.98) .. (586.52,116.92) -- cycle ;
\draw    (527.71,85.52) -- (574.02,121.81) ;
\draw    (599.1,121.81) -- (640.59,85.52) ;
\draw    (596.17,138.31) -- (623,160) ;
\draw    (572.27,137.88) .. controls (447.35,186.57) and (643.19,233.27) .. (586.77,141.66) ;
\draw  [fill={rgb, 255:red, 255; green, 255; blue, 255 }  ,fill opacity=1 ] (313.51,121.11) .. controls (321.37,121.05) and (327.79,126.74) .. (327.86,133.83) .. controls (327.93,140.92) and (321.61,146.72) .. (313.75,146.78) .. controls (305.9,146.84) and (299.48,141.14) .. (299.41,134.05) .. controls (299.34,126.96) and (305.66,121.17) .. (313.51,121.11) -- cycle ;
\draw    (255.71,88.52) -- (301.22,126.18) ;
\draw    (325.87,126.18) -- (366.65,88.52) ;
\draw    (323.98,143.3) -- (370,183) ;
\draw  [fill={rgb, 255:red, 255; green, 255; blue, 255 }  ,fill opacity=1 ] (50.92,117.09) .. controls (57.07,117.04) and (62.09,121.34) .. (62.14,126.68) .. controls (62.2,132.03) and (57.26,136.4) .. (51.11,136.44) .. controls (44.97,136.49) and (39.94,132.19) .. (39.89,126.84) .. controls (39.84,121.5) and (44.78,117.13) .. (50.92,117.09) -- cycle ;
\draw    (5.71,92.52) -- (41.31,120.91) ;
\draw    (60.59,120.91) -- (92.49,92.52) ;
\draw    (59.11,133.82) -- (95.11,163.75) ;
\draw  [fill={rgb, 255:red, 255; green, 255; blue, 255 }  ,fill opacity=1 ] (101.01,158.34) .. controls (107.16,158.29) and (112.18,162.59) .. (112.23,167.94) .. controls (112.29,173.28) and (107.35,177.65) .. (101.2,177.7) .. controls (95.06,177.74) and (90.03,173.44) .. (89.98,168.1) .. controls (89.93,162.75) and (94.87,158.38) .. (101.01,158.34) -- cycle ;
\draw    (109.2,175.07) -- (145.2,205) ;
\draw [->]   (390,130) -- (500,130) ;
\draw  [fill={rgb, 255:red, 255; green, 255; blue, 255 }  ,fill opacity=1 ] (475.28,98.92) .. controls (477.47,98.91) and (479.25,100.55) .. (479.27,102.6) .. controls (479.29,104.64) and (477.53,106.31) .. (475.35,106.33) .. controls (473.16,106.35) and (471.37,104.7) .. (471.36,102.66) .. controls (471.34,100.61) and (473.09,98.94) .. (475.28,98.92) -- cycle ;
\draw    (459.19,89) -- (471.86,100) ;
\draw    (478.72,100) -- (490.07,89) ;
\draw    (478.19,105) -- (491,116.78) ;
\draw    (471.38,105) .. controls (437.22,119.78) and (490.78,133.76) .. (475.35,106.33) ;
\draw [<-]   (130,130) -- (240,130) ;
\draw  [fill={rgb, 255:red, 255; green, 255; blue, 255 }  ,fill opacity=1 ] (197.47,91.93) .. controls (199.07,91.91) and (200.38,93.38) .. (200.39,95.21) .. controls (200.4,97.04) and (199.12,98.53) .. (197.52,98.55) .. controls (195.92,98.56) and (194.61,97.09) .. (194.6,95.26) .. controls (194.59,93.44) and (195.87,91.94) .. (197.47,91.93) -- cycle ;
\draw    (185.71,83.52) -- (194.97,93.23) ;
\draw    (199.99,93.23) -- (208.29,83.52) ;
\draw  (199.6,97.65) -- (208.97,107.89) ;
\draw  [fill={rgb, 255:red, 255; green, 255; blue, 255 }  ,fill opacity=1 ] (210.5,106.04) .. controls (212.1,106.02) and (213.41,107.49) .. (213.42,109.32) .. controls (213.44,111.15) and (212.15,112.64) .. (210.55,112.66) .. controls (208.95,112.67) and (207.65,111.2) .. (207.63,109.38) .. controls (207.62,107.55) and (208.9,106.05) .. (210.5,106.04) -- cycle ;
\draw    (212.63,111.76) -- (222,122) ;
\draw   (344.41,221.98) .. controls (342.36,206.06) and (354.24,191.4) .. (370.94,189.25) .. controls (387.65,187.1) and (402.86,198.26) .. (404.91,214.19) .. controls (406.96,230.11) and (395.08,244.76) .. (378.37,246.91) .. controls (361.66,249.07) and (346.46,237.9) .. (344.41,221.98) -- cycle ;
\draw  [fill={rgb, 255:red, 255; green, 255; blue, 255 }  ,fill opacity=1 ] (369.85,179.02) .. controls (375.52,178.4) and (380.61,182.48) .. (381.22,188.13) .. controls (381.82,193.79) and (377.71,198.87) .. (372.04,199.49) .. controls (366.37,200.11) and (361.28,196.03) .. (360.67,190.37) .. controls (360.06,184.72) and (364.17,179.64) .. (369.85,179.02) -- cycle ;
\draw  [fill={rgb, 255:red, 255; green, 255; blue, 255 }  ,fill opacity=1 ] (375.89,234.38) .. controls (381.56,233.76) and (386.65,237.84) .. (387.26,243.5) .. controls (387.87,249.15) and (383.76,254.23) .. (378.08,254.85) .. controls (372.41,255.47) and (367.32,251.39) .. (366.71,245.74) .. controls (366.11,240.08) and (370.21,235) .. (375.89,234.38) -- cycle ;
\draw   (603.36,188.27) .. controls (601.78,176.31) and (610.9,165.31) .. (623.73,163.69) .. controls (636.57,162.07) and (648.24,170.46) .. (649.82,182.41) .. controls (651.4,194.37) and (642.27,205.38) .. (629.44,206.99) .. controls (616.61,208.61) and (604.93,200.23) .. (603.36,188.27) -- cycle ;
\draw  [fill={rgb, 255:red, 255; green, 255; blue, 255 }  ,fill opacity=1 ] (622.89,156) .. controls (627.25,155.54) and (631.16,158.6) .. (631.63,162.85) .. controls (632.09,167.09) and (628.94,170.91) .. (624.58,171.38) .. controls (620.22,171.84) and (616.31,168.77) .. (615.84,164.53) .. controls (615.38,160.28) and (618.53,156.47) .. (622.89,156) -- cycle ;
\draw  [fill={rgb, 255:red, 255; green, 255; blue, 255 }  ,fill opacity=1 ] (627.53,197.58) .. controls (631.89,197.12) and (635.8,200.18) .. (636.27,204.43) .. controls (636.73,208.67) and (633.58,212.49) .. (629.22,212.95) .. controls (624.86,213.42) and (620.95,210.35) .. (620.49,206.11) .. controls (620.02,201.86) and (623.18,198.04) .. (627.53,197.58) -- cycle ;
\draw   (123.8,231.27) .. controls (122.23,219.31) and (131.35,208.31) .. (144.18,206.69) .. controls (157.01,205.07) and (168.69,213.46) .. (170.27,225.41) .. controls (171.85,237.37) and (162.72,248.38) .. (149.89,249.99) .. controls (137.06,251.61) and (125.38,243.23) .. (123.8,231.27) -- cycle ;
\draw  [fill={rgb, 255:red, 255; green, 255; blue, 255 }  ,fill opacity=1 ] (143.34,199) .. controls (147.7,198.54) and (151.61,201.6) .. (152.08,205.85) .. controls (152.54,210.09) and (149.39,213.91) .. (145.03,214.38) .. controls (140.67,214.84) and (136.76,211.77) .. (136.29,207.53) .. controls (135.83,203.28) and (138.98,199.47) .. (143.34,199) -- cycle ;
\draw  [fill={rgb, 255:red, 255; green, 255; blue, 255 }  ,fill opacity=1 ] (147.98,240.58) .. controls (152.34,240.12) and (156.25,243.18) .. (156.72,247.43) .. controls (157.18,251.67) and (154.03,255.49) .. (149.67,255.95) .. controls (145.31,256.42) and (141.4,253.35) .. (140.94,249.11) .. controls (140.47,244.86) and (143.63,241.04) .. (147.98,240.58) -- cycle ;

\draw (308,127) node [anchor=north west][inner sep=0.75pt]    {$v$};
\draw (404,100) node [anchor=north west][inner sep=0.75pt]    {$( d'_{v})_{1} =$};
\draw (135,100) node [anchor=north west][inner sep=0.75pt]    {$( d_{v})_{1} =$};

\end{tikzpicture}
\]
\end{example} 
  
Codegeneracy maps "delete" arity $2$ vertices (a vertex $v$ with $|\nb(v)|=2$).  

\begin{definition}\label{def: deg}
A \emph{codegeneracy map}  $s_v: G\rightarrow G'$ is a graphical map defined by ``blowing-up'' a vertex $v$ in $G$ by $\updownarrow$. 
\end{definition} 

\begin{example} 
In the graphical map below, the vertex $v$ has arity two. The codegeneracy map $s_v: G\rightarrow G'$ is the identity on half edges and $(s_w)_1=\medwhitestar_w$ at all vertices except $w=v$ where $(s_v)_1=\updownarrow$. 
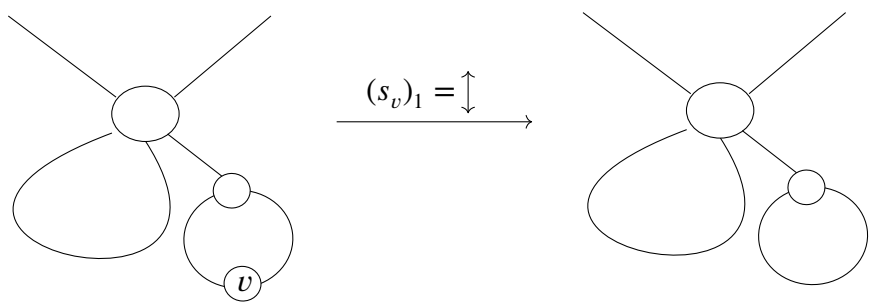
\[
\begin{tikzpicture}[x=0.75pt,y=0.75pt,yscale=-.85,xscale=.85]

\draw  [fill={rgb, 255:red, 255; green, 255; blue, 255 }  ,fill opacity=1 ] (165.36,118.97) .. controls (176.4,118.9) and (185.42,126.14) .. (185.52,135.16) .. controls (185.61,144.18) and (176.74,151.55) .. (165.7,151.63) .. controls (154.67,151.71) and (145.64,144.46) .. (145.55,135.44) .. controls (145.46,126.42) and (154.33,119.05) .. (165.36,118.97) -- cycle ;
\draw    (84.16,77.52) -- (148.1,125.43) ;
\draw    (182.73,125.43) -- (240.01,77.52) ;
\draw    (178.68,147.21) -- (215.72,175.84) ;
\draw    (145.68,146.63) .. controls (-26.79,210.92) and (243.61,272.56) .. (165.7,151.63) ;
\draw   (188.6,213.15) .. controls (186.42,197.37) and (199.02,182.84) .. (216.74,180.7) .. controls (234.45,178.57) and (250.58,189.64) .. (252.76,205.42) .. controls (254.93,221.21) and (242.34,235.74) .. (224.62,237.87) .. controls (206.9,240.01) and (190.78,228.94) .. (188.6,213.15) -- cycle ;
\draw  [fill={rgb, 255:red, 255; green, 255; blue, 255 }  ,fill opacity=1 ] (215.58,170.56) .. controls (221.59,169.95) and (226.99,173.99) .. (227.64,179.6) .. controls (228.28,185.2) and (223.92,190.24) .. (217.91,190.85) .. controls (211.89,191.47) and (206.49,187.42) .. (205.84,181.82) .. controls (205.2,176.21) and (209.56,171.17) .. (215.58,170.56) -- cycle ;
\draw  [fill={rgb, 255:red, 255; green, 255; blue, 255 }  ,fill opacity=1 ] (223.5,225.38) .. controls (229.54,225.67) and (234.28,230.47) .. (234.09,236.11) .. controls (233.89,241.75) and (228.84,246.09) .. (222.8,245.8) .. controls (216.76,245.51) and (212.01,240.71) .. (212.21,235.07) .. controls (212.4,229.44) and (217.46,225.1) .. (223.5,225.38) -- cycle ;
\draw  [fill={rgb, 255:red, 255; green, 255; blue, 255 }  ,fill opacity=1 ] (506.36,116.97) .. controls (517.4,116.9) and (526.42,124.14) .. (526.52,133.16) .. controls (526.61,142.18) and (517.74,149.55) .. (506.7,149.63) .. controls (495.67,149.71) and (486.64,142.46) .. (486.55,133.44) .. controls (486.46,124.42) and (495.33,117.05) .. (506.36,116.97) -- cycle ;
\draw    (425.16,75.52) -- (489.1,123.43) ;
\draw    (523.73,123.43) -- (581.01,75.52) ;
\draw    (519.68,145.21) -- (556.72,173.84) ;
\draw    (486.68,144.63) .. controls (314.21,208.92) and (584.61,270.56) .. (506.7,149.63) ;
\draw   (529.6,211.15) .. controls (527.42,195.37) and (540.02,180.84) .. (557.74,178.7) .. controls (575.45,176.57) and (591.58,187.64) .. (593.76,203.42) .. controls (595.93,219.21) and (583.34,233.74) .. (565.62,235.87) .. controls (547.9,238.01) and (531.78,226.94) .. (529.6,211.15) -- cycle ;
\draw  [fill={rgb, 255:red, 255; green, 255; blue, 255 }  ,fill opacity=1 ] (556.58,168.56) .. controls (562.59,167.95) and (567.99,171.99) .. (568.64,177.6) .. controls (569.28,183.2) and (564.92,188.24) .. (558.91,188.85) .. controls (552.89,189.47) and (547.49,185.42) .. (546.84,179.82) .. controls (546.2,174.21) and (550.56,169.17) .. (556.58,168.56) -- cycle ;
\draw [->]   (279,140) -- (395,140) ;
\draw [<->]   (357,110) -- (357,135) ;

\draw (122,7.4) node [anchor=north west][inner sep=0.75pt]    {$$};
\draw (295,112.4) node [anchor=north west][inner sep=0.75pt]    {$( s_{v})_{1} =$};
\draw (218,230) node [anchor=north west][inner sep=0.75pt]    {$v$};

\end{tikzpicture}\]
\end{example}

\begin{remark}
In practice, when working with graphical maps it is often enough to study these elementary maps (and isomorphisms).   We also note that there is an equivalent, purely combinatorial, definition of graphical maps which does not require reference to graph substitution (see Theorem A~\cite{phil_graphs}). 
\end{remark} 

\subsection{Modular dendroidal sets and the nerve theorem}\label{sec: nerve theorem}
As we have (hopefully) motivated with pictures, the graphical category $\mathbf{U}$ is closely related to modular operads. In particular, every object of $\mathbf{U}$ freely generates a modular operad (Definition 2.7 \cite{hry2}). 

\begin{definition}\label{def:modular_op_from_graph}
The modular operad $\left<G\right>$ generated by a graph $G$ is the free modular operad whose: 
\begin{itemize} 
\item set of colours is the set of half edges $A$; 
\item a collection of $\Sigma_n$-sets is $E (a_1,\ldots, a_n) = \begin{cases} \{v\} \ \text{if} \ (a_1,\ldots, a_n)=\partial(\medwhitestar_v) \\
\emptyset \ \text{otherwise}. \end{cases}$
\item $\left<G\right> =F(E)$
\end{itemize} 
\end{definition}

\begin{example} 
To describe the modular operad $\left<G\right>$ generated by the graph $G$ in Figure~\ref{ex:graph_generating_modular} we recall that $\partial(\medwhitestar_v) = (a_1^{\dagger},a_2^{\dagger}, c_2, c_1)$ and $\partial(\medwhitestar_w) = (b_1^{\dagger},b_2^{\dagger},b_3^{\dagger}, c_4, c_3)$. The underlying collection of $\left<G\right>$ consists of two one-point sets:
\[\begin{array}{cc}
   E(a_1^{\dagger},a_2^{\dagger}, c_2, c_1) =\{v\},   &  
  E (b_1^{\dagger},b_2^{\dagger},b_3^{\dagger}, c_4, c_3)=\{w\}.
\end{array}\] The graph $G$ provides gluing instructions.

\begin{figure}[h!]

\begin{tikzpicture}[x=0.75pt,y=0.75pt,yscale=-.75,xscale=.75]

\draw  [fill={rgb, 255:red, 255; green, 255; blue, 255 }  ,fill opacity=1 ] (172.97,112.07) .. controls (181.25,112) and (188.02,118.66) .. (188.09,126.94) .. controls (188.16,135.22) and (181.5,142) .. (173.22,142.07) .. controls (164.93,142.14) and (158.16,135.48) .. (158.09,127.19) .. controls (158.02,118.91) and (164.68,112.14) .. (172.97,112.07) -- cycle ;
\draw    (113,68) -- (143.15,96.92) -- (162,115) ;
\draw    (184,115) -- (227,71) ;
\draw  [fill={rgb, 255:red, 255; green, 255; blue, 255 }  ,fill opacity=1 ] (231.97,185.07) .. controls (240.25,185) and (247.02,191.66) .. (247.09,199.94) .. controls (247.16,208.22) and (240.5,215) .. (232.22,215.07) .. controls (223.93,215.14) and (217.16,208.48) .. (217.09,200.19) .. controls (217.02,191.91) and (223.68,185.14) .. (231.97,185.07) -- cycle ;
\draw    (232.22,215.07) -- (232,273) ;
\draw    (173.22,142.07) .. controls (177,158) and (183,192) .. (217.09,200.19) ;
\draw    (188.09,126.94) .. controls (209,127) and (234,162) .. (231.97,185.07) ;
\draw [color={rgb, 255:red, 0; green, 0; blue, 0 }  ,draw opacity=1 ]   (246,208) -- (285,234) ;
\draw [color={rgb, 255:red, 0; green, 0; blue, 0 }  ,draw opacity=1 ]   (247,193) -- (286,170) ;

\draw (165,120) node [anchor=north west][inner sep=0.75pt] [font=\footnotesize]    {$v$};
\draw (225,195) node [anchor=north west][inner sep=0.75pt] [font=\footnotesize]    {$w$};
\draw (105,75) node [anchor=north west][inner sep=0.75pt]   [font=\footnotesize]  {$a_{1}^{\dagger }$};
\draw (135,103) node [anchor=north west][inner sep=0.75pt] [font=\footnotesize]    {$a_{1}$};
\draw (244,215.4) node [anchor=north west][inner sep=0.75pt] [font=\footnotesize]    {$b_{2}$};
\draw (264,235.4) node [anchor=north west][inner sep=0.75pt]   [font=\footnotesize]  {$b_{2}^{\dagger }$};
\draw (214,253.4) node [anchor=north west][inner sep=0.75pt]   [font=\footnotesize]  {$b_{1}^{\dagger }$};
\draw (215,220.4) node [anchor=north west][inner sep=0.75pt] [font=\footnotesize]    {$b_{1}$};
\draw (215,75) node [anchor=north west][inner sep=0.75pt]   [font=\footnotesize]  {$a_{2}^{\dagger }$};
\draw (195,103) node [anchor=north west][inner sep=0.75pt][font=\footnotesize]     {$a_{2}$};
\draw (252,189.4) node [anchor=north west][inner sep=0.75pt] [font=\footnotesize]    {$b_{3}$};
\draw (277,173.4) node [anchor=north west][inner sep=0.75pt]   [font=\footnotesize]  {$b_{3}^{\dagger }$};
\draw (188,193.4) node [anchor=north west][inner sep=0.75pt][font=\footnotesize]     {$c_{1}$};
\draw (210,165) node [anchor=north west][inner sep=0.75pt] [font=\footnotesize]    {$c_{2}$};
\draw (185,135) node [anchor=north west][inner sep=0.75pt]  [font=\footnotesize]   {$c_{4}$};
\draw (155,145) node [anchor=north west][inner sep=0.75pt] [font=\footnotesize]    {$c_{3}$};
\end{tikzpicture}

\caption{}\label{ex:graph_generating_modular}
\end{figure}
\end{example}

\begin{prop}[Proposition 2.25 \cite{hry2}] The assignment $G\mapsto \left<G\right>$ defines a faithful functor $\mathbf{U}\rightarrow \ModOp$ which is injective on isomorphism classes of objects. 
\end{prop} 

\begin{exercise}
Write $\ModOp_{\fC}$ for the category of modular operads with $\fC$-colours. Show that $\ModOp_{\emptyset}$ is equivalent to the category of sets and that $\left<\medwhitestar_0\right>$ is a generator. Here $G=\medwhitestar_0$ is the modular operad freely generated by an isolated vertex. See Example 2.20~\cite{hry2} for a hint. 
\end{exercise}

We note that the functor $$J:\mathbf{U} \longrightarrow \ModOp$$ is not full. To see this one can consider the graphs $G$ and $G'$ in Figure~\ref{ex:not_full}. There is a map of modular operads from $\left<G\right>$ to $\left<G'\right>$ which sends each $v_i$ to $v$ and each $w_j$ to $w$ but there is no graphical map $G\rightarrow G'$ which has this behavior. 

\begin{figure}[h!]
\[\begin{tikzpicture}[x=0.75pt,y=0.75pt,yscale=-.75,xscale=.75]

\draw  [fill={rgb, 255:red, 255; green, 255; blue, 255 }  ,fill opacity=1 ] (79.84,89.98) .. controls (88.85,90.25) and (96.04,98.87) .. (95.9,109.22) .. controls (95.75,119.58) and (88.33,127.76) .. (79.32,127.49) .. controls (70.31,127.22) and (63.12,118.6) .. (63.26,108.24) .. controls (63.41,97.89) and (70.83,89.71) .. (79.84,89.98) -- cycle ;
\draw  [fill={rgb, 255:red, 255; green, 255; blue, 255 }  ,fill opacity=1 ] (78.68,198.72) .. controls (87.57,198.98) and (94.66,208.33) .. (94.5,219.58) .. controls (94.34,230.84) and (87,239.75) .. (78.11,239.48) .. controls (69.21,239.22) and (62.13,229.87) .. (62.28,218.62) .. controls (62.44,207.36) and (69.78,198.45) .. (78.68,198.72) -- cycle ;
\draw  [fill={rgb, 255:red, 255; green, 255; blue, 255 }  ,fill opacity=1 ] (187.51,97.53) .. controls (196.06,97.79) and (202.9,105.36) .. (202.77,114.44) .. controls (202.64,123.52) and (195.61,130.67) .. (187.05,130.42) .. controls (178.49,130.16) and (171.66,122.59) .. (171.79,113.51) .. controls (171.92,104.43) and (178.95,97.28) .. (187.51,97.53) -- cycle ;
\draw  [fill={rgb, 255:red, 255; green, 255; blue, 255 }  ,fill opacity=1 ] (186.43,192) .. controls (195.33,192.27) and (202.43,201.33) .. (202.28,212.24) .. controls (202.12,223.15) and (194.78,231.78) .. (185.88,231.52) .. controls (176.97,231.25) and (169.88,222.19) .. (170.03,211.27) .. controls (170.18,200.36) and (177.53,191.73) .. (186.43,192) -- cycle ;
\draw    (62.28,218.62) .. controls (19.84,216.04) and (18.79,112.01) .. (63.26,108.24) ;
\draw    (171,205) .. controls (150.28,182.36) and (117.69,128.3) .. (95.88,110.56) ;
\draw    (94.51,219.29) .. controls (111.21,205.1) and (106.68,206.08) .. (132.58,176.55) ;
\draw    (203.16,106.64) .. controls (246.98,94.85) and (249.78,220.85) .. (202.51,217.9) ;
\draw    (146.18,158.83) .. controls (161.07,141.1) and (155.89,148) .. (175.96,123.38) ;
\draw   (482,164.8) .. controls (482,132.96) and (508.86,107.15) .. (542,107.15) .. controls (575.14,107.15) and (602,132.96) .. (602,164.8) .. controls (602,196.64) and (575.14,222.45) .. (542,222.45) .. controls (508.86,222.45) and (482,196.64) .. (482,164.8) -- cycle ;
\draw  [fill={rgb, 255:red, 255; green, 255; blue, 255 }  ,fill opacity=1 ] (542.43,86.74) .. controls (553.66,86.96) and (562.56,96.27) .. (562.33,107.55) .. controls (562.09,118.82) and (552.79,127.78) .. (541.57,127.56) .. controls (530.34,127.34) and (521.44,118.02) .. (521.67,106.75) .. controls (521.91,95.47) and (531.21,86.51) .. (542.43,86.74) -- cycle ;
\draw  [fill={rgb, 255:red, 255; green, 255; blue, 255 }  ,fill opacity=1 ] (540.31,197.17) .. controls (551.53,197.4) and (560.44,206.71) .. (560.2,217.99) .. controls (559.96,229.26) and (550.67,238.22) .. (539.44,238) .. controls (528.22,237.78) and (519.31,228.46) .. (519.55,217.18) .. controls (519.79,205.91) and (529.08,196.95) .. (540.31,197.17) -- cycle ;

\draw (69,103) node [anchor=north west][inner sep=0.75pt][font=\footnotesize]     {$v_{1}$};
\draw (69,209) node [anchor=north west][inner sep=0.75pt] [font=\footnotesize]    {$w_{1}$};
\draw (176,105) node [anchor=north west][inner sep=0.75pt] [font=\footnotesize]    {$v_{2}$};
\draw (175,203) node [anchor=north west][inner sep=0.75pt][font=\footnotesize]     {$w_{2}$};
\draw (535,100) node [anchor=north west][inner sep=0.75pt] [font=\footnotesize]    {$v$};
\draw (533,210) node [anchor=north west][inner sep=0.75pt] [font=\footnotesize]    {$w$};
\end{tikzpicture}\]
    \caption{}
    \label{ex:not_full}
\end{figure}

\subsubsection{Modular dendroidal  sets} 
Using the functor $\mathbf{U}\rightarrow\ModOp$ we can define the \emph{nerve} of a modular operad: \[\begin{tikzcd}\ModOp \arrow[r, "N"] &
\mathbf{Set}^{\bU^{op}}.
\end{tikzcd}\] The category of presheaves $\mathbf{Set}^{\bU^{op}}$ is called the category of \emph{modular dendroidal sets}.\footnote{In \cite{hry1},\cite{hry2}, we avoid naming the category of set-valued $\mathbf{U}$-presheaves largely because we used the term ``graphical sets'' in \cite{hry15} when modeling $\infty$-properads. The name modular dendroidal sets, suggested to us by Ieke Moerdijk, follows the convention in the literature by using the term ``modular'' or ``cyclic'' modify the term ``operad''.} Our goal for the remainder of this first lecture is to define a subcategory of modular dendroidal  sets which satisfy a strict Segal condition. 

\begin{definition}\label{def:modular_dendroidal_sets}
The category of \emph{modular dendroidal sets} $\mathbf{Set}^{\bU^{op}}$ is the category whose objects are functors $X: \mathbf{U}^{op}\rightarrow \textbf{Set}$. Morphisms in $\mathbf{Set}^{\bU^{op}}$ are natural transformations.
\end{definition} 

In these notes we follow some notation conventions from the theory of dendroidal sets (\cite{Moerdijk_notes}, \cite{MW_Kan}, \cite{MW1}). Given an $X\in\mathbf{Set}^{\bU^{op}}$ we will write $X_G$ for the evaluation of the presheaf $X$ at a graph $G\in\mathbf{U}$. Similarly, for every morphism $\varphi: G\rightarrow G'$ in $\bU$, there's map $\varphi^*: X_{G'}\rightarrow X_{G}$ in $\mathbf{Set}^{\bU^{op}}$.  

There are several examples of modular dendroidal sets which warrant special notation.  For any graph $G$, the \emph{representable} presheaf \[\mathbf{U}[G]:=\mathbf{U}(-;G)\] is given by  \[\mathbf{U}[G]_{H}:=\mathbf{U}(H,G),\] where $H$ ranges over all graphs $H\in \bU$.  We think of an element $x\in X_{G}$ as a \emph{decoration} of shape $G$, similar to those we depicted in Figure~\ref{operations modular}. The Yoneda Lemma tells us that ``the set of all decorations of $G$'' can be identified with the set of maps out of the representable $\mathbf{U}[G]$ as there is a bijection \[X_G = \textbf{Set}^{\bU^{op}}(\bU[G], X).\] 

\medskip
 
For any graph $G$ with at least two vertices, we define a presheaf $X^1_{G}$ which captures the ``local decoration data'' of $G$. 
If $G$ is a graph with at least two vertices, each internal edge between vertices $v$ and $w$ produces a diagram of embeddings:
\[\begin{tikzcd} \medwhitestar_{v}\ar[dr] & \arrow[l] \updownarrow \ar[d] \arrow[r]& \medwhitestar_{w} \ar[dl]\\
&G\end{tikzcd}\] 
in $\bU$. The maps in this diagram consist of edge and star inclusions from \eqref{star_inclusion}. 

\begin{definition}\label{def: ribbon}
Let $X$ be a modular dendroidal  set and $G$ be a graph with at least $2$ vertices. The \emph{corolla ribbon} or \emph{spine} of $X$ at $G$ is defined by: \[ X^1_{G} = \lim\limits_{ \medwhitestar_v\leftarrow \updownarrow\rightarrow \medwhitestar_w} \left( \begin{tikzcd} X_{\medwhitestar_{v}} \arrow[dr] & & \arrow[dl] X_{\medwhitestar_{w}} \\&X_{\updownarrow}. & \end{tikzcd}\right)\\ \] Here, the limit ranges over all edge inclusions and stars of $G$.
\end{definition}

\begin{example}
In the graph $G$ in Figure~\ref{example:ribbon} the internal edges $e_1=[c_3,c_1]$, $e_2=[c_4,c_2]$ and $e_3=[b_2,b_3]$ are highlighted in red. To simplify the edge inclusion diagrams we will write $v$ for $\medwhitestar_v$ and $w$ for $\medwhitestar_w$: 

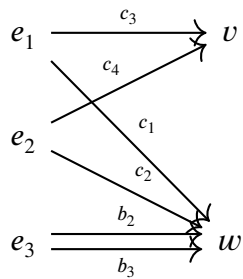
\[\begin{tikzcd}
e_1 \arrow[rr, "c_3"]\arrow[ddrr, "c_1"] & & v\\
e_2 \arrow[urr, "c_4"]\arrow[drr, "c_2"]\\
e_3 \arrow[rr, shift left =.1cm, "b_2"]\arrow[rr, shift right =.1cm, swap, "b_3"]{}{}&& w
\end{tikzcd}\]

The presheaf $X^1_{G}$ is therefore the limit over the diagram: 

\[\begin{tikzcd}
X_{e_1}  & & \arrow[ll, "c_3^*"] X_{v} \arrow[dll, "c_4^*"]\\
X_{e_2} \\
X_{e_3} && \arrow[ll, shift left =.1cm, "b_2^*"]\arrow[ll, shift right =.1cm, swap, "b_3^*"]\arrow[uull, "c_1^*", swap, very near end]X_{w} \arrow[ull, "c_2^*", swap, near end]
\end{tikzcd}\]
\end{example}
\begin{figure}
\[\begin{tikzpicture}[x=0.75pt,y=0.75pt,yscale=-1,xscale=1]

\draw  [fill={rgb, 255:red, 255; green, 255; blue, 255 }  ,fill opacity=1 ] (102.97,90.07) .. controls (111.25,90) and (118.02,96.66) .. (118.09,104.94) .. controls (118.16,113.22) and (111.5,120) .. (103.22,120.07) .. controls (94.93,120.14) and (88.16,113.48) .. (88.09,105.19) .. controls (88.02,96.91) and (94.68,90.14) .. (102.97,90.07) -- cycle ;
\draw    (43,46) -- (73.15,74.92) -- (92,93) ;
\draw    (114,93) -- (157,49) ;
\draw  [fill={rgb, 255:red, 255; green, 255; blue, 255 }  ,fill opacity=1 ] (161.97,163.07) .. controls (170.25,163) and (177.02,169.66) .. (177.09,177.94) .. controls (177.16,186.22) and (170.5,193) .. (162.22,193.07) .. controls (153.93,193.14) and (147.16,186.48) .. (147.09,178.19) .. controls (147.02,169.91) and (153.68,163.14) .. (161.97,163.07) -- cycle ;
\draw    (162.22,193.07) -- (162,251) ;
\draw [color={rgb, 255:red, 208; green, 2; blue, 27 }  ,draw opacity=1 ]   (103.22,120.07) .. controls (107,136) and (113,170) .. (147.09,178.19) ;
\draw [color={rgb, 255:red, 208; green, 2; blue, 27 }  ,draw opacity=1 ]   (118.09,104.94) .. controls (139,105) and (164,140) .. (161.97,163.07) ;
\draw [color={rgb, 255:red, 208; green, 2; blue, 27 }  ,draw opacity=1 ]   (177,169) .. controls (240,79) and (253,302) .. (173,192) ;

\draw (99,100) node [anchor=north west][inner sep=0.75pt]    {$v$};
\draw (155,175) node [anchor=north west][inner sep=0.75pt]    {$w$};
\draw (36,52.4) node [anchor=north west][inner sep=0.75pt]    {$a_{1}^{\dagger }$};
\draw (70,83) node [anchor=north west][inner sep=0.75pt]    {$a_{1}$};
\draw (178,187.4) node [anchor=north west][inner sep=0.75pt]    {$b_{2}$};
\draw (144,231.4) node [anchor=north west][inner sep=0.75pt]    {$b_{1}^{\dagger }$};
\draw (146,198.4) node [anchor=north west][inner sep=0.75pt]    {$b_{1}$};
\draw (144,54.4) node [anchor=north west][inner sep=0.75pt]    {$a_{2}^{\dagger }$};
\draw (122,83) node [anchor=north west][inner sep=0.75pt]    {$a_{2}$};
\draw (179,163) node [anchor=north west][inner sep=0.75pt]    {$b_{3}$};
\draw (120,170) node [anchor=north west][inner sep=0.75pt]    {$c_{1}$};
\draw (144,143.4) node [anchor=north west][inner sep=0.75pt]    {$c_{2}$};
\draw (120,110) node [anchor=north west][inner sep=0.75pt]    {$c_{4}$};
\draw (86,123) node [anchor=north west][inner sep=0.75pt]    {$c_{3}$};
\end{tikzpicture}\]

\caption{}\label{example:ribbon}
\end{figure}
\begin{exercise}\label{exercise: limit is prod}
In the case when $X_{\updownarrow}=\ast$, show that $X^1_{G}=\prod\limits_{v\in V(G)} X_{\medwhitestar_v}$. 
\end{exercise}

\begin{definition}\label{def: Segal map}
For a graph $G$ with at least two vertices the Segal map is the map: 
\[\begin{tikzcd} X_G \arrow[r] & X^1_{G}\subseteq \prod\limits_{v\in V(G)} X_{\medwhitestar_v} \end{tikzcd}\] induced by the embeddings $\medwhitestar_v\hookrightarrow G$. 
\end{definition}

The intuition is that the Segal map says compares the ``decorations of the graph $G$'' with the ``decorations at each vertex''.   

\begin{definition}\label{def:strict_Segal}
A modular dendroidal  set $X\in\mathbf{Set}^{\bU^{op}}$ is strictly \emph{Segal} if the Segal map is a bijection for each $G$ in $\bU$. 
\end{definition}

Modular operads can be identified with the (strictly) Segal modular dendroidal  sets via the following construction. 

\begin{definition} Let $\bP$ be a discrete modular operad and let $G$ be any graph in $\bU$. The modular \emph{nerve} functor \[\begin{tikzcd} N: \ModOp \arrow[r] & \mathbf{Set}^{\mathbf{U}^{op}} \end{tikzcd} \] is defined by \[ N\bP_{G} = \ModOp(\left<G\right>, \bP).\]
\end{definition} 

\begin{exercise} 
Given a graph $G\in\bU$, we now have two ways to assign an object in $\mathbf{Set}^{\mathbf{U}^{op}}$ to $G$: we can take the representable presheaf $\bU[G]$ or we can take the nerve of the free modular operad $\left<G\right>$,  $N\left<G\right>$. The representable $\bU[G]$ is a sub-object of $N\left<G\right>$ (since $J:\bU\rightarrow \ModOp$ is faithful) but they nearly never coincide. 

\begin{enumerate}
    \item Let $G$ be the loop with one node and show $\bU[G]\subset N \left<G\right>$.
    \item Show that we have $\bU[\star_0]= N \left<\star_0\right>$. 
\end{enumerate}

\end{exercise} 

\medskip

For any graph $G$, we picture the set $N\bP_{G}$ as the set of $\bP$-decorations of the graph $G$ (Figure~\ref{operations modular}). In particular, if $\bP$ is a $\mathfrak{C}$-coloured modular operad, then the set $$N\bP_{\updownarrow} :=\ModOp(\left<\updownarrow\right>, \bP) = \mathfrak{C}.$$ 
\begin{exercise}
For any $n$, check that the set 
$$N\bP_{\medwhitestar_n} := \ModOp(\left<\medwhitestar_n\right>, \bP) = \bP(c_1,\ldots c_n).$$ Note that the symmetric group actions on $\bP(c_1,\ldots c_n)$ are captured by the isomorphisms of the graph $\medwhitestar_n$. 
\end{exercise} 

The following theorem says that the strictly Segal modular dendroidal  sets are precisely those which live in the essential image of the nerve functor. We will not give the full proof here, but to provide some intuition for the idea, consider that given a modular operad $\bP$, the fact that $N$ is a functor means that we have morphisms such as \[\begin{tikzcd}  N\bP_{\medwhitestar_{n_1}}\times N\bP_{\medwhitestar_{n_2}} := \bP(c_1,\ldots c_{n_1}) \times\bP(c_1,\ldots c_{n_2}) \arrow[r, "\circ_{ij}"]& \bP(c_1,\ldots,\hat{c_i},\ldots, \hat{c_j},\ldots c_{n_1+n_2})=:N\bP_{\medwhitestar_{G}} \end{tikzcd}\] given by composition operations. Here, the graph $G$ is the graph with a single internal edge given by $e=[c_i,c_j^{\dagger}]$. Similar maps are given for contraction operations. The key to a proof of Theorem~\ref{thm: nerve}, is then showing that compositions of contraction operations of $\bP$ are precisely the composite of the dashed maps in \[\begin{tikzcd} \prod\limits_{v\in V(G)} N\bP_{\medwhitestar_v}\arrow[r, shift left =.2cm, dotted]& N\bP_G \arrow[r] \arrow[l, shift left =.25, "Segal"] & N\bP_{\medwhitestar_G}\end{tikzcd}\] with the natural maps $N\bP_{G}\rightarrow N\bP_{\medwhitestar_G}$ induced by the graphical maps $\medwhitestar_{G}\rightarrow G$. 

\begin{example}
The graphical map $\medwhitestar_{n-2} \to G$ obtained by substituting $G$ for the vertex $v$ is depicted in Figure~\ref{fig: contraction} induces the map map $N\bP_G \to N\bP_{\medwhitestar_{n-2}}$ in the diagram:
\[\begin{tikzcd}
N\bP_{\medwhitestar_n} \arrow[d]\arrow["\xi", r, dashed] & N\bP_{\medwhitestar_{n-2}}\\ 
N\bP_{G}. \arrow[ru]            
\end{tikzcd}\] The dashed map is precisely the application of the nerve functor to the contraction operation in $\bP$.  
\begin{figure}[h!]
\[
\begin{tikzpicture}[x=0.75pt,y=0.75pt,yscale=-0.7,xscale=0.7]

\draw   (170.53,38.51) .. controls (170.53,34.69) and (173.63,31.59) .. (177.45,31.59) .. controls (181.28,31.59) and (184.37,34.69) .. (184.37,38.51) .. controls (184.37,42.33) and (181.28,45.43) .. (177.45,45.43) .. controls (173.63,45.43) and (170.53,42.33) .. (170.53,38.51) -- cycle ;
\draw    (181.43,32.87) .. controls (206.2,-1.37) and (224.47,54.29) .. (183.98,40.94) ;
\draw   (170.53,38.51) .. controls (170.53,34.69) and (173.63,31.59) .. (177.45,31.59) .. controls (181.28,31.59) and (184.37,34.69) .. (184.37,38.51) .. controls (184.37,42.33) and (181.28,45.43) .. (177.45,45.43) .. controls (173.63,45.43) and (170.53,42.33) .. (170.53,38.51) -- cycle ;
\draw    (181.43,44.34) -- (191.63,72.1) ;
\draw    (154.84,67.28) -- (173.2,43.61) ;
\draw    (142.22,38.57) -- (170.53,38.51) ;
\draw    (157.57,10.07) -- (174.78,32.14) ;
\draw   (170.53,38.51) .. controls (170.53,34.69) and (173.63,31.59) .. (177.45,31.59) .. controls (181.28,31.59) and (184.37,34.69) .. (184.37,38.51) .. controls (184.37,42.33) and (181.28,45.43) .. (177.45,45.43) .. controls (173.63,45.43) and (170.53,42.33) .. (170.53,38.51) -- cycle ;
\draw  [color={rgb, 255:red, 4; green, 146; blue, 194 }  ,draw opacity=1 ][dash pattern={on 4.5pt off 4.5pt}] (140.74,40.1) .. controls (140.74,20.42) and (156.69,4.47) .. (176.37,4.47) .. controls (196.05,4.47) and (212,20.42) .. (212,40.1) .. controls (212,59.77) and (196.05,75.73) .. (176.37,75.73) .. controls (156.69,75.73) and (140.74,59.77) .. (140.74,40.1) -- cycle ;

\draw    (131.71,96.83) -- (150.43,153.24) ;
\draw    (225.5,78.74) -- (284.5,78.74) ;
\draw [shift={(286.5,78.74)}, rotate = 180] [color={rgb, 255:red, 0; green, 0; blue, 0 }  ][line width=0.75]    (10.93,-3.29) .. controls (6.95,-1.4) and (3.31,-0.3) .. (0,0) .. controls (3.31,0.3) and (6.95,1.4) .. (10.93,3.29)   ;
\draw   (384.07,83.11) .. controls (384.07,74.11) and (391.36,66.82) .. (400.36,66.82) .. controls (409.35,66.82) and (416.64,74.11) .. (416.64,83.11) .. controls (416.64,92.1) and (409.35,99.39) .. (400.36,99.39) .. controls (391.36,99.39) and (384.07,92.1) .. (384.07,83.11) -- cycle ;
\draw    (409.71,69.83) .. controls (468,-10.76) and (511,120.24) .. (415.71,88.83) ;
\draw    (69.13,150.81) -- (112.36,95.11) ;
\draw    (39.43,83.24) -- (106.07,83.11) ;
\draw    (80.43,18.24) -- (116.07,68.11) ;
\draw  [color={rgb, 255:red, 4; green, 146; blue, 194 }  ,draw opacity=1 ][dash pattern={on 4.5pt off 4.5pt}] (106.07,83.11) .. controls (106.07,74.11) and (113.36,66.82) .. (122.36,66.82) .. controls (131.35,66.82) and (138.64,74.11) .. (138.64,83.11) .. controls (138.64,92.1) and (131.35,99.39) .. (122.36,99.39) .. controls (113.36,99.39) and (106.07,92.1) .. (106.07,83.11) -- cycle ;
\draw   (384.07,83.11) .. controls (384.07,74.11) and (391.36,66.82) .. (400.36,66.82) .. controls (409.35,66.82) and (416.64,74.11) .. (416.64,83.11) .. controls (416.64,92.1) and (409.35,99.39) .. (400.36,99.39) .. controls (391.36,99.39) and (384.07,92.1) .. (384.07,83.11) -- cycle ;
\draw    (409.71,96.83) -- (428.43,153.24) ;
\draw    (347.13,150.81) -- (390.36,95.11) ;
\draw    (317.43,83.24) -- (384.07,83.11) ;
\draw    (358.43,18.24) -- (394.07,68.11) ;
\draw   (384.07,83.11) .. controls (384.07,74.11) and (391.36,66.82) .. (400.36,66.82) .. controls (409.35,66.82) and (416.64,74.11) .. (416.64,83.11) .. controls (416.64,92.1) and (409.35,99.39) .. (400.36,99.39) .. controls (391.36,99.39) and (384.07,92.1) .. (384.07,83.11) -- cycle ;
\draw (115,78) node [anchor=north west][inner sep=0.75pt]   [align=left] {$\displaystyle v$};
\end{tikzpicture}
\]

\caption{}\label{fig: contraction}
\end{figure}

\end{example}

\begin{theorem}\label{thm: nerve}\cite[Theorem 3.6]{hry2}
 The nerve functor \[\begin{tikzcd} N: \ModOp \arrow[r] & \mathbf{Set}^{\mathbf{U}^{op}} \end{tikzcd} \]  is fully faithful. Moreover, for any $X\in  \mathbf{Set}^{\mathbf{U}^{op}}$, the following statements are equivalent: 
\begin{enumerate} 
\item There exists a modular operad $\bP$ and an isomorphism $X\cong N\bP$.
\item $X$ satisfies the strict Segal condition.
\end{enumerate} 
\end{theorem}

\begin{remark} 
There are several related constructions and results in the literature.  In \cite{Sophie_phdthesis,Sophie_21}, Raynor presents a slightly larger graphical category (the category of graphical species) together with a monad for modular operads which allows her to avoid issues with the nodeless loop in Remark~\ref{rmk: loop}. This has the formal advantage that her graphical category embeds fully in the category of modular operads, but the practical disadvantage that the resulting corresponding construction of modular $\infty$-operads (\cite[Section 8.4]{Sophie_21}) is somewhat opaque. 

We also note that our definition of Segal objects (eg: Definition~\ref{def: Segal map} and Definition~\ref{def: weak Segal 1}) have been generalized by several authors. In particular, the definition of Segal objects in Example 3.11 of Chu–Haugseng \cite{MR4256131} and Berger's unital hypermoment categories \cite{berger2021moment} agrees with the one we have given here. 
\end{remark}

\subsection{Further Directions}\label{sec:further directions 1}
In the next lecture we will describe how weakening the Segal condition gives us a model for cyclic and modular $\infty$-operads.  At the time of these lectures, there are many open questions one would want to see answered before we can say that we have a comprehensive understanding of what a modular $\infty$-operad should be. The following are a few open problems.

In Definition~\ref{def: outer coface} and \ref{def: inner coface}, we defined the notion of (inner and outer) coface maps of $\bU$. Given a coface map $\delta$ with codomain $G$, one can define the \textbf{horn} $\Lambda^{\delta} [G]$ which is a sub-object of the representable object $\bU[G]$. A \emph{strict inner Kan} graphical set is a presheaf $X\in \mathbf{Sets}^{\mathbf{U}^{op}}$ such that every diagram \[\begin{tikzcd} \Lambda^{\delta} [G] \arrow[d] \arrow[r] & X\\ \bU[G]\arrow[ur, dotted]\end{tikzcd}\] with $\delta$ an inner coface map admits a unique filler. Michelle Strumila shows in her PhD thesis that : 

\begin{theorem}[Strumila]
The nerve functor \[\begin{tikzcd} N: \ModOp \arrow[r] & \mathbf{Set}^{\mathbf{U}^{op}} \end{tikzcd} \]  is fully faithful. Moreover, the following statements are equivalent for $X\in  \mathbf{Set}^{\mathbf{U}^{op}}.$ 
\begin{enumerate} 
\item There exists a modular operad $\bP$ and an isomorphism $X\cong N\bP$.
\item $X$ satisfies the strict Segal condition.
\item $X$ is strict inner Kan. 
\end{enumerate} 
\end{theorem} 

If one relaxes the inner Kan condition you arrive at a model for \emph{quasi-modular operads}. 

\begin{problem}\label{problem: fibrant kan}
Following the example of dendroidal sets \cite{cm_quasi_operad} find a model category structure in which the weak inner Kan graphical sets are the fibrant objects. 
\end{problem}

The involution on colour sets in Definition~\ref{def: modular} allows us to consider wheeled properads as a subcategory of $\ModOp$. In \cite{phil_graphs}, Hackney makes this explicit at the level of graphical categories, identifying a slice category $\bU/_{o}$ with the graphical category for wheeled properads defined in \cite{hry_hha}, \cite{hry15}. In particular, he shows that the adjunction \[\begin{tikzcd} \mathbf{WPrd} \arrow[r, shift left =.15cm] & \arrow[l, shift left =.15cm] \ModOp\end{tikzcd}\] can be well understood via graphical presheaves. Similar adjunctions between modular operads, cyclic operads, and operads can all be described via adjunctions of graphical categories. 

\begin{problem}
Assuming a solution to Problem~\ref{problem: fibrant kan}, use the adjunctions of graphical categories described in \cite{phil_graphs} to define a model category structure on the category of graphical sets from \cite{hry15}, \cite{hry_hha} in which the quasi-wheeled properads are the fibrant objects. 
\end{problem}

\section{Lecture 2: A weak Segal model for modular $\infty$-operads} 


For the remainder of this lecture series we will simplify Definition~\ref{def: modular} and focus on the \emph{one-coloured} modular operads. Just to refresh our memory, a symmetric sequence $\bP=\{\bP(n)\}$ consists of a sequence of sets $\bP(n)$ each of which is equipped with a right $\Sigma_n$-action.\footnote{We follow the convention of calling a one-coloured collection a symmetric sequence.} A \emph{modular operad} $\bP$ consists of a symmetric sequence $\bP=\{\bP(n)\}$ together with:
\begin{enumerate}
\item A distinguished \emph{unit} element $1\in \bP(1)$; 
\item A family of \emph{compositions} \[\begin{tikzcd} \bP(n)\times \bP(m) \arrow[r, "\circ_{ij}"]& \bP(n+m-2);\end{tikzcd}\]
\item A family of \emph{contraction} operations \[\begin{tikzcd} \bP(n) \arrow[r, "\xi_{ij}"]& \bP(n-2).\end{tikzcd}\]
\end{enumerate} 
Moreover, we require the compositions and contractions satisfy a series of axioms (eg: \cite[Definition A1]{batanin2021koszul}) assuring that compositions are associative, unital and equivariant, contractions are associative and equivariant and that the two operations are compatible. 

\medskip 
Theorem~\ref{thm: nerve} tells us is that, given a discrete modular operad $\bP$, we can construct a set-valued presheaf $N\bP\in (\textbf{Set}^{\bU^{op}})_{Segal}$ where $$N\bP_{\medwhitestar_n}= \bP(n)$$ in which the Segal maps \[\begin{tikzcd} N\bP_G \arrow[r] & \prod\limits_{v\in V(G)} N\bP_{\medwhitestar_v} \end{tikzcd}\] are bijections. In other words, modular operad compositions and contractions of $\bP$ are modeled by graphical maps
 \[\begin{tikzcd} \prod\limits_{v\in V(G)} N\bP_{\medwhitestar_v}\arrow[r, shift left =.2cm, dotted]& N\bP_G \arrow[r] \arrow[l, shift left =.25, "Segal"] & N\bP_{\medwhitestar_G}. \end{tikzcd}\]


If our goal is to now model modular $\infty$-operads, i.e. modular operads where operations are defined ``up to coherent homotopy''. This means we will want to replace our Segal map with a homotopy equivalence. In this second lecture we will introduce space-valued presheaves $\mathbf{sSet}^{\bU^{op}}$ and describe a corresponding notion of a \emph{weak} Segal condition on modular dendroidal spaces. At the end of this lecture, we include some brief notes about variations on the graphical category $\bU$ which can give genus graded modular $\infty$-operads, cyclic $\infty$-operads, etc.

\subsection{Modular dendroidal  spaces}
Throughout, we write $\textbf{sSet}$ for category of simplicial sets equipped with the standard Kan-Quillen model structure. We often abuse terminology and refer to simplicial sets as ``spaces''. 

\begin{definition} \label{def: modular spaces}
The category of \emph{modular dendroidal  spaces} is the category of space-valued $\bU$-presheaves denoted by $\mathbf{sSet}^{\bU^{op}}$.  
\end{definition} 

As in the previous lecture, for any $X\in\mathbf{sSet}^{\bU^{op}}$ we write $X_G$ for the evaluation of $X$ at a graph $G\in \bU$.  We consider the representable presheaf in $\mathbf{Set}^{\bU^{op}}$: \[\bU[G]_H:=\bU(H, G)\] as an object in $\mathbf{sSet}^{\bU^{op}}$ via the inclusion $\mathbf{Set}^{\bU^{op}} \hookrightarrow \mathbf{sSet}^{\bU^{op}}$.

\begin{exercise} 
The Yoneda Lemma says that a map $x: \bU[G]\rightarrow X$ in $\textbf{Set}^{\bU^{op}}$ is equivalent to an element $x \in X_G$. Show that every $X\in \textbf{Set}^{\bU^{op}}$ is, up to isomorphism, a colimit of representables \[X \cong \colim \bU[G]\]  where the colimit is indexed by the maps $\bU[G]\rightarrow X$.
\end{exercise} 
 
\subsubsection{Segal cores}
In the previous lecture, we introduced the Segal maps via a limit construction (Definition~\ref{def: ribbon}). To describe the weak Segal maps, we will use a dual construction called the \emph{Segal core}. 

To understand the Segal core construction, it can be useful to revisit the definition of a graph. Recall that a graph $G$ as a diagram in the category of finite sets in the shape of 
\[\begin{tikzcd} \calI : = & \arrow[loop left]{l}{i}\bullet & \arrow[l, swap, "s"] \bullet \arrow[r, "t"] & \bullet \end{tikzcd} \] where the arrow $s$ is sent to a monomorphism and the generating endomorphism $i$ is a free involution. If a graph $G$ has at least one vertex, we can choose an orientation for each internal edge, and present $G$ as a coequalizer in $\textbf{FinSet}^{\calI}$:\footnote{Note this is \emph{not} a coequalizer in $\bU$ as these objects don't exist in $\bU$ (they are not connected).}
\[\begin{tikzcd}  \coprod\limits_{e\in \textrm{iE}}\updownarrow \arrow[r, shift right =.1cm] \arrow[r, shift left =.1cm] & \coprod\limits_{v\in V}\medwhitestar_v \arrow[r]& G. \end{tikzcd}\]

\begin{exercise}
Write the graph $G$ from Figure~\ref{example:ribbon} as a coequalizer: 
\[\begin{tikzcd}  \{e_1,e_2,e_3\} \arrow[r, shift right =.1cm] \arrow[r, shift left =.1cm] & \medwhitestar_v\cup \medwhitestar_w \arrow[r]& G. \end{tikzcd}\]
\end{exercise}

\begin{definition}\label{def: Segal core}
The \emph{Segal core} of a graph $G$ is the coequalizer in $\mathbf{Set}^{\bU^{op}}$: 
\[\begin{tikzcd}  \coprod\limits_{e\in \textrm{iE}}\mathbf{U}[\updownarrow] \arrow[r, shift right =.15cm] \arrow[r, shift left =.15cm] & \coprod\limits_{v\in V}\mathbf{U}[\medwhitestar_v] \arrow[r]& \mathbf{Sc}[G]  \end{tikzcd}\] in $\mathbf{Set}^{\bU^{op}}$. 
\end{definition} 

The Segal core comes with a natural map $\mathbf{Sc}[G] \rightarrow \mathbf{U}[G]$ induced by the embeddings $\medwhitestar_v\rightarrow G$. In the case that $G=\updownarrow$ we declare that the map $\mathbf{Sc}[G]\rightarrow \mathbf{U}[G]$ to be the identity map on $\mathbf{U}[G]$.

\begin{exercise}
Check that the Segal core is precisely the colimit so that \[\mathbf{Set}^{\bU^{op}}(\mathbf{Sc}[G], X) = X^{1}_{G} =  \lim\limits_{ \medwhitestar_v\leftarrow \updownarrow\rightarrow \medwhitestar_w} \left( \begin{tikzcd} X_{\medwhitestar_{v}} \arrow[dr] & & \arrow[dl] X_{\medwhitestar_{w}} \\&X_{\updownarrow} & \end{tikzcd}\right) .\]  
\end{exercise}

\medskip
As we saw in Exercise~\ref{exercise: limit is prod}, whenever $X_{\updownarrow}=\mathbf{Set}^{\bU^{op}}(\bU[\updownarrow], X)=\ast$, we can identify \[ \mathbf{Set}^{\bU^{op}}(\mathbf{Sc}[G], X)  = \prod\limits_{v\in V(G)} X_{\medwhitestar_v}.\] This implies that the Segal map from Definition~\ref{def: Segal map} is equivalent to the map of sets \[\begin{tikzcd} X_G=\mathbf{Set}^{\bU^{op}}(\mathbf{U}[G], X)  \arrow[rr] && \mathbf{Set}^{\bU^{op}}(\mathbf{Sc}[G], X) =\prod\limits_{v\in V(G)} X_{\medwhitestar_v} .\end{tikzcd}\]

\medskip 

As with the representable presheaves $\bU[G]$, we consider the Segal core $\mathbf{Sc}[G]$ as an object in $\mathbf{sSet}^{\bU^{op}}$ via the inclusion $\mathbf{Set}^{\bU^{op}} \hookrightarrow \mathbf{sSet}^{\bU^{op}}$. This leads to the following definition: 

\begin{definition}\label{def: weak Segal 1}
A modular dendroidal space $X\in \textbf{sSet}^{\bU^{op}}$ is \emph{Segal} if: 
\begin{itemize} 
\item $X_{\updownarrow}=\ast$;
\item  for all $G\in\bU$, the Segal map \[\begin{tikzcd} \map(\mathbf{U}[G], X)  \arrow[r] & \map(\mathbf{Sc}[G], X)\end{tikzcd}\] is a weak equivalence. 
\end{itemize} 
\end{definition} 

Here $\map(X,Y)$ is the derived mapping space, which is well-defined as long as we can equip the category $\mathbf{sSet}^{\bU^{op}}$ with a class of weak equivalences. The category $\mathbf{sSet}^{\bU^{op}}$ admits several model category structures including the \emph{projective model structure} and a \emph{Reedy model structure}.  In either case, weak equivalences are defined entrywise. In other words, $f: X\rightarrow Y$ in $\mathbf{sSet}^{\bU^{op}}$ is a weak equivalence if $f^*: X_{G} \rightarrow Y_{G}$ is a weak equivalence of $\mathbf{sSet}$ for every $G\in\bU$. 

\begin{remark} 
The assumption that $X_{\updownarrow}= \ast$ is not required for Definition~\ref{def: weak Segal 1} but is required in Theorem~\ref{thm: Segal are fibrant}. We have included the assumption here for consistency throughout these notes. 
\end{remark} 

Definition~\ref{def: weak Segal 1} is a perfectly fine definition for modular dendroidal Segal spaces. In our intended applications, however, we will want to demonstrate that a specific presheaf $X\in\mathbf{sSet}^{\bU^{op}}$ is Segal and this simplifies significantly when one uses the Reedy model structure on the category $\textbf{sSet}^{\bU^{op}}$. 

\subsection{Generalized Reedy Categories} 
The notion of a generalized Reedy category was introduced in \cite[Definition 1.1]{bm_reedy}.\footnote{As an interesting historical note, we noticed while preparing these notes that the initial results on generalized Reedy categories were presented at the CRM program on Homotopy Theory
and Higher Categories in 2008.} 

\begin{definition} Let $\mathbb{R}$ be a small category. A dualizable generalized Reedy structure on $\mathbb{R}$ consists of two wide subcategories \[\mathbb{R}^{+} \quad \text{and} \quad \mathbb{R}^{-}\] together with a \emph{degree function} on objects $\textnormal{ob}(\mathbb{R}) \rightarrow \mathbb{N}$ satisfying:

\begin{enumerate} 
\item non-invertible morphisms in $\mathbb{R}^{+}$ (respectively $\mathbb{R}^{-}$) raise (respectively lower
degree). Isomorphisms preserve degree. 

\item $\mathbb{R}^{+} \bigcap \mathbb{R}^{-} = \textnormal{Iso}(\mathbb{R})$

\item Every morphism $f$ factors as $f = gh$ such that $g\in \mathbb{R}^{+}$ and $h\in \mathbb{R}^{-} $. Moreover, this factorization is unique up to isomorphism.

\item If $\theta f = f$ for any isomorphism $\theta$ and $f \in \mathbb{R}^{-} $ then $\theta$ is an identity. 

\item If $f \theta = f$ for any isomorphism $\theta$ and $f \in \mathbb{R}^{+} $ then $\theta$ is an identity. 
\end{enumerate} 

\end{definition} 

The subcategory $\mathbb{R}^{+}$ is commonly called the ‘direct category’ and $\mathbb{R}^{-}$ the ‘inverse category.’  A category $\mathbb{R}$ that satisfies axioms $(1)-(4)$ is a generalised Reedy category. If, in addition, $\mathbb{R}$ satisfies axiom $(5)$ then $\mathbb{R}$ is said to be dualizable, which implies that $\mathbb{R}^{op}$ is also a generalised Reedy category. 

\begin{example} 
The simplicial category $\Delta$ is a Reedy category in which every isomorphism is an identity. 
\end{example} 

\begin{example} 
Other examples of generalized Reedy categories include the \emph{dendroidal category} $\Omega$, finite sets, and pointed finite sets. 
\end{example} 

The main use of Reedy categories is that one can use latching and matching objects to lift morphisms from $\mathbb{R}$ to the diagram category $\bE^{\mathbb{R}}$ by induction on the degree of objects.  For any object $r\in\mathbb{R}$, the category $\mathbb{R}^{+}(r)$ is the full subcategory of $\mathbb{R}^{+}\downarrow r$ consisting of non-invertible maps with target $r$. Similarly, the category $\mathbb{R}^{-}(r)$ is the full subcategory of $r \downarrow \mathbb{R}^{-}$ consisting of the non-invertible maps $\alpha: r\rightarrow s$. 

\begin{definition} 
Let $X$ be a diagram in $\bE^{\mathbb{R}}$
\begin{itemize} 
\item The latching object $L_rX =\colim_{\mathbb{R}^{+}(r)} X$;
\item The matching object $M_rX =\lim_{\mathbb{R}^{-}(r)} X$. 
\end{itemize} 
\end{definition} 

If $\bE$ is a cofibrantly generated model category.  We say that a morphism $f: X\rightarrow Y$ in $\bE^{\mathbb{R}}$ is: 
\begin{itemize} 
\item a Reedy cofibration if $X_r\bigcup_{L_rX} L_rY \rightarrow Y_r$ is a cofibration in $\bE^{\Aut(r)}$ for all $r\in \mathbb{R}$;
\item a Reedy weak equivalence if $X_r\rightarrow Y_r$ in $\bE^{\Aut(r)}$ for all $r\in \mathbb{R}$;
\item a Reedy fibration if $X_r\rightarrow M_rX\times_{M_rY} Y_r$ in $\bE^{\Aut(r)}$ for all $r\in \mathbb{R}$. 
\end{itemize} 

\begin{theorem}\cite{bm_reedy}
If $\mathbb{R}$ is a dualizable generalized Reedy category and $\bE$ is a nice enough model category, then $\bE^{\mathbb{R}^{op}}$ admits a cofibrantly generated model category structure with level-wise weak equivalences. 
\end{theorem}

\subsubsection{The Reedy structure on $\bU$}
The graphical category $\bU$ has many nice factorization properties (eg. Remark 1.8 \cite{hry2}) and is, in particular, a generalized Reedy category. We define the \emph{degree of a graph} $G$ to be \[\textrm{deg}(G)= |V|+ |iE|.\] Then the degree function $\textrm{deg}:\textrm{ob}(\bU) \rightarrow \mathbb{N}$. 

\begin{theorem}\label{thm: reedy}\cite[Theorem 2.22]{hry1}
The graphical category $\bU$ is a (dualizable) generalised Reedy category. The wide subcategory $\bU^{-}$ is generated by the codegeneracy maps (Definition~\ref{def: deg}) and the wide subcategory $\bU^{+}$ is generated by the inner and outer coface maps (Definition~\ref{def: inner coface} and Definition~\ref{def: outer coface}). 
\end{theorem} 

Applying Theorem 1.6 of \cite{bm_reedy}, we have the following corollary: 

\begin{cor} 
The diagram category $\textbf{sSet}^{\bU^{op}}$ admits a model category structure with the Reedy fibrations, Reedy cofibrations, and entrywise weak equivalences.
\end{cor}

\begin{exercise}[Hard-ish]
In Proposition 3.5 of \cite{hry1} we show that Segal cores are cofibrant in the Reedy model structure on $\textbf{sSet}^{\bU^{op}}$. Give an example of a graph $G$ in which the Segal core of $G$ fails to be cofibrant in the projective model structure. 
\end{exercise} 

The advantage of a Reedy model structure on $\textbf{sSet}^{\bU^{op}}$ is that homotopy limits of Reedy fibrant diagrams are just limits. Revisiting our Definition~\ref{def: weak Segal 1} we now have: 

\begin{definition}\label{def: weak Segal 2} 
A modular dendroidal  space $X\in \textbf{sSet}^{\bU^{op}}$ is \emph{Segal} if: 

\begin{enumerate}
\item $X_{\updownarrow}=\ast$;
\item $X$ is Reedy fibrant;
\item for all $G\in\bU$, the Segal map \[\begin{tikzcd} X_{G}= \map(\mathbf{U}[G], X)  \arrow[rr] && \map(\mathbf{Sc}[G], X) = \prod\limits_{v\in V(G)} X_{\medwhitestar_{v}}\end{tikzcd}\] is a weak equivalence in $\mathbf{sSet}$. 
\end{enumerate} 
\end{definition} 

In the final lecture of this series, we will give our motivating example of a modular $\infty$-operad. We note, however, that given any (one-coloured) modular operad in $\mathbf{sSet}$, $\bP$, the nerve $N\bP$ is a Segal modular operad in the sense of Definition~\ref{def: weak Segal 1}. Moreover, if $N\bP$ is Reedy fibrant, then the Segal map \[\begin{tikzcd} N\bP_G \rightarrow \prod_{v\in V} N\bP_{\medwhitestar_v} \end{tikzcd}\] is an isomorphism for every $G$.  Thus, up to fibrant replacement, every (one-coloured) modular operad gives rise to a modular $\infty$-operad.  We conclude our description of modular $\infty$-operads by pointing out that there is a classification of modular dendroidal Segal spaces as the fibrant objects in a localization of the Reedy model category structure on graphical spaces. 

\begin{theorem}\cite[Theorem 3.8; Proposition 3.19]{hry1}\label{thm: Segal are fibrant}
The category $\mathbf{sSet}^{\bU^{op}}$ admits a cofibrantly generated model structure whose fibrant objects are the Segal modular operads.
\end{theorem} 

\begin{remark} 
While we have discussed space-valued presheaves everything contained in this second lecture about modular dendroidal  spaces still makes sense for presheaves in any Cartesian monoidal model category $\mathbf{C}$. 
\end{remark}

\subsection{Variations on the graphical category $\bU$ and open problems}\label{sec:open problems 2}
The definition of modular operad often comes equipped with an additional 
``genus'' grading (eg: \cite{gk_modular}, \cite{modular_truncated}, 
\cite{dotsenko2021deformation}). A \emph{graded} modular operad consists of a 
bi-graded sequence $\bP=\{\bP(g,n)\}$, in which each $\bP(g,n)$ is equipped with 
a right action of the symmetric group $\Sigma_n$, together with units, 
composition and contraction maps. Often one also requires that the underlying 
collection of $\bP=\{\bP(g,n)\}$ satisfy a \emph{stability condition}, i.e. 
\[\bP(g,n) = \emptyset \quad \text{whenever} \quad 2g+n -2\leq 0.\] One can 
define the corresponding genus graded version of the graphical category $\bU$ and a corresponding stable version of modular $\infty$-operads. 

\begin{definition} 
Let $G$ be a graph: 
\begin{enumerate} 
\item A genus function for $G$ is a function $g:V(G)\rightarrow \mathbb{N}$.
\item The total genus of a pair $(G, g:V\rightarrow \mathbb{N})$ is given by: \[g(G)=\beta_1(G)+ \Sigma_{v\in V} g(v)\]
where $\beta_1(G)$ is the Betti number of $G$. 
\item A pair $(G, g)$ is called \textbf{stable} if $G$ is connected and for every vertex v: \[2g(v)+ |\nb(v)| - 2 > 0.\]
\end{enumerate} 

\end{definition} 

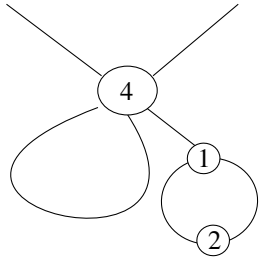
\begin{figure}

\[\begin{tikzpicture}[x=0.75pt,y=0.75pt,yscale=-.75,xscale=.75]

\draw  [fill={rgb, 255:red, 255; green, 255; blue, 255 }  ,fill opacity=1 ] (316.36,107.97) .. controls (327.4,107.9) and (336.42,115.14) .. (336.52,124.16) .. controls (336.61,133.18) and (327.74,140.55) .. (316.7,140.63) .. controls (305.67,140.71) and (296.64,133.46) .. (296.55,124.44) .. controls (296.46,115.42) and (305.33,108.05) .. (316.36,107.97) -- cycle ;
\draw    (235.16,66.52) -- (299.1,114.43) ;
\draw    (333.73,114.43) -- (391.01,66.52) ;
\draw    (329.68,136.21) -- (366.72,164.84) ;
\draw    (296.68,135.63) .. controls (124.21,199.92) and (394.61,261.56) .. (316.7,140.63) ;
\draw   (339.6,202.15) .. controls (337.42,186.37) and (350.02,171.84) .. (367.74,169.7) .. controls (385.45,167.57) and (401.58,178.64) .. (403.76,194.42) .. controls (405.93,210.21) and (393.34,224.74) .. (375.62,226.87) .. controls (357.9,229.01) and (341.78,217.94) .. (339.6,202.15) -- cycle ;
\draw  [fill={rgb, 255:red, 255; green, 255; blue, 255 }  ,fill opacity=1 ] (366.58,159.56) .. controls (372.59,158.95) and (377.99,162.99) .. (378.64,168.6) .. controls (379.28,174.2) and (374.92,179.24) .. (368.91,179.85) .. controls (362.89,180.47) and (357.49,176.42) .. (356.84,170.82) .. controls (356.2,165.21) and (360.56,160.17) .. (366.58,159.56) -- cycle ;
\draw  [fill={rgb, 255:red, 255; green, 255; blue, 255 }  ,fill opacity=1 ] (372.98,214.44) .. controls (379,213.83) and (384.4,217.88) .. (385.04,223.48) .. controls (385.69,229.09) and (381.33,234.13) .. (375.31,234.74) .. controls (369.3,235.35) and (363.9,231.3) .. (363.25,225.7) .. controls (362.61,220.1) and (366.97,215.06) .. (372.98,214.44) -- cycle ;

\draw (310,116.4) node [anchor=north west][inner sep=0.75pt] [font=\footnotesize]   {$4$};
\draw (362,161.4) node [anchor=north west][inner sep=0.75pt] [font=\footnotesize]   {$1$};
\draw (369,216.4) node [anchor=north west][inner sep=0.75pt] [font=\footnotesize]   {$2$};

\end{tikzpicture}\]
\caption{An example of a graph with genus}\label{fig:genus graph}
\end{figure}

For example, in Figure~\ref{fig:genus graph}, the graph $G$ has genus \[g(G)=\beta_1(G)+ \Sigma_{v\in V} g(v) = 2+4+1+2=9.\]  Given an embedding $f:H \rightarrow G$, then we can define the \emph{genus of $f$} \[g(f) := \beta_1(H) + \Sigma_{v\in V(H)} g(f(v)).\]

\begin{remark}
The exceptional edge admits only one genus function, and $G(\updownarrow) =\beta_1(\updownarrow) =0.$ This graph trivially satisfies the stability condition.
\end{remark}

\begin{definition} 

The stable graphical category $\bU_{st}$ has: 
\begin{itemize} 
\item Objects: stable graphs $(G, g)$ 
\item Morphisms: $(G,g)\rightarrow (G',g')$ are graphical maps $\varphi: G\rightarrow G'$ which make the diagram: 
\[\begin{tikzcd}
V(G) \arrow[dr, "g", swap]\arrow[rr, "\varphi_1"] && \Emb(G')\arrow[dl, "g' "]\\
& \mathbb{N} & 
 \end{tikzcd}\] commute. 
\end{itemize} 
\end{definition} 

The stability condition ensures that in the stable graphical category $\bU_{st}$ there are no codegeneracy maps, because a genus $0$ vertex with $|\nb(v)|=2$ cannot be stable and thus the substitution of the edge into an arity $2$ vertex is not in our category. This makes the following theorem an immediate corollary to Theorem~\ref{thm: reedy}.  
\begin{theorem} 
$\bU_{st}$ is a generalized Reedy category. 
\end{theorem} 

We can therefore define a stable version of a modular $\infty$-operad as follows: 
\begin{definition} 
There is model structure on $\textbf{sSet}^{\bU_{st}^{op}}$ in which $X\in \textbf{sSet}^{\bU_{st}^{op}}$ is fibrant if: 

\begin{itemize} 
\item $X_{\updownarrow}=\ast$;
\item $X$ is Reedy fibrant;
\item  for all $G\in\bU$, the Segal map \[\begin{tikzcd} X_{(G,g)}= \map^{h}(\mathbf{U}[(G,g)], X)  \arrow[r] & \map^h(\mathbf{Sc}[(G,g)], X) \end{tikzcd}\] is a weak equivalence. 
\end{itemize} 
\end{definition} 
 
\subsubsection{Cyclic operads}
As we mentioned in the first lecture, modular operads are cyclic 
operads with contraction operations. In the literature there are actually 
various notions of cyclic operad, and we define various subcategories of 
our graphical category $\bU$ which correspond to the reader's desired notion of cyclic operad. In particular, there is nested sequence of subcategories: \[\bU_{cyc}\subset \bU_0 \subset \bU\] defined as follows: 
\begin{enumerate}
\item The category $\bU_0$ is the full subcategory of $\bU$ whose objects are all \emph{simply connected} graphs. The category $\bU_0$ corresponds to \emph{augmented cyclic operads}. 
\item The category $\bU_{cyc}$ is the full subcategory of $\bU$ whose objects are all \emph{simply connected} graphs with \emph{non-empty boundary}. The category $\bU_{cyc}$ corresponds to \emph{cyclic operads}. 
\end{enumerate}

\begin{exercise}
Show that $\bU_0$ and $\bU_{cyc}$ are sieves of $\bU$. In other words if $\varphi:G \rightarrow T$ is in $\bU$ with $T\in \bU_0$ (respectively, $\bU_{cyc}$) then $G\in \bU_0$ (respectively, $\bU_{cyc}$). 
\end{exercise} 

We also note that the category $\bU_{cyc}$ is itself related to other categories in the literature: 

\begin{enumerate}
\item In \cite{walde} Walde introduces a category $\Omega_{cyc}$, which is a non-symmetric version of $\bU_{cyc}$. That is: $\bU_{cyc}$ is equivalent to a category $\bU'_{cyc}$ in which every object has a specified cyclic ordering and $\Omega_{cyc}$ is the wide subcategory of $\bU'_{cyc}$ where maps preserve the ordering. 

\item There is another category of Segal cyclic operads $\Xi$ in \cite{hry_cyc}.  This is a graphical category in which models coloured cyclic operads where the involution on colour sets is always trivial. In practice, this category has the same objects as $\bU_{cyc}$ but slightly different morphisms. 
\end{enumerate} 

There are well understood (and very useful) adjunctions between operads cyclic operads and modular operads. One would hope that these same relationships hold between $\infty$-versions of all these objects. This inspires the following open problems. 

\begin{problem}
Show there are Quillen adjunctions \[\begin{tikzcd} \textbf{sSet}^{\Xi^{op}}  \arrow[r, shift left = .1cm] & \textbf{sSet}^{\bU_{cyc}^{op}} \arrow[l, shift left=.1cm]  \arrow[r, shift right=.1cm] & \arrow[l, shift right =.1cm] \textbf{sSet}^{\Omega^{op}}.\end{tikzcd}\] Is the adjunction 
 \[\begin{tikzcd} \textbf{sSet}^{\Xi^{op}}  \arrow[r, shift left = .1cm] & \textbf{sSet}^{\bU_{cyc}^{op}} \arrow[l, shift left=.1cm] \end{tikzcd}\] a Quillen equivalence ? The interested reader may want to look at \cite[Proposition 8.5]{hry_cyc} where we establish a Quillen adjunction  \[\begin{tikzcd} \textbf{sSet}^{\Xi^{op}}  \arrow[r, shift right = .1cm] & \arrow[l, shift right =.1cm] \textbf{sSet}^{\Omega^{op}}\end{tikzcd}\] for inspiration. 
\end{problem}

\begin{problem}
 Work of Barwick \cite{barwick},  Hirschhorne and Volic \cite{hv}  characterizes when $F:\mathbb{R} \rightarrow \mathbb{S}$ between \emph{strict} Reedy categories result in Quillen adjunctions between diagram categories. Is there a similar characterization for generalized Reedy categories? 
\end{problem}

\section{Lecture 3: Lego-Teichm\"uller theory and modular operads} 
In this final lecture, we introduce the genus graded modular operad built from surfaces which, after profinite completion, is related to the ideal Teichm\"uller towers in our introduction. Throughout this final lecture we will often make use of the fact that we have adjunctions
\begin{equation}\label{adjunction}
    \begin{tikzcd} \textbf{Operad} \arrow[r, shift left = .15cm] & \textbf{Cyc} \arrow[l, shift left=.15cm]  \arrow[r, shift left=.15cm] & \arrow[l, shift left =.15cm, "\tau_{0}"]\ModOp  .\end{tikzcd}
\end{equation} The adjunction between operads and cyclic operads is explicitly described in Section 3 of \cite{DCH_Dwyer_Kan}. The map $\tau_{0}:\ModOp\rightarrow \textbf{Cyc}$ is ``truncation at genus $0$'' or, equivalently, forgetting all contraction operations. This functor is actually a special case of a family of adjunctions 
\[\begin{tikzcd} \ModOp_{k}  \arrow[r, shift left = .15cm] & \ModOp \arrow[l, shift left=.15cm,  "\tau_{k}"] \end{tikzcd}\] where the map $\tau_{k}$ is truncation at genus $k$, meaning we forget all operations of genus $\geq k$. These functors are a straightforward generalization of those in \cite[Section 4.1; 8.4]{modular_truncated} or \cite[Section 2.6]{ward_modular}.

\subsection{Profinite completion of modular operads in groupoids} 

In a Cartesian monoidal category $\bE$, an \emph{inverse system}, $\calI$, is a collection of objects $\{X_i\}_{i\in I}$ in $\bE$, together with maps $\phi_{ij}: X_i \to X_j$ for $i\geq j$, such that:
\begin{enumerate}
    \item \quad$\phi_{ii}: X_i \to X_i$ is the identity $\textrm{id}_{X_i}$;
    \item \quad$\phi_{ij} \circ \phi_{jk} = \phi_{ik}$ for $i\geq j\geq k$. 
\end{enumerate}
The limit over $\calI$ is then given by
$$\lim\limits_{\longleftarrow} X_i = \Big\{ (x_i) \in \prod\limits_{i\in I} X_i \;|\; \phi_{ij}(X_i) = X_j \;\forall i \geq j\Big\}.$$

\begin{example}
A group $G$ has an associated inverse system of finite index subgroups $\{G/N_i\}_{i\in\calI}$ where the $N_i$ runs over all normal subgroups of $G$ and the maps $$\phi_{ij}: G / N_i \to G / N_j$$ are the natural projections. 
\end{example}

\begin{definition}\label{def: profinite group} Given a finite group $G$, the \emph{profinite completion} of $G$ is the limit \[ \widehat{G} = \lim\limits_{\longleftarrow} G/N\] inverse system of finite index subgroups. 
\end{definition} 

\begin{example}
Let $G = \mathbb{Z}$ which has finite index subgroups $\mathbb{Z}/n\mathbb{Z}$. The profinite completion is the group of \emph{profinite integers} $\widehat{\mathbb{Z}} = \lim\limits_{\longleftarrow} \mathbb{Z}/n\mathbb{Z}$. Elements of $\widehat{\mathbb{Z}}$ are given by a sequence $(a_n)_{n\in\mathbb{N}}$ with $a_n \in \mathbb{Z}/n\mathbb{Z}$, and the structure maps for $n \geq m$ are the projections
$$\phi_{nm}: \mathbb{Z}/n\mathbb{Z} \to \mathbb{Z}/m\mathbb{Z}$$
which means that $a_n \equiv a_m \,\textrm{mod}\, m$ whenever $m\,|\,n$. 

\end{example}




Since profinite completion is a limit, there is a natural map $G\mapsto \widehat{G}$. These maps assemble into a functor $\widehat{(-)}:\textbf{Grp} \rightarrow \widehat{\textbf{Grp}}$ from the category of groups to the category of profinite groups. This functor is the left adjoint in an adjunction
\begin{equation}\label{completion adjunction grp}
\begin{tikzcd}\textbf{Grp} \arrow[r, shift left =.15cm, "\widehat{(-)}"] & \arrow[l, shift left =.15cm, "|-|"] \widehat{\textbf{Grp}}. \end{tikzcd}\end{equation} The right adjoint sends a profinite group to its underlying discrete group. 

We can generalize profinite completion of groups to define profinite completion of groupoids. In \cite{Horel_groupoid}, Horel extends the adjunction \eqref{completion adjunction grp} to the category of \emph{ groupoids}:   
\begin{equation}\label{completion adjunction grpd}
\begin{tikzcd}\textbf{Gpd} \arrow[r, shift left =.15cm, "\widehat{(-)}"] & \arrow[l, shift left =.15cm, "|-|"] \widehat{\textbf{Gpd}}. \end{tikzcd}\end{equation}

The category of profinite groupoids $\widehat{\textbf{Gpd}}$ admits a model category structure (\cite[Theorem 4.12]{Horel_groupoid}) and the adjunction \eqref{completion adjunction grpd} is a Quillen adjunction. In addition, profinite completion of groupoids is monoidal in the following sense:

\begin{prop}\cite[Proposition 4.23]{Horel_groupoid}\label{prof_monoidal}
Let $\bC$ and $\bD$ be two groupoids with finitely many objects. Then the map 
$$\widehat{\bC \times \bD} \to \widehat{\bC} \times \widehat{\bD}$$
induced by the projection maps $\widehat{\bC \times \bD} \to \widehat{\bC}$ and $\widehat{\bC \times \bD} \to \widehat{\bD}$ is an isomorphism.
\end{prop}

Proposition~\ref{prof_monoidal} enables us to define the \emph{profinite completion} of (cyclic and modular) operads in groupoids, so long as every entry of the (cyclic and modular) operad is a groupoid with a finite set of objects. 

\begin{definition}\label{def:profinte modular operad}
Let $\bP=\{\bP(n)\}$ be a modular operad in groupoids, in which each groupoid $\bP(n)$ has a finite set of objects. The profinite completion of $\bP$ is the modular operad \[\widehat{\bP} := \{\widehat{\bP}(n)\}\] where the profinite completion functor is applied entrywise. Composition operations are defined via the dashed lines:
\[\begin{tikzcd}
\widehat{\bP}(n) \times\widehat{\bP}(m) \arrow[d, "\cong"'] \arrow[r, dashed, "\circ_{ij}"] & \widehat{\bP}(n+m-2) \\
\widehat{\bP(n) \times \bP(m)}.\arrow[ur, swap, "\widehat{\circ}_{ij}" ] 
\end{tikzcd}\] Here the map $\widehat{\circ}_{ij}$ is the result of applying the profinite completion functor in \eqref{completion adjunction grpd} to the $\circ_{ij}$-composition maps of $\bP$.  The contraction operations $$\widehat{\xi}_{ij}: \widehat{\bP}(n) \to \widehat{\bP}(n-2)$$ result of applying the profinite completion functor \eqref{completion adjunction grpd} to the $\xi_{ij}$-composition maps of $\bP$.
\end{definition}

\begin{remark}
In Definition~\ref{def:profinte modular operad} we have restricted ourselves to one-coloured modular operads, but this is not necessary. Definition~\ref{def:profinte modular operad} holds for any $\fC$-coloured or genus-graded modular operad in groupoids, so long as each entry of modular operad only has finitely many objects. 
\end{remark}

The nerve theorem from the first lecture (Theorem~\ref{thm: nerve}) generalizes to groupoid-valued $\bU$-presheaves. We can therefore use the nerve functor to identify a modular operad $\widehat{\bP}$ in $\widehat{\textbf{Gpd}}$ with a presheaf $N\widehat{\bP}:\bU^{op}\rightarrow \widehat{\textbf{Gpd}}$ in which every graph $G$ gives an isomorphism of profinite groupoids \[\begin{tikzcd} N\widehat{\bP}_G \arrow[r]& \prod_{v\in V} N\widehat{\bP}_{\medwhitestar_v}. \end{tikzcd}\]

\subsection{Profinite completion of modular operads in spaces}
To every groupoid $\bG$ we can associate a space by taking the classifying space $B\bG$. This fits into an adjunction \[\begin{tikzcd} \mathbf{Gpd}\arrow[r, "B", swap, shift right=.15cm] & \arrow[l, "\pi_0", swap, shift right=.15cm] \mathbf{sSet}. \end{tikzcd}\] Both the classifying space and fundamental groupoid functors preserves products and thus we can lift this to an adjunction \[\begin{tikzcd} \ModOp(\mathbf{Gpd})\arrow[r, "B", swap, shift right=.15cm] & \arrow[l, "\pi_0", swap, shift right=.15cm] \ModOp(\mathbf{sSet}). \end{tikzcd}\] 

A \emph{profinite space} is a simplicial object in profinite sets. The category of profinite spaces is denoted $\widehat{\textbf{sSet}}$ and is equipped with the model category structure from Quick~\cite{quick_profinite_spaces} (See also the discussion in Section 3 of \cite{bhr1}). The profinite completion functor of spaces
\begin{equation}\label{completion adjunction space}
\begin{tikzcd}\textbf{sSet} \arrow[r, shift left =.15cm, "\widehat{(-)}"] & \arrow[l, shift left =.15cm, "|-|"] \widehat{\textbf{sSet}}. \end{tikzcd}\end{equation}is not as well-behaved as the profinite completion of groupoids. In particular, it is rarely the case that $\widehat{X\times Y}\simeq \widehat{X}\times \widehat{Y}$. In \cite{bhr1}, we establish a criteria which allows us to profinitely complete a small family of modular operads. 

\begin{definition}\label{def: good}
A discrete group $G$ is said to be \emph{good} if for any finite abelian group $M$ equipped with a $G$-action the map $G\rightarrow \widehat{G}$ induces an isomorphism $$H^{i}(\widehat{G},M) \rightarrow H^{i}(G,M).$$
\end{definition} 

\begin{prop}\label{prof spaces}\cite[Proposition 3.9]{bhr1}
Let $X$ and $Y$ be two connected spaces whose homotopy groups are good. Then the map \[\begin{tikzcd} \widehat{X\times Y} \arrow[r] & \widehat{X} \times \widehat{Y}\end{tikzcd}\] is a weak equivalence of profinite spaces. 
\end{prop}

In the case that every space of $B\bP=\{B\bP(n)\}$ satisfies the conditions of the proposition then 
\[\begin{tikzcd} N\widehat{B\bP}_G \rightarrow \prod_{v\in V} N\widehat{B\bP}_{\medwhitestar_v} \end{tikzcd}\] is a weak equivalence for all graphs and $N\widehat{B\bP}$ is a modular $\infty$-operad (\cite[Proposition 5.1]{bhr1}; \cite{BR22}). 

\begin{remark}
In a recent paper by Blom and Moerdijk \cite{blom2021profinite} they provide a more complete characterization of profinte topological operads. In particular, they provide a fibrantly generated model structure on the category of dendroidal
profinite sets which characterizes profinte operads. 
\begin{problem}
Extend the work of Blom and Moerdijk to characterize profinite cyclic and modular operads. 
\end{problem}
\end{remark}

\subsection{Operads and mapping class groups}\label{sec: mapping class}
Let $\Sigma$ be a surface of genus $g$ with $n$-boundary components. We say that such a surface of \emph{type} $(g,n)$. The boundary components of $\Sigma$ will always be equipped with an \emph{ordering} $\rho: \mathbb{Z}/n\mathbb{Z} \to \pi_0(\partial \Sigma)$. Moreover, we require that each boundary, $\partial_i$, be equipped with a \emph{collar}. In other words, for $1\leq i\leq n$ and some fixed $\varepsilon > 0$,  there is an embedding $\phi_i: S^1 \times [0, \varepsilon) \to \Sigma$ onto a neighborhood of $\partial_i$ such that $\phi_i(S^1 \times \{0\}) = \partial_i$.

\begin{definition}
The \emph{mapping class group} of a surface of type $(g,n)$ is the group of isotopy classes of orientation preserving diffeomorphisms which fix collars pointwise: $$\Gamma^{g}_{n} = \pi_0(\textrm{Diff}^+(\Sigma, \partial \Sigma)).$$ 
\end{definition}

\begin{figure}[h!]
\[
\begin{tikzpicture}[x=0.75pt,y=0.75pt,yscale=-0.8,xscale=0.8]

\draw [color={rgb, 255:red, 128; green, 128; blue, 128 }  ,draw opacity=1 ] [dash pattern={on 4.5pt off 4.5pt}]  (195.13,254.79) .. controls (185.51,254.68) and (185.65,193.21) .. (195.03,193.45) ;
\draw    (195.13,254.79) .. controls (204.99,254.89) and (204.89,193.42) .. (195.03,193.45) ;

\draw [color={rgb, 255:red, 4; green, 146; blue, 194 }  ,draw opacity=1 ]   (148.08,225.71) -- (202.43,225.71) ;
\draw [color={rgb, 255:red, 128; green, 128; blue, 128 }  ,draw opacity=1 ] [dash pattern={on 4.5pt off 4.5pt}]  (140.39,256.07) .. controls (130.77,255.95) and (130.92,193.81) .. (140.29,194.06) ;
\draw    (140.39,256.07) .. controls (150.25,256.18) and (150.15,194.03) .. (140.29,194.06) ;

\draw    (124.57,186.33) .. controls (154.88,202.56) and (171.93,204.68) .. (209.18,187.03) ;
\draw    (134.04,190.56) .. controls (158.03,177.15) and (173.82,174.33) .. (200.34,191.97) ;

\draw    (54.58,190.06) .. controls (43.99,100.17) and (278.25,94.24) .. (279.1,190.06) ;
\draw    (279.1,188.81) .. controls (278.82,276.85) and (56.08,283.37) .. (54.58,190.06) ;
\draw [color={rgb, 255:red, 237; green, 130; blue, 14 }  ,draw opacity=1 ] [dash pattern={on 4.5pt off 4.5pt}]  (168.36,257.37) .. controls (158.75,257.27) and (158.89,199.17) .. (168.27,199.4) ;
\draw [color={rgb, 255:red, 237; green, 130; blue, 14 }  ,draw opacity=1 ]   (168.36,257.37) .. controls (178.22,257.47) and (178.13,199.37) .. (168.27,199.4) ;

\draw [color={rgb, 255:red, 4; green, 146; blue, 194 }  ,draw opacity=1 ] [dash pattern={on 4.5pt off 4.5pt}]  (482.09,198.65) .. controls (493.73,194.94) and (485.44,255.52) .. (498.39,255.66) ;
\draw [color={rgb, 255:red, 4; green, 146; blue, 194 }  ,draw opacity=1 ]   (522.51,226.45) .. controls (508.87,225.92) and (514.69,250.7) .. (498.39,255.66) ;
\draw [color={rgb, 255:red, 4; green, 146; blue, 194 }  ,draw opacity=1 ]   (468.41,226.53) .. controls (476.34,226.65) and (474.69,198.95) .. (482.09,198.65) ;
\draw    (299.89,190.22) -- (346.64,190.22) ;
\draw [shift={(348.64,190.22)}, rotate = 180] [color={rgb, 255:red, 0; green, 0; blue, 0 }  ][line width=0.75]    (10.93,-3.29) .. controls (6.95,-1.4) and (3.31,-0.3) .. (0,0) .. controls (3.31,0.3) and (6.95,1.4) .. (10.93,3.29)   ;

\draw [color={rgb, 255:red, 128; green, 128; blue, 128 }  ,draw opacity=1 ] [dash pattern={on 4.5pt off 4.5pt}]  (515.94,254.17) .. controls (506.32,254.06) and (506.47,192.59) .. (515.84,192.83) ;
\draw    (515.94,254.17) .. controls (525.8,254.28) and (525.7,192.81) .. (515.84,192.83) ;

\draw [color={rgb, 255:red, 128; green, 128; blue, 128 }  ,draw opacity=1 ] [dash pattern={on 4.5pt off 4.5pt}]  (461.2,255.45) .. controls (451.58,255.34) and (451.73,193.19) .. (461.1,193.44) ;
\draw    (461.2,255.45) .. controls (471.06,255.56) and (470.96,193.41) .. (461.1,193.44) ;

\draw    (445.38,185.71) .. controls (475.69,201.94) and (492.74,204.06) .. (529.99,186.41) ;
\draw    (454.85,189.94) .. controls (478.84,176.53) and (494.63,173.71) .. (521.15,191.35) ;

\draw    (375.39,189.44) .. controls (364.8,99.55) and (599.06,93.62) .. (599.91,189.44) ;
\draw    (599.91,188.19) .. controls (599.63,276.23) and (376.89,282.75) .. (375.39,189.44) ;

\draw [color={rgb, 255:red, 4; green, 146; blue, 194 }  ,draw opacity=1 ]   (148.08,225.71) .. controls (138.63,226.43) and (123.52,226.43) .. (115.33,220) .. controls (107.15,213.58) and (101.04,196.11) .. (101.04,188.44) .. controls (101.04,180.77) and (107,163.23) .. (121.32,157.75) .. controls (135.63,152.27) and (197.67,152.27) .. (209.6,157.75) .. controls (221.53,163.23) and (228.69,174.19) .. (228.69,189.54) .. controls (228.69,204.88) and (227.96,205.03) .. (222.14,214.32) .. controls (216.32,223.62) and (213.76,225) .. (202.43,225.71) ;
\draw [color={rgb, 255:red, 4; green, 146; blue, 194 }  ,draw opacity=1 ]   (468.17,226.45) .. controls (458.72,227.16) and (443.61,227.16) .. (435.42,220.74) .. controls (427.23,214.32) and (421.12,196.85) .. (421.12,189.18) .. controls (421.12,181.5) and (427.09,163.97) .. (441.41,158.48) .. controls (455.72,153) and (517.76,153) .. (529.69,158.48) .. controls (541.62,163.97) and (548.78,174.93) .. (548.78,190.27) .. controls (548.78,205.62) and (548.05,205.76) .. (542.23,215.06) .. controls (536.41,224.35) and (533.85,225.74) .. (522.51,226.45) ;

\draw (313.06,164.44) node [anchor=north west][inner sep=0.75pt]    {$a$};
\draw (163.06,267.56) node [anchor=north west][inner sep=0.75pt]    {$\alpha $};
\end{tikzpicture}
\]
    \caption{A Dehn twist around a curve $\alpha$}
    \label{twist}
\end{figure}
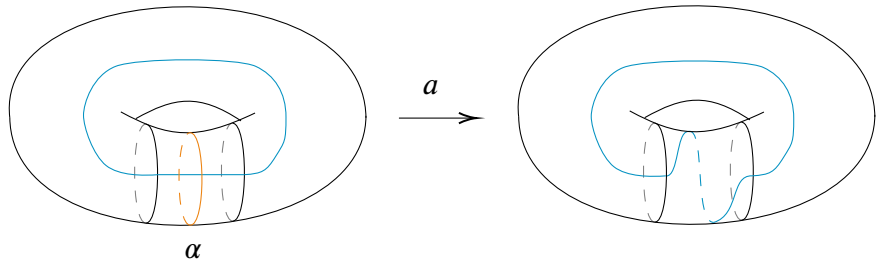

A theorem of Hatcher and Thurston (\cite{HT}) shows that the mapping class group $\Gamma_{g}^{n}$ has a finite presentation: 
$$\Gamma_{g}^{n} = \left< b_1,\ldots, b_n, a_1,\ldots, a_k \mid (C), (B), (D), (L)\right>.$$ The generators of $\Gamma_{g}^{n}$ are Dehn twists along a chosen set of simple closed curves on a surface of type $(g,n)$. Given a curve (i.e. embedded circle) $\alpha$ on a surface, a \emph{Dehn twist} is a diffeomorphism which acts on a neighborhood of $\alpha$, $N_{\alpha}:=S^1 \times [0,1]$, by a full twist: $$a(\theta, t) = (\theta + 2\pi t, t).$$  See Figure~\ref{twist} for an example.   For the purposes of these lecture notes, we will write a Roman letter $a$ to represent a Dehn twist around $\alpha$ as an element of the mapping class group. 

For a surface of type $(g,n)$, the generators of the mapping class consist of Dehn twists around each boundary, $b_1,\ldots, b_n$, as well as Dehn twists, $a_1,\ldots, a_k$, for each curve in a \emph{pants decomposition} of the surface. 

\begin{definition}\label{def:pants decomposition}
Let $\Sigma$ be a surface of type $(g,n)$. A \emph{pants decomposition} of $\Sigma$ is a finite collection of disjoint simple curves (modulo isotopy) which cut $\Sigma$ into surfaces of type $(0,3)$.
\end{definition}

An example of a pants decomposition is depicted in Figure~\ref{fig:pants decomposition} where we have depicted the curves of the decomposition in blue. 

\begin{remark}
Though they play a crucial role in the proofs of some of the theorems we state below, we will not explicitly use the relations from the presentation of the mapping class group in this lecture.  We refer the curious reader to \cite{HT} or \cite[Theorem 1]{hls} for a full description. 
\end{remark}

Shortly, we will describe how the Grothendieck-Teichm\"uller and Nakumara-Schneps groups act on the (profinite) mapping class groups. To do this, it is useful to add slightly more structure to the pants decomposition of a given surface.   For a surface $\Sigma$ of type $(g,n)$, a \emph{quilt} is a choice of two distinct points on each curve and boundary component of a pants decomposition of $\Sigma$, together with a set of disjoint lines between these points. The disjoint lines, or \emph{seams}, cut each pair of pants in a decomposition into two hexagonal patches. In Figure~\ref{quilt} the curves providing the pants decomposition are in blue and the seams of the quilt are in orange. A \emph{quilted pants decomposition} of a surface $\Sigma$ is a pants decomposition of $\Sigma$ together with a choice of quilt.

\begin{figure}[h!]
\begin{subfigure}[b]{0.3\textwidth}
\[
\begin{tikzpicture}[x=0.75pt,y=0.75pt,yscale=-0.8,xscale=0.8]
\draw   (21.61,58.02) .. controls (19.08,55.83) and (24.58,45.33) .. (33.89,34.56) .. controls (43.2,23.79) and (52.8,16.83) .. (55.33,19.02) .. controls (57.86,21.21) and (52.36,31.71) .. (43.05,42.48) .. controls (33.74,53.25) and (24.14,60.21) .. (21.61,58.02) -- cycle ;
\draw   (50.04,155.17) .. controls (47.22,156.96) and (38.75,148.66) .. (31.13,136.63) .. controls (23.51,124.6) and (19.63,113.4) .. (22.45,111.61) .. controls (25.28,109.82) and (33.74,118.12) .. (41.36,130.15) .. controls (48.98,142.18) and (52.87,153.38) .. (50.04,155.17) -- cycle ;
\draw   (185.01,33.39) .. controls (184.99,30.05) and (196.52,27.27) .. (210.76,27.19) .. controls (224.99,27.11) and (236.55,29.75) .. (236.57,33.09) .. controls (236.59,36.44) and (225.06,39.21) .. (210.83,39.29) .. controls (196.59,39.38) and (185.03,36.73) .. (185.01,33.39) -- cycle ;
\draw  [color={rgb, 255:red, 4; green, 146; blue, 194 }  ,draw opacity=1 ] (283.15,88.17) .. controls (283.1,77.78) and (299.43,69.27) .. (319.64,69.17) .. controls (339.85,69.06) and (356.27,77.39) .. (356.33,87.77) .. controls (356.38,98.16) and (340.05,106.66) .. (319.84,106.77) .. controls (299.64,106.88) and (283.21,98.55) .. (283.15,88.17) -- cycle ;
\draw   (238.05,151.14) .. controls (238.06,154.48) and (226.52,157.22) .. (212.29,157.24) .. controls (198.05,157.27) and (186.5,154.58) .. (186.49,151.24) .. controls (186.49,147.9) and (198.03,145.16) .. (212.26,145.14) .. controls (226.5,145.11) and (238.05,147.8) .. (238.05,151.14) -- cycle ;
\draw    (21.61,58.02) .. controls (39.53,61.98) and (46.94,97.56) .. (22.45,111.61) ;
\draw    (186.49,151.24) .. controls (188.72,138.64) and (91.62,136.42) .. (50.04,155.17) ;
\draw    (265.06,110.47) .. controls (253.12,108.16) and (230.4,135.83) .. (238.05,151.14) ;
\draw    (236.57,33.09) .. controls (233.19,49.7) and (255.43,64.52) .. (265.84,67.75) ;
\draw    (185.01,33.39) .. controls (177.6,41.54) and (177.82,45.01) .. (136.94,44.19) .. controls (96.06,43.37) and (65.22,28.64) .. (55.33,19.02) ;
\draw    (79.05,87.6) .. controls (99.34,98.46) and (110.75,99.88) .. (135.69,88.07) ;
\draw    (85.39,90.43) .. controls (101.46,81.46) and (112.02,79.57) .. (129.77,91.38) ;

\draw [color={rgb, 255:red, 4; green, 146; blue, 194 }  ,draw opacity=1 ] [dash pattern={on 4.5pt off 4.5pt}]  (107.48,83.43) .. controls (102.47,83.39) and (102.16,41.45) .. (107.04,41.58) ;
\draw [color={rgb, 255:red, 4; green, 146; blue, 194 }  ,draw opacity=1 ]   (107.48,83.43) .. controls (112.62,83.47) and (112.18,41.53) .. (107.04,41.58) ;

\draw [color={rgb, 255:red, 4; green, 146; blue, 194 }  ,draw opacity=1 ] [dash pattern={on 4.5pt off 4.5pt}]  (36.07,96.4) .. controls (35.78,90.92) and (87.78,87.13) .. (87.96,92.48) ;
\draw [color={rgb, 255:red, 4; green, 146; blue, 194 }  ,draw opacity=1 ]   (36.07,96.4) .. controls (36.38,102.01) and (88.38,98.09) .. (87.96,92.48) ;

\draw [color={rgb, 255:red, 4; green, 146; blue, 194 }  ,draw opacity=1 ] [dash pattern={on 4.5pt off 4.5pt}]  (108.97,142.71) .. controls (102.95,142.67) and (102.62,96.3) .. (108.49,96.44) ;
\draw [color={rgb, 255:red, 4; green, 146; blue, 194 }  ,draw opacity=1 ]   (108.97,142.71) .. controls (115.15,142.75) and (114.66,96.38) .. (108.49,96.44) ;

\draw [color={rgb, 255:red, 4; green, 146; blue, 194 }  ,draw opacity=1 ] [dash pattern={on 4.5pt off 4.5pt}]  (152.14,141.6) .. controls (146.62,141.44) and (145.73,43.44) .. (151.12,43.81) ;
\draw [color={rgb, 255:red, 4; green, 146; blue, 194 }  ,draw opacity=1 ]   (152.14,141.6) .. controls (157.81,141.76) and (156.78,43.75) .. (151.12,43.81) ;

\draw [color={rgb, 255:red, 4; green, 146; blue, 194 }  ,draw opacity=1 ] [dash pattern={on 4.5pt off 4.5pt}]  (264.62,111.22) .. controls (259.61,111.18) and (259.29,68.13) .. (264.17,68.27) ;
\draw [color={rgb, 255:red, 4; green, 146; blue, 194 }  ,draw opacity=1 ]   (264.62,111.22) .. controls (269.76,111.26) and (269.31,68.21) .. (264.17,68.27) ;
\draw    (265.84,67.75) .. controls (291.74,71.93) and (297.67,48.21) .. (318.43,45.99) ;
\draw    (265.06,110.47) .. controls (291,110.47) and (301.38,126.78) .. (322.87,124.56) ;
\draw    (318.43,45.99) .. controls (381.43,44.51) and (402.18,119.37) .. (322.87,124.56) ;
\draw    (291.78,85.37) .. controls (312.06,96.24) and (323.48,97.66) .. (348.41,85.85) ;
\draw    (298.12,88.21) .. controls (314.18,79.23) and (324.74,77.34) .. (342.5,89.15) ;

\draw [color={rgb, 255:red, 4; green, 146; blue, 194 }  ,draw opacity=1 ]   (173.89,42.28) .. controls (170.4,59.44) and (190.94,71.19) .. (190.94,91.94) .. controls (190.94,112.7) and (173.15,123.82) .. (176.12,144.57) ;
\draw [color={rgb, 255:red, 4; green, 146; blue, 194 }  ,draw opacity=1 ]   (243.57,52.66) .. controls (241.34,77.86) and (233.93,74.16) .. (233.93,94.91) .. controls (233.93,115.66) and (241.34,117.89) .. (238.38,135.67) ;
\draw [color={rgb, 255:red, 4; green, 146; blue, 194 }  ,draw opacity=1 ] [dash pattern={on 4.5pt off 4.5pt}]  (243.57,52.66) .. controls (228.74,43.77) and (226.52,51.92) .. (207.99,51.92) .. controls (189.46,51.92) and (192.42,42.28) .. (173.89,42.28) ;
\draw [color={rgb, 255:red, 4; green, 146; blue, 194 }  ,draw opacity=1 ] [dash pattern={on 4.5pt off 4.5pt}]  (238.38,135.67) .. controls (224.3,136.42) and (227.26,130.49) .. (208.73,131.23) .. controls (190.2,131.97) and (196.87,143.09) .. (176.12,144.57) ;
\end{tikzpicture}
\]
    \caption{An example of a pants decomposition.}
    \label{fig:pants decomposition}
\end{subfigure}
  \hfill
\begin{subfigure}[b]{0.3\textwidth}
\[\begin{tikzpicture}[x=0.75pt,y=0.75pt,yscale=-0.5,xscale=0.5]

\draw   (93.55,223.79) .. controls (83.86,223.67) and (76.36,192.21) .. (76.8,153.53) .. controls (77.25,114.85) and (85.47,83.59) .. (95.16,83.72) .. controls (104.86,83.84) and (112.36,115.3) .. (111.91,153.98) .. controls (111.47,192.66) and (103.25,223.92) .. (93.55,223.79) -- cycle ;
\draw   (489.03,223.79) .. controls (479.33,223.67) and (471.83,192.21) .. (472.28,153.53) .. controls (472.72,114.85) and (480.94,83.59) .. (490.64,83.72) .. controls (500.34,83.84) and (507.83,115.3) .. (507.39,153.98) .. controls (506.94,192.66) and (498.72,223.92) .. (489.03,223.79) -- cycle ;
\draw    (95.17,83.72) .. controls (197.96,73.86) and (200.04,28.21) .. (289.54,28.21) .. controls (379.04,28.21) and (383.21,71.68) .. (490.64,83.72) ;
\draw    (489.03,223.79) .. controls (381.13,239.04) and (381.12,278.89) .. (291.62,278.16) .. controls (202.12,277.42) and (195.88,239.04) .. (93.55,223.79) ;
\draw    (224.13,149.39) .. controls (273.28,174.44) and (300.93,177.7) .. (361.35,150.48) ;
\draw    (232.18,152.79) .. controls (271.1,132.11) and (308.53,127.54) .. (351.54,154.76) ;

\draw [color={rgb, 255:red, 4; green, 146; blue, 194 }  ,draw opacity=1 ] [dash pattern={on 4.5pt off 4.5pt}]  (289.8,279.29) .. controls (273.6,279.22) and (272.88,170.71) .. (288.67,171.03) ;
\draw [color={rgb, 255:red, 4; green, 146; blue, 194 }  ,draw opacity=1 ]   (289.8,279.29) .. controls (306.4,279.35) and (305.28,170.85) .. (288.67,171.03) ;

\draw [color={rgb, 255:red, 4; green, 146; blue, 194 }  ,draw opacity=1 ] [dash pattern={on 4.5pt off 4.5pt}]  (289.8,136.46) .. controls (273.6,136.39) and (272.88,27.88) .. (288.67,28.2) ;
\draw [color={rgb, 255:red, 4; green, 146; blue, 194 }  ,draw opacity=1 ]   (289.8,136.46) .. controls (306.4,136.52) and (305.28,28.02) .. (288.67,28.2) ;

\draw [color={rgb, 255:red, 237; green, 130; blue, 14 }  ,draw opacity=1 ] [dash pattern={on 4.5pt off 4.5pt}]  (454.5,230) .. controls (441.5,229.96) and (416,177.8) .. (415,147.8) .. controls (414,117.8) and (393,45.8) .. (375,47.8) ;
\draw [color={rgb, 255:red, 237; green, 130; blue, 14 }  ,draw opacity=1 ]   (110.21,120.95) .. controls (136.86,121.03) and (178,92.8) .. (203,79.8) .. controls (228,66.8) and (264,55.8) .. (301.21,60.95) ;
\draw  [color={rgb, 255:red, 237; green, 130; blue, 14 }  ,draw opacity=1 ][fill={rgb, 255:red, 237; green, 130; blue, 14 }  ,fill opacity=1 ] (107.19,120.95) .. controls (107.19,119.38) and (108.54,118.11) .. (110.21,118.11) .. controls (111.87,118.11) and (113.22,119.38) .. (113.22,120.95) .. controls (113.22,122.53) and (111.87,123.8) .. (110.21,123.8) .. controls (108.54,123.8) and (107.19,122.53) .. (107.19,120.95) -- cycle ;
\draw  [color={rgb, 255:red, 237; green, 130; blue, 14 }  ,draw opacity=1 ][fill={rgb, 255:red, 237; green, 130; blue, 14 }  ,fill opacity=1 ] (106.19,188.95) .. controls (106.19,187.38) and (107.54,186.11) .. (109.21,186.11) .. controls (110.87,186.11) and (112.22,187.38) .. (112.22,188.95) .. controls (112.22,190.53) and (110.87,191.8) .. (109.21,191.8) .. controls (107.54,191.8) and (106.19,190.53) .. (106.19,188.95) -- cycle ;
\draw  [color={rgb, 255:red, 237; green, 130; blue, 14 }  ,draw opacity=1 ][fill={rgb, 255:red, 237; green, 130; blue, 14 }  ,fill opacity=1 ] (298.19,60.95) .. controls (298.19,59.38) and (299.54,58.11) .. (301.21,58.11) .. controls (302.87,58.11) and (304.22,59.38) .. (304.22,60.95) .. controls (304.22,62.53) and (302.87,63.8) .. (301.21,63.8) .. controls (299.54,63.8) and (298.19,62.53) .. (298.19,60.95) -- cycle ;
\draw  [color={rgb, 255:red, 237; green, 130; blue, 14 }  ,draw opacity=1 ][fill={rgb, 255:red, 237; green, 130; blue, 14 }  ,fill opacity=1 ] (298.19,98.95) .. controls (298.19,97.38) and (299.54,96.11) .. (301.21,96.11) .. controls (302.87,96.11) and (304.22,97.38) .. (304.22,98.95) .. controls (304.22,100.53) and (302.87,101.8) .. (301.21,101.8) .. controls (299.54,101.8) and (298.19,100.53) .. (298.19,98.95) -- cycle ;
\draw  [color={rgb, 255:red, 237; green, 130; blue, 14 }  ,draw opacity=1 ][fill={rgb, 255:red, 237; green, 130; blue, 14 }  ,fill opacity=1 ] (471.19,189.95) .. controls (471.19,188.38) and (472.54,187.11) .. (474.21,187.11) .. controls (475.87,187.11) and (477.22,188.38) .. (477.22,189.95) .. controls (477.22,191.53) and (475.87,192.8) .. (474.21,192.8) .. controls (472.54,192.8) and (471.19,191.53) .. (471.19,189.95) -- cycle ;
\draw  [color={rgb, 255:red, 237; green, 130; blue, 14 }  ,draw opacity=1 ][fill={rgb, 255:red, 237; green, 130; blue, 14 }  ,fill opacity=1 ] (298.19,249.95) .. controls (298.19,248.38) and (299.54,247.11) .. (301.21,247.11) .. controls (302.87,247.11) and (304.22,248.38) .. (304.22,249.95) .. controls (304.22,251.53) and (302.87,252.8) .. (301.21,252.8) .. controls (299.54,252.8) and (298.19,251.53) .. (298.19,249.95) -- cycle ;
\draw  [color={rgb, 255:red, 237; green, 130; blue, 14 }  ,draw opacity=1 ][fill={rgb, 255:red, 237; green, 130; blue, 14 }  ,fill opacity=1 ] (298.19,209.95) .. controls (298.19,208.38) and (299.54,207.11) .. (301.21,207.11) .. controls (302.87,207.11) and (304.22,208.38) .. (304.22,209.95) .. controls (304.22,211.53) and (302.87,212.8) .. (301.21,212.8) .. controls (299.54,212.8) and (298.19,211.53) .. (298.19,209.95) -- cycle ;
\draw  [color={rgb, 255:red, 237; green, 130; blue, 14 }  ,draw opacity=1 ][fill={rgb, 255:red, 237; green, 130; blue, 14 }  ,fill opacity=1 ] (471.19,128.95) .. controls (471.19,127.38) and (472.54,126.11) .. (474.21,126.11) .. controls (475.87,126.11) and (477.22,127.38) .. (477.22,128.95) .. controls (477.22,130.53) and (475.87,131.8) .. (474.21,131.8) .. controls (472.54,131.8) and (471.19,130.53) .. (471.19,128.95) -- cycle ;
\draw [color={rgb, 255:red, 237; green, 130; blue, 14 }  ,draw opacity=1 ]   (301.21,249.95) .. controls (274,253.8) and (237,239.8) .. (205,222.8) .. controls (173,205.8) and (147,190.8) .. (109.21,188.95) ;
\draw [color={rgb, 255:red, 237; green, 130; blue, 14 }  ,draw opacity=1 ]   (301.21,209.95) .. controls (264.8,204) and (206.2,191) .. (206.5,150) .. controls (206.8,109) and (262.8,102) .. (301.21,98.95) ;
\draw [color={rgb, 255:red, 237; green, 130; blue, 14 }  ,draw opacity=1 ]   (301.21,60.95) .. controls (343,61) and (348,40) .. (375,47.8) ;
\draw [color={rgb, 255:red, 237; green, 130; blue, 14 }  ,draw opacity=1 ]   (454.5,230) .. controls (467.82,230.04) and (467.5,212) .. (474.21,189.95) ;
\draw [color={rgb, 255:red, 237; green, 130; blue, 14 }  ,draw opacity=1 ]   (301.21,249.95) .. controls (343,250) and (447.21,121.15) .. (474.21,128.95) ;
\draw [color={rgb, 255:red, 237; green, 130; blue, 14 }  ,draw opacity=1 ]   (301.87,98.96) .. controls (338.24,105.12) and (381.5,111) .. (381.5,151) .. controls (381.5,191) and (339.63,207.13) .. (301.21,209.95) ;
\end{tikzpicture}\]

\caption{An example of a quilted pants decomposition.}\label{quilt}
\end{subfigure}
\end{figure}

\subsubsection{The modular operad of quilted surfaces}
The reader may have noticed by now, that a pants decomposition of a surface $\Sigma$ is equivalent to giving a prescription for building $\Sigma$ via modular operad operations such as those pictured in Figure~\ref{comp and contract}.  The second author and L. Bonatto have modified the surface operad of Tillmann \cite{Till} and Wahl \cite{W2} to define a modular operad in groupoids whose objects are surfaces of type $(g,n)$ built from a \emph{standard quilted pair of pants}. Our choice of the standard quilted pair of pants is pictured on the left hand side of Figure~\ref{half_Dehn}. 

\begin{definition}\label{def:groupoids} 
Define a collection of genus graded groupoids $\bS(g,n)$:
\begin{itemize} 
\item $\bS(0,0)=\bS(1,0)=\emptyset$;
\item $\bS(0,2)$ is the groupoid whose only object is homotopic to the circle $S^1$ (the ``thin cyclinder'') and whose morphisms are given by the integers.\footnote{The groupoid $\bS(0,2)$ is equivalent to the groupoid $\calS_{0,1,1}$ in \cite[3.1.1]{W2}.} 
\item For all other $g,n\geq 0$, $\bS(g,n)$ is the groupoid whose objects are surfaces of type $(g,n)$ which are built from the the standard pair of pants. Morphisms are isotopy classes of orientation preserving diffeomorphisms which fix the boundary collars pointwise. 
\end{itemize} 
\end{definition} 

The symmetric group $\Sigma_n$ acts freely on the groupoid $\bS(g,n)$ by permuting the labels of boundaries. Moreover, we can define composition and contraction functors
 \[\begin{tikzcd}
    \bS(g,n)\times \bS(h,m)\arrow[r, "\circ_{ij}"] & \bS(g+h,n+k-2)
    \end{tikzcd}\]
\[\begin{tikzcd}
\bS(g,n) \arrow[r, "\xi_{ij}"] &\bS(g+1,n-2)
\end{tikzcd}\] at the level of objects by gluing of surfaces. On morphisms these functors act by inclusion of Dehn twists.  The following proposition appears in \cite{BR22}:

\begin{prop}
The composition and contraction operations are well-defined, associative, equivariant and unital operations and thus $$\bS =\{\bS(g,n)\}$$ assembles into a genus graded modular operad in groupoids. 
\end{prop}

The groupoids $\bS(g,n)$ are homotopy approximations of the mapping class groups in the sense that:$$B\bS(g,n)\simeq B\Gamma^{g}_{n}.$$  Aside from this point, the usefulness of the $\bS(g,n)$ is that they are \emph{finite}\footnote{We are using the term finite here in the sense that our groupoids $\bS(g,n)$ all have finitely many objects.} groupoid approximations of a contractible $2$-dimensional simplicial complex called the \underline{S}eamed \underline{H}atcher-\underline{T}hurston complex, denoted by $\mathcal{SHT}$ in \cite{ns} (see also the related $\mathcal{HT}$ complex in \cite{hls},\cite{HT}). In short, a finite number of the points of $\mathcal{SHT}$ are objects in our groupoids. Morphisms in $\bS(g,n)$ are then $1$-cells (or composites of $1$-cells) from the complex $\mathcal{SHT}$.  A consequence of this fact is that morphisms in $\bS(g,n)$ are generated by three types of elementary morphisms: half-Dehn twists, $A$-moves (``associative moves''), and $S$-moves (``simple moves''). Examples of these diffeomorphisms are depicted in Figure~\ref{morphisms}.

\begin{figure}[h!]
\begin{subfigure}[b]{0.4\textwidth}
\[
\begin{tikzpicture}[x=0.75pt,y=0.75pt,yscale=-.75,xscale=.75]

\draw   (15.51,222.14) .. controls (13.19,222.11) and (11.39,214.89) .. (11.5,206.02) .. controls (11.6,197.14) and (13.57,189.97) .. (15.9,189.99) .. controls (18.22,190.02) and (20.02,197.24) .. (19.91,206.12) .. controls (19.8,215) and (17.83,222.17) .. (15.51,222.14) -- cycle ;
\draw    (15.9,190) .. controls (31.49,190.08) and (56.88,165.71) .. (85.79,165.71) .. controls (114.69,165.71) and (138.29,188.31) .. (137.7,207.92) .. controls (137.11,227.54) and (111.74,247.76) .. (85.2,247.16) .. controls (58.65,246.57) and (28.06,221.89) .. (15.51,222.14) ;
\draw    (59.77,203.67) .. controls (75.92,212.39) and (85,213.53) .. (104.85,204.05) ;
\draw    (64.82,205.95) .. controls (77.6,198.75) and (86.01,197.23) .. (100.14,206.71) ;

\draw [color={rgb, 255:red, 4; green, 146; blue, 194 }  ,draw opacity=1 ] [dash pattern={on 4.5pt off 4.5pt}]  (83.37,247.38) .. controls (77.93,247.35) and (77.68,210.41) .. (82.99,210.52) ;
\draw [color={rgb, 255:red, 4; green, 146; blue, 194 }  ,draw opacity=1 ]   (83.37,247.38) .. controls (88.95,247.4) and (88.57,210.46) .. (82.99,210.52) ;

\draw  [color={rgb, 255:red, 4; green, 146; blue, 194 }  ,draw opacity=1 ] (250.8,205.8) .. controls (250.75,195.8) and (267.27,187.6) .. (287.7,187.49) .. controls (308.14,187.38) and (324.74,195.4) .. (324.8,205.4) .. controls (324.85,215.4) and (308.33,223.6) .. (287.9,223.71) .. controls (267.46,223.82) and (250.85,215.8) .. (250.8,205.8) -- cycle ;
\draw   (220.8,222.14) .. controls (218.48,222.11) and (216.68,214.89) .. (216.79,206.02) .. controls (216.9,197.14) and (218.87,189.97) .. (221.19,189.99) .. controls (223.51,190.02) and (225.31,197.24) .. (225.2,206.12) .. controls (225.09,215) and (223.12,222.17) .. (220.8,222.14) -- cycle ;
\draw    (221.19,190) .. controls (236.78,190.08) and (262.17,165.71) .. (291.08,165.71) .. controls (319.98,165.71) and (343.58,188.31) .. (342.99,207.92) .. controls (342.4,227.54) and (317.03,247.76) .. (290.49,247.16) .. controls (263.94,246.57) and (233.35,221.89) .. (220.8,222.14) ;
\draw    (265.06,203.67) .. controls (281.21,212.39) and (290.29,213.53) .. (310.14,204.05) ;
\draw    (270.11,205.95) .. controls (282.89,198.75) and (291.3,197.23) .. (305.43,206.71) ;

\draw    (153.93,206.61) -- (197.98,206.81) ;
\draw [shift={(199.98,206.82)}, rotate = 180.26] [color={rgb, 255:red, 0; green, 0; blue, 0 }  ][line width=0.75]    (10.93,-3.29) .. controls (6.95,-1.4) and (3.31,-0.3) .. (0,0) .. controls (3.31,0.3) and (6.95,1.4) .. (10.93,3.29)   ;
\draw   (216.82,308.89) .. controls (215.02,307.41) and (218.07,300.64) .. (223.63,293.75) .. controls (229.2,286.87) and (235.17,282.49) .. (236.97,283.97) .. controls (238.77,285.45) and (235.72,292.23) .. (230.15,299.11) .. controls (224.59,305.99) and (218.62,310.37) .. (216.82,308.89) -- cycle ;
\draw   (235.51,376.54) .. controls (233.56,377.81) and (228.12,372.78) .. (223.36,365.31) .. controls (218.59,357.84) and (216.32,350.76) .. (218.27,349.49) .. controls (220.23,348.23) and (225.67,353.26) .. (230.43,360.73) .. controls (235.19,368.2) and (237.47,375.28) .. (235.51,376.54) -- cycle ;
\draw   (285.62,284.68) .. controls (287.23,282.99) and (293.68,286.62) .. (300.02,292.79) .. controls (306.35,298.95) and (310.18,305.32) .. (308.57,307) .. controls (306.96,308.69) and (300.51,305.05) .. (294.17,298.89) .. controls (287.84,292.72) and (284.01,286.36) .. (285.62,284.68) -- cycle ;
\draw   (307.01,347.69) .. controls (309.14,348.62) and (308.01,355.98) .. (304.49,364.12) .. controls (300.97,372.25) and (296.39,378.09) .. (294.26,377.16) .. controls (292.13,376.22) and (293.26,368.87) .. (296.78,360.73) .. controls (300.3,352.59) and (304.88,346.75) .. (307.01,347.69) -- cycle ;
\draw    (216.82,308.89) .. controls (229.14,317.06) and (228.55,339.65) .. (218.28,349.49) ;
\draw    (294.26,377.16) .. controls (283.41,367) and (247.42,365.21) .. (235.51,376.54) ;
\draw    (285.62,284.67) .. controls (276.33,293.87) and (246.24,291.5) .. (236.97,283.97) ;
\draw    (307.01,347.69) .. controls (295.8,343.81) and (297.57,313.49) .. (308.57,307) ;

\draw [color={rgb, 255:red, 4; green, 146; blue, 194 }  ,draw opacity=1 ] [dash pattern={on 4.5pt off 4.5pt}]  (263.19,368.64) .. controls (256.3,368.54) and (255.66,290.1) .. (262.37,290.38) ;
\draw [color={rgb, 255:red, 4; green, 146; blue, 194 }  ,draw opacity=1 ]   (263.19,368.64) .. controls (270.25,368.74) and (269.43,290.31) .. (262.37,290.38) ;

\draw   (40.43,308.3) .. controls (38.63,306.82) and (41.68,300.04) .. (47.25,293.16) .. controls (52.82,286.28) and (58.79,281.9) .. (60.59,283.38) .. controls (62.39,284.86) and (59.34,291.64) .. (53.77,298.52) .. controls (48.2,305.4) and (42.23,309.78) .. (40.43,308.3) -- cycle ;
\draw   (59.13,375.95) .. controls (57.17,377.21) and (51.73,372.18) .. (46.97,364.71) .. controls (42.21,357.24) and (39.94,350.16) .. (41.89,348.9) .. controls (43.84,347.63) and (49.29,352.66) .. (54.05,360.13) .. controls (58.81,367.6) and (61.08,374.68) .. (59.13,375.95) -- cycle ;
\draw   (109.23,284.08) .. controls (110.85,282.4) and (117.29,286.03) .. (123.63,292.19) .. controls (129.97,298.36) and (133.8,304.72) .. (132.19,306.41) .. controls (130.57,308.09) and (124.13,304.46) .. (117.79,298.29) .. controls (111.45,292.13) and (107.62,285.76) .. (109.23,284.08) -- cycle ;
\draw   (130.62,347.09) .. controls (132.75,348.03) and (131.62,355.38) .. (128.1,363.52) .. controls (124.58,371.66) and (120,377.5) .. (117.88,376.56) .. controls (115.75,375.63) and (116.87,368.27) .. (120.39,360.14) .. controls (123.91,352) and (128.49,346.16) .. (130.62,347.09) -- cycle ;
\draw    (40.43,308.3) .. controls (52.75,316.46) and (52.16,339.06) .. (41.89,348.9) ;
\draw    (117.88,376.56) .. controls (107.02,366.4) and (71.04,364.62) .. (59.13,375.95) ;
\draw    (109.23,284.08) .. controls (99.94,293.28) and (69.86,290.9) .. (60.59,283.38) ;
\draw    (130.62,347.09) .. controls (119.41,343.22) and (121.18,312.9) .. (132.18,306.41) ;

\draw [color={rgb, 255:red, 4; green, 146; blue, 194 }  ,draw opacity=1 ] [dash pattern={on 4.5pt off 4.5pt}]  (49.43,327.83) .. controls (49.66,320.89) and (123.35,321.81) .. (122.95,328.57) ;
\draw [color={rgb, 255:red, 4; green, 146; blue, 194 }  ,draw opacity=1 ]   (49.43,327.83) .. controls (49.2,334.95) and (122.88,335.69) .. (122.95,328.57) ;

\draw    (153.93,330.86) -- (197.98,331.06) ;
\draw [shift={(199.98,331.07)}, rotate = 180.26] [color={rgb, 255:red, 0; green, 0; blue, 0 }  ][line width=0.75]    (10.93,-3.29) .. controls (6.95,-1.4) and (3.31,-0.3) .. (0,0) .. controls (3.31,0.3) and (6.95,1.4) .. (10.93,3.29)   ;

\draw (150,178) node [anchor=north west][inner sep=0.75pt]   [align=left][font=\footnotesize] {S-move};
\draw (151,302) node [anchor=north west][inner sep=0.75pt]   [align=left][font=\footnotesize] {A-move};
\end{tikzpicture}\]
    \caption{Elementary morphisms in $\bS(g,n)$}
    \label{AandS}
\end{subfigure}
\hfill
\begin{subfigure}[b]{0.5\textwidth}
\[
\begin{tikzpicture}[x=0.75pt,y=0.75pt,yscale=-0.55,xscale=0.55]

\draw   (85.43,61.85) .. controls (85.43,56.73) and (104.07,52.57) .. (127.07,52.57) .. controls (150.06,52.57) and (168.71,56.73) .. (168.71,61.85) .. controls (168.71,66.98) and (150.06,71.14) .. (127.07,71.14) .. controls (104.07,71.14) and (85.43,66.98) .. (85.43,61.85) -- cycle ;
\draw   (24.42,156.24) .. controls (28.4,153) and (43.4,164.82) .. (57.93,182.64) .. controls (72.46,200.47) and (81.02,217.54) .. (77.04,220.78) .. controls (73.07,224.02) and (58.07,212.2) .. (43.54,194.37) .. controls (29.01,176.55) and (20.45,159.48) .. (24.42,156.24) -- cycle ;
\draw   (181.97,220.81) .. controls (178.15,217.39) and (187.49,200.73) .. (202.84,183.6) .. controls (218.19,166.48) and (233.72,155.37) .. (237.54,158.79) .. controls (241.36,162.21) and (232.01,178.87) .. (216.67,195.99) .. controls (201.32,213.12) and (185.79,224.23) .. (181.97,220.81) -- cycle ;
\draw    (24.42,156.24) .. controls (44.06,145.29) and (85.1,86.77) .. (85.43,61.85) ;
\draw    (77.04,220.78) .. controls (101.82,208.37) and (161.11,209.89) .. (181.97,220.81) ;
\draw    (237.54,158.79) .. controls (205.19,141.49) and (165.67,98.17) .. (168.71,61.85) ;
\draw  [color={rgb, 255:red, 237; green, 130; blue, 14 }  ,draw opacity=1 ][fill={rgb, 255:red, 237; green, 130; blue, 14 }  ,fill opacity=1 ] (107.62,71.1) .. controls (107.62,69.91) and (108.64,68.94) .. (109.91,68.94) .. controls (111.17,68.94) and (112.2,69.91) .. (112.2,71.1) .. controls (112.2,72.3) and (111.17,73.27) .. (109.91,73.27) .. controls (108.64,73.27) and (107.62,72.3) .. (107.62,71.1) -- cycle ;
\draw  [color={rgb, 255:red, 237; green, 130; blue, 14 }  ,draw opacity=1 ][fill={rgb, 255:red, 237; green, 130; blue, 14 }  ,fill opacity=1 ] (210.98,172.18) .. controls (210.98,170.99) and (212.01,170.02) .. (213.27,170.02) .. controls (214.53,170.02) and (215.56,170.99) .. (215.56,172.18) .. controls (215.56,173.38) and (214.53,174.35) .. (213.27,174.35) .. controls (212.01,174.35) and (210.98,173.38) .. (210.98,172.18) -- cycle ;
\draw  [color={rgb, 255:red, 237; green, 130; blue, 14 }  ,draw opacity=1 ][fill={rgb, 255:red, 237; green, 130; blue, 14 }  ,fill opacity=1 ] (188.18,198.78) .. controls (188.18,197.59) and (189.2,196.62) .. (190.47,196.62) .. controls (191.73,196.62) and (192.76,197.59) .. (192.76,198.78) .. controls (192.76,199.98) and (191.73,200.95) .. (190.47,200.95) .. controls (189.2,200.95) and (188.18,199.98) .. (188.18,198.78) -- cycle ;
\draw  [color={rgb, 255:red, 237; green, 130; blue, 14 }  ,draw opacity=1 ][fill={rgb, 255:red, 237; green, 130; blue, 14 }  ,fill opacity=1 ] (65.06,194.98) .. controls (65.06,193.79) and (66.08,192.82) .. (67.35,192.82) .. controls (68.61,192.82) and (69.64,193.79) .. (69.64,194.98) .. controls (69.64,196.18) and (68.61,197.15) .. (67.35,197.15) .. controls (66.08,197.15) and (65.06,196.18) .. (65.06,194.98) -- cycle ;
\draw  [color={rgb, 255:red, 237; green, 130; blue, 14 }  ,draw opacity=1 ][fill={rgb, 255:red, 237; green, 130; blue, 14 }  ,fill opacity=1 ] (42.26,168.38) .. controls (42.26,167.19) and (43.28,166.22) .. (44.55,166.22) .. controls (45.81,166.22) and (46.84,167.19) .. (46.84,168.38) .. controls (46.84,169.58) and (45.81,170.55) .. (44.55,170.55) .. controls (43.28,170.55) and (42.26,169.58) .. (42.26,168.38) -- cycle ;
\draw  [color={rgb, 255:red, 237; green, 130; blue, 14 }  ,draw opacity=1 ][fill={rgb, 255:red, 237; green, 130; blue, 14 }  ,fill opacity=1 ] (141.82,71.1) .. controls (141.82,69.91) and (142.84,68.94) .. (144.11,68.94) .. controls (145.37,68.94) and (146.4,69.91) .. (146.4,71.1) .. controls (146.4,72.3) and (145.37,73.27) .. (144.11,73.27) .. controls (142.84,73.27) and (141.82,72.3) .. (141.82,71.1) -- cycle ;
\draw [color={rgb, 255:red, 237; green, 130; blue, 14 }  ,draw opacity=1 ]   (44.55,168.38) .. controls (78.59,150.61) and (106.71,92.85) .. (109.91,71.1) ;
\draw [color={rgb, 255:red, 237; green, 130; blue, 14 }  ,draw opacity=1 ]   (213.27,172.18) .. controls (175.11,157.45) and (140.91,92.85) .. (144.11,71.1) ;
\draw [color={rgb, 255:red, 237; green, 130; blue, 14 }  ,draw opacity=1 ]   (67.35,194.98) .. controls (97.59,177.97) and (159.91,178.73) .. (190.47,198.78) ;
\draw    (265.93,160.08) -- (310.43,160.28) ;
\draw [shift={(312.43,160.29)}, rotate = 180.26] [color={rgb, 255:red, 0; green, 0; blue, 0 }  ][line width=0.75]    (10.93,-3.29) .. controls (6.95,-1.4) and (3.31,-0.3) .. (0,0) .. controls (3.31,0.3) and (6.95,1.4) .. (10.93,3.29)   ;
\draw   (398.43,62.85) .. controls (398.43,57.73) and (417.07,53.57) .. (440.07,53.57) .. controls (463.06,53.57) and (481.71,57.73) .. (481.71,62.85) .. controls (481.71,67.98) and (463.06,72.14) .. (440.07,72.14) .. controls (417.07,72.14) and (398.43,67.98) .. (398.43,62.85) -- cycle ;
\draw   (337.42,157.24) .. controls (341.4,154) and (356.4,165.82) .. (370.93,183.64) .. controls (385.46,201.47) and (394.02,218.54) .. (390.04,221.78) .. controls (386.07,225.02) and (371.07,213.2) .. (356.54,195.37) .. controls (342.01,177.55) and (333.45,160.48) .. (337.42,157.24) -- cycle ;
\draw   (494.97,221.81) .. controls (491.15,218.39) and (500.49,201.73) .. (515.84,184.6) .. controls (531.19,167.48) and (546.72,156.37) .. (550.54,159.79) .. controls (554.36,163.21) and (545.01,179.87) .. (529.67,196.99) .. controls (514.32,214.12) and (498.79,225.23) .. (494.97,221.81) -- cycle ;
\draw    (337.42,157.24) .. controls (357.06,146.29) and (398.1,87.77) .. (398.43,62.85) ;
\draw    (390.04,221.78) .. controls (414.82,209.37) and (474.11,210.89) .. (494.97,221.81) ;
\draw    (550.54,159.79) .. controls (518.19,142.49) and (478.67,99.17) .. (481.71,62.85) ;
\draw  [color={rgb, 255:red, 237; green, 130; blue, 14 }  ,draw opacity=1 ][fill={rgb, 255:red, 237; green, 130; blue, 14 }  ,fill opacity=1 ] (420.62,72.1) .. controls (420.62,70.91) and (421.64,69.94) .. (422.91,69.94) .. controls (424.17,69.94) and (425.2,70.91) .. (425.2,72.1) .. controls (425.2,73.3) and (424.17,74.27) .. (422.91,74.27) .. controls (421.64,74.27) and (420.62,73.3) .. (420.62,72.1) -- cycle ;
\draw  [color={rgb, 255:red, 237; green, 130; blue, 14 }  ,draw opacity=1 ][fill={rgb, 255:red, 237; green, 130; blue, 14 }  ,fill opacity=1 ] (523.98,173.18) .. controls (523.98,171.99) and (525.01,171.02) .. (526.27,171.02) .. controls (527.53,171.02) and (528.56,171.99) .. (528.56,173.18) .. controls (528.56,174.38) and (527.53,175.35) .. (526.27,175.35) .. controls (525.01,175.35) and (523.98,174.38) .. (523.98,173.18) -- cycle ;
\draw  [color={rgb, 255:red, 237; green, 130; blue, 14 }  ,draw opacity=1 ][fill={rgb, 255:red, 237; green, 130; blue, 14 }  ,fill opacity=1 ] (501.18,199.78) .. controls (501.18,198.59) and (502.2,197.62) .. (503.47,197.62) .. controls (504.73,197.62) and (505.76,198.59) .. (505.76,199.78) .. controls (505.76,200.98) and (504.73,201.95) .. (503.47,201.95) .. controls (502.2,201.95) and (501.18,200.98) .. (501.18,199.78) -- cycle ;
\draw  [color={rgb, 255:red, 237; green, 130; blue, 14 }  ,draw opacity=1 ][fill={rgb, 255:red, 237; green, 130; blue, 14 }  ,fill opacity=1 ] (378.06,195.98) .. controls (378.06,194.79) and (379.08,193.82) .. (380.35,193.82) .. controls (381.61,193.82) and (382.64,194.79) .. (382.64,195.98) .. controls (382.64,197.18) and (381.61,198.15) .. (380.35,198.15) .. controls (379.08,198.15) and (378.06,197.18) .. (378.06,195.98) -- cycle ;
\draw  [color={rgb, 255:red, 237; green, 130; blue, 14 }  ,draw opacity=1 ][fill={rgb, 255:red, 237; green, 130; blue, 14 }  ,fill opacity=1 ] (355.26,169.38) .. controls (355.26,168.19) and (356.28,167.22) .. (357.55,167.22) .. controls (358.81,167.22) and (359.84,168.19) .. (359.84,169.38) .. controls (359.84,170.58) and (358.81,171.55) .. (357.55,171.55) .. controls (356.28,171.55) and (355.26,170.58) .. (355.26,169.38) -- cycle ;
\draw  [color={rgb, 255:red, 237; green, 130; blue, 14 }  ,draw opacity=1 ][fill={rgb, 255:red, 237; green, 130; blue, 14 }  ,fill opacity=1 ] (454.82,72.1) .. controls (454.82,70.91) and (455.84,69.94) .. (457.11,69.94) .. controls (458.37,69.94) and (459.4,70.91) .. (459.4,72.1) .. controls (459.4,73.3) and (458.37,74.27) .. (457.11,74.27) .. controls (455.84,74.27) and (454.82,73.3) .. (454.82,72.1) -- cycle ;
\draw [color={rgb, 255:red, 237; green, 130; blue, 14 }  ,draw opacity=1 ]   (357.55,169.38) .. controls (366.43,161) and (366.43,147) .. (363.43,131) ;
\draw [color={rgb, 255:red, 237; green, 130; blue, 14 }  ,draw opacity=1 ]   (526.27,173.18) .. controls (488.11,158.45) and (419.71,93.85) .. (422.91,72.1) ;
\draw [color={rgb, 255:red, 237; green, 130; blue, 14 }  ,draw opacity=1 ]   (380.35,195.98) .. controls (410.59,178.97) and (472.91,179.73) .. (503.47,199.78) ;
\draw [color={rgb, 255:red, 237; green, 130; blue, 14 }  ,draw opacity=1 ] [dash pattern={on 4.5pt off 4.5pt}]  (491.43,101) .. controls (500,115.43) and (453,108.43) .. (429,110.43) .. controls (405,112.43) and (375,118.43) .. (363.43,131) ;
\draw [color={rgb, 255:red, 237; green, 130; blue, 14 }  ,draw opacity=1 ]   (457.11,74.27) .. controls (467,88) and (483,83.43) .. (491.43,101) ;

\draw (274,128.4) node [anchor=north west][inner sep=0.75pt]    {$a^{1/2}$};
\draw (181,52.4) node [anchor=north west][inner sep=0.75pt]    {$\alpha $};
\draw (495,54.4) node [anchor=north west][inner sep=0.75pt]    {$\alpha $};
\end{tikzpicture}\]
    \caption{A half Dehn twist of the standard pair of pants.}
    \label{half_Dehn}
\end{subfigure}
\caption{}\label{morphisms}
\end{figure}

\begin{remark}
The half-Dehn twists are diffeomorphisms which change the quilt on a pair of pants relative to a specific curve. It is beyond the scope of these lectures to explicitly describe the effect of an $A$-move or $S$-move on the quilting of a surface, but the important part is that effect of an $A$-move or $S$-move is strictly defined. In particular, an $A$-move or $S$-move which takes a curve $\alpha$ to a curve $\beta$, changes a quilt by a uniquely determined half-Dehn twist around $\alpha$, $a^{n/2}$, $n\in\mathbb{Z}$. Full details will appear in \cite{BR22}, but see also the notion of a ``quilt adjustment'' in \cite{ns}.  
\end{remark}

The $2$-cells of the $\mathcal{SHT}$ complex can then be used to show that our generating morphisms in $\bS(g,n)$ satisfy the following local relations:

\begin{itemize}[leftmargin=3em]
\item[(3A)] 
 The loops $\beta_1$, $\beta_2$ and $\beta_3$ in Figure~\ref{3A} represent all possible pants decompositions of a surface of type $(0,4)$. If we let $A_{b_i,b_j}$ denote the $A$-move which takes curve $\beta_i \to \beta_j$, then $A_{b_1,b_2}A_{b_2,b_3}A_{b_3,b_1}=1$ in $\bS(0,4)$.
 
\item[(5A)] Similar to above, the loops $\beta_1,\ldots,\beta_5$ in Figure~\ref{fig:5A} pairwise comprise all possible pants decompositions of a surface of type $(0,5)$. Then the equation: 
\begin{equation*}
  A_{b_3,b_4} A_{b_1,b_2}A_{b_4,b_5}A_{b_2,b_3}A_{b_5,b_1} =1
\end{equation*}holds in $\bS(0,5)$
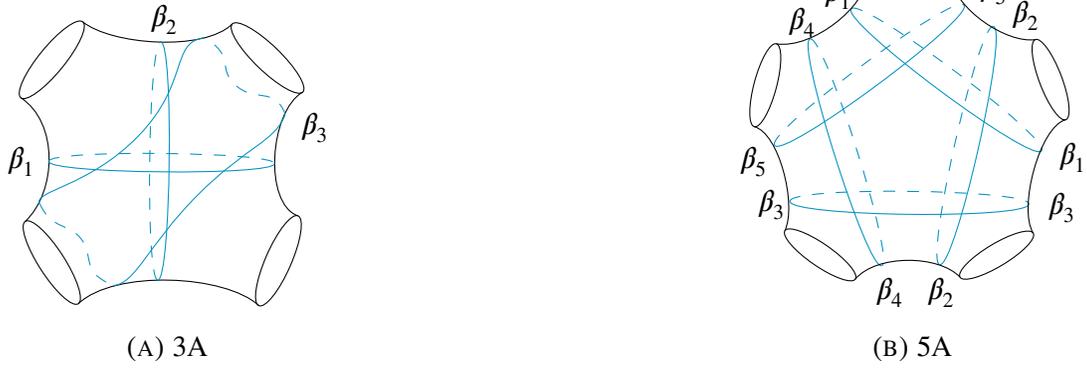
\begin{figure}[h!]
\centering
\begin{subfigure}[b]{0.4\textwidth}
\[\begin{tikzpicture}[x=0.75pt,y=0.75pt,yscale=-0.8,xscale=0.8]

\draw   (58.34,143.12) .. controls (54.85,140.29) and (60.76,127.35) .. (71.55,114.21) .. controls (82.33,101.06) and (93.9,92.7) .. (97.39,95.52) .. controls (100.88,98.35) and (94.96,111.29) .. (84.18,124.44) .. controls (73.39,137.58) and (61.82,145.95) .. (58.34,143.12) -- cycle ;
\draw   (94.56,272.35) .. controls (90.78,274.77) and (80.23,265.16) .. (71.01,250.89) .. controls (61.78,236.62) and (57.37,223.1) .. (61.16,220.68) .. controls (64.95,218.26) and (75.5,227.87) .. (84.72,242.14) .. controls (93.94,256.41) and (98.35,269.93) .. (94.56,272.35) -- cycle ;
\draw   (191.64,96.86) .. controls (194.77,93.64) and (207.26,100.58) .. (219.54,112.36) .. controls (231.82,124.14) and (239.24,136.29) .. (236.11,139.51) .. controls (232.99,142.73) and (220.5,135.79) .. (208.21,124.01) .. controls (195.93,112.24) and (188.51,100.08) .. (191.64,96.86) -- cycle ;
\draw   (233.08,217.23) .. controls (237.21,219.02) and (235.02,233.07) .. (228.2,248.61) .. controls (221.38,264.16) and (212.51,275.31) .. (208.38,273.53) .. controls (204.26,271.74) and (206.45,257.69) .. (213.27,242.15) .. controls (220.09,226.6) and (228.96,215.45) .. (233.08,217.23) -- cycle ;
\draw    (58.33,143.12) .. controls (82.2,158.72) and (81.06,201.88) .. (61.16,220.68) ;
\draw    (208.39,273.53) .. controls (187.36,254.12) and (117.64,250.71) .. (94.56,272.35) ;
\draw    (191.64,96.86) .. controls (173.64,114.43) and (115.35,109.89) .. (97.39,95.52) ;
\draw    (233.08,217.23) .. controls (211.36,209.83) and (214.79,151.91) .. (236.11,139.51) ;

\draw [color={rgb, 255:red, 4; green, 146; blue, 194 }  ,draw opacity=1 ] [dash pattern={on 4.5pt off 4.5pt}]  (75.82,183.03) .. controls (76.21,175.36) and (218.97,176.96) .. (218.26,184.44) ;
\draw [color={rgb, 255:red, 4; green, 146; blue, 194 }  ,draw opacity=1 ]   (75.82,183.03) .. controls (75.42,190.89) and (218.18,192.31) .. (218.26,184.44) ;

\draw [color={rgb, 255:red, 4; green, 146; blue, 194 }  ,draw opacity=1 ] [dash pattern={on 4.5pt off 4.5pt}]  (144.67,257.55) .. controls (136.5,257.11) and (138.59,107.3) .. (146.54,108.06) ;
\draw [color={rgb, 255:red, 4; green, 146; blue, 194 }  ,draw opacity=1 ]   (144.67,257.55) .. controls (153.04,257.99) and (154.92,108.17) .. (146.54,108.06) ;

\draw [color={rgb, 255:red, 4; green, 146; blue, 194 }  ,draw opacity=1 ] [dash pattern={on 4.5pt off 4.5pt}]  (225.06,153.5) .. controls (218.33,135.85) and (206.29,146.2) .. (193.95,135.49) .. controls (181.62,124.77) and (193.95,116.51) .. (173.78,105.58) ;
\draw [color={rgb, 255:red, 4; green, 146; blue, 194 }  ,draw opacity=1 ]   (116.5,260.96) .. controls (130.9,259.9) and (138.47,235.52) .. (169.57,204.42) .. controls (200.68,173.32) and (218.33,171.15) .. (225.06,153.5) ;
\draw [color={rgb, 255:red, 4; green, 146; blue, 194 }  ,draw opacity=1 ]   (173.78,105.58) .. controls (152.76,108.11) and (164.53,134.17) .. (132.58,166.95) .. controls (100.64,199.74) and (65.33,201.42) .. (70.38,209.82) ;
\draw [color={rgb, 255:red, 4; green, 146; blue, 194 }  ,draw opacity=1 ] [dash pattern={on 4.5pt off 4.5pt}]  (70.38,209.82) .. controls (78.36,226.95) and (89.71,221.23) .. (97.28,233) .. controls (104.84,244.77) and (95.58,251.52) .. (116.5,260.96) ;

\draw (48.41,172.12) node [anchor=north west][inner sep=0.75pt]    {$\beta _{1}$};
\draw (138.36,83.01) node [anchor=north west][inner sep=0.75pt]    {$\beta _{2}$};
\draw (233.35,151.11) node [anchor=north west][inner sep=0.75pt]    {$\beta _{3}$};

\end{tikzpicture}\]
    \caption{3A}
    \label{3A}
\end{subfigure}
\hfill
\begin{subfigure}[b]{0.4\textwidth}
\[
\begin{tikzpicture}[x=0.75pt,y=0.75pt,yscale=-0.7,xscale=0.7]
\draw [color={rgb, 255:red, 4; green, 146; blue, 194 }  ,draw opacity=1 ] [dash pattern={on 4.5pt off 4.5pt}]  (140.46,69.11) .. controls (147.55,59.86) and (284.89,162.1) .. (277.66,170.89) ;
\draw [color={rgb, 255:red, 4; green, 146; blue, 194 }  ,draw opacity=1 ]   (140.46,69.11) .. controls (133.19,78.61) and (270.7,180.61) .. (277.66,170.89) ;

\draw [color={rgb, 255:red, 4; green, 146; blue, 194 }  ,draw opacity=1 ] [dash pattern={on 4.5pt off 4.5pt}]  (96.35,207) .. controls (96.86,196.06) and (269.34,198.38) .. (268.45,209.03) ;
\draw [color={rgb, 255:red, 4; green, 146; blue, 194 }  ,draw opacity=1 ]   (96.35,207) .. controls (95.84,218.22) and (268.32,220.25) .. (268.45,209.03) ;

\draw [color={rgb, 255:red, 4; green, 146; blue, 194 }  ,draw opacity=1 ] [dash pattern={on 4.5pt off 4.5pt}]  (110.84,90.65) .. controls (121.3,87.61) and (172.77,250.63) .. (162.47,253.22) ;
\draw [color={rgb, 255:red, 4; green, 146; blue, 194 }  ,draw opacity=1 ]   (110.84,90.65) .. controls (100.13,93.78) and (151.87,256.71) .. (162.47,253.22) ;

\draw [color={rgb, 255:red, 4; green, 146; blue, 194 }  ,draw opacity=1 ] [dash pattern={on 4.5pt off 4.5pt}]  (204.16,252.47) .. controls (193.13,249.48) and (232.23,78.41) .. (242.9,81.71) ;
\draw [color={rgb, 255:red, 4; green, 146; blue, 194 }  ,draw opacity=1 ]   (204.16,252.47) .. controls (215.48,255.53) and (254.3,84.39) .. (242.9,81.71) ;

\draw [color={rgb, 255:red, 4; green, 146; blue, 194 }  ,draw opacity=1 ] [dash pattern={on 4.5pt off 4.5pt}]  (85.22,166.97) .. controls (78.71,157.18) and (215.73,56.58) .. (221.76,66.35) ;
\draw [color={rgb, 255:red, 4; green, 146; blue, 194 }  ,draw opacity=1 ]   (85.22,166.97) .. controls (91.89,177) and (228.74,76.15) .. (221.76,66.35) ;

\draw   (151.03,47.63) .. controls (150.99,43.13) and (164.71,39.42) .. (181.69,39.33) .. controls (198.66,39.24) and (212.46,42.81) .. (212.5,47.3) .. controls (212.55,51.8) and (198.83,55.51) .. (181.85,55.6) .. controls (164.88,55.69) and (151.08,52.12) .. (151.03,47.63) -- cycle ;
\draw   (295.13,148.55) .. controls (290.94,150.18) and (282.62,138.65) .. (276.55,122.79) .. controls (270.48,106.94) and (268.96,92.77) .. (273.15,91.14) .. controls (277.34,89.52) and (285.65,101.05) .. (291.72,116.9) .. controls (297.79,132.75) and (299.32,146.92) .. (295.13,148.55) -- cycle ;
\draw   (218.88,261.48) .. controls (216.55,257.64) and (226.47,247.45) .. (241.04,238.74) .. controls (255.61,230.02) and (269.3,226.08) .. (271.63,229.92) .. controls (273.95,233.77) and (264.03,243.95) .. (249.46,252.66) .. controls (234.9,261.38) and (221.2,265.33) .. (218.88,261.48) -- cycle ;
\draw   (92.73,228.12) .. controls (95.25,224.4) and (108.65,229.15) .. (122.66,238.74) .. controls (136.67,248.33) and (145.99,259.11) .. (143.47,262.83) .. controls (140.95,266.56) and (127.55,261.8) .. (113.54,252.22) .. controls (99.53,242.63) and (90.21,231.84) .. (92.73,228.12) -- cycle ;
\draw   (86.51,94.18) .. controls (90.86,95.31) and (90.83,109.53) .. (86.45,125.93) .. controls (82.07,142.33) and (75,154.71) .. (70.65,153.57) .. controls (66.3,152.43) and (66.33,138.21) .. (70.71,121.81) .. controls (75.09,105.41) and (82.16,93.04) .. (86.51,94.18) -- cycle ;
\draw    (70.65,153.57) .. controls (88.43,158.14) and (101.43,209.14) .. (92.73,228.12) ;
\draw    (143.47,262.83) .. controls (158.43,245.14) and (203.43,245.14) .. (218.88,261.48) ;
\draw    (271.63,229.92) .. controls (260.43,212.14) and (276.43,153.14) .. (295.13,148.55) ;
\draw    (86.51,94.18) .. controls (113.43,98.14) and (151.43,66.14) .. (151.03,47.63) ;
\draw    (212.5,47.3) .. controls (215.43,69.14) and (252.43,94.14) .. (273.15,91.14) ;

\draw (58.41,166.12) node [anchor=north west][inner sep=0.75pt]    {$\beta _{5}$};
\draw (288.41,163.12) node [anchor=north west][inner sep=0.75pt]    {$\beta _{1}$};
\draw (254.41,61.12) node [anchor=north west][inner sep=0.75pt]    {$\beta _{2}$};
\draw (231.41,43.12) node [anchor=north west][inner sep=0.75pt]    {$\beta _{5}$};
\draw (280.41,199.12) node [anchor=north west][inner sep=0.75pt]    {$\beta _{3}$};
\draw (119.41,46.12) node [anchor=north west][inner sep=0.75pt]    {$\beta _{1}$};
\draw (93.41,66.12) node [anchor=north west][inner sep=0.75pt]    {$\beta _{4}$};
\draw (71.41,197.12) node [anchor=north west][inner sep=0.75pt]    {$\beta _{3}$};
\draw (193.41,260.12) node [anchor=north west][inner sep=0.75pt]    {$\beta _{2}$};
\draw (156.41,260.12) node [anchor=north west][inner sep=0.75pt]    {$\beta _{4}$};
\end{tikzpicture}\]
\caption{5A}\label{fig:5A}
\end{subfigure}
\caption{3A and 5A}\label{3A5A}
\end{figure}

\item[(3S)]
The curves $\beta_1$, $\beta_2$ and $\beta_3$ in Figure~\ref{fig:3S} represent all possible pants decompositions of a surface of type $(1,1)$. If we write $S_{b_i,b_j}$ for an $S$-move mapping a curve $\beta_i \to \beta_j$, then the equation \begin{equation*}
    S_{b_1,b_2}S_{b_2,b_3}S_{b_3,b_1}=1
\end{equation*}holds in $\bS(1,1)$.

\begin{figure}[h]
\[
\begin{tikzpicture}[x=0.75pt,y=0.75pt,yscale=-0.7,xscale=0.7]

\draw [color={rgb, 255:red, 4; green, 146; blue, 194 }  ,draw opacity=1 ] [dash pattern={on 4.5pt off 4.5pt}]  (203.37,207.38) .. controls (197.93,207.27) and (197.31,134.61) .. (202.61,134.88) ;
\draw [color={rgb, 255:red, 4; green, 146; blue, 194 }  ,draw opacity=1 ]   (203.37,207.38) .. controls (208.95,207.48) and (208.19,134.82) .. (202.61,134.88) ;

\draw  [color={rgb, 255:red, 4; green, 146; blue, 194 }  ,draw opacity=1 ] (125.02,121.8) .. controls (124.86,91.31) and (161.27,66.4) .. (206.36,66.16) .. controls (251.45,65.91) and (288.13,90.43) .. (288.3,120.92) .. controls (288.46,151.41) and (252.05,176.33) .. (206.96,176.57) .. controls (161.87,176.81) and (125.19,152.29) .. (125.02,121.8) -- cycle ;
\draw   (67.84,156.99) .. controls (63.18,156.93) and (59.58,142.45) .. (59.79,124.66) .. controls (60.01,106.86) and (63.96,92.48) .. (68.61,92.54) .. controls (73.27,92.59) and (76.87,107.07) .. (76.66,124.87) .. controls (76.45,142.66) and (72.5,157.04) .. (67.84,156.99) -- cycle ;
\draw    (68.62,92.54) .. controls (99.88,92.71) and (150.78,43.86) .. (208.73,43.86) .. controls (266.68,43.86) and (313.99,89.15) .. (312.81,128.48) .. controls (311.63,167.82) and (260.77,208.34) .. (207.55,207.15) .. controls (154.33,205.96) and (93,156.49) .. (67.84,156.98) ;
\draw    (156.57,119.96) .. controls (188.95,137.43) and (207.16,139.71) .. (246.95,120.72) ;
\draw    (166.69,124.52) .. controls (192.32,110.09) and (209.18,107.05) .. (237.51,126.04) ;

\draw [color={rgb, 255:red, 4; green, 146; blue, 194 }  ,draw opacity=1 ] [dash pattern={on 4.5pt off 4.5pt}]  (166.69,124.52) .. controls (194,60.71) and (333.43,108.86) .. (302.43,159.86) ;
\draw [color={rgb, 255:red, 4; green, 146; blue, 194 }  ,draw opacity=1 ]   (166.69,124.52) .. controls (160,181.71) and (264.43,201.86) .. (302.43,159.86) ;

\draw (195.41,218.12) node [anchor=north west][inner sep=0.75pt]    {$\beta _{1}$};
\draw (97.41,111.12) node [anchor=north west][inner sep=0.75pt]    {$\beta _{2}$};
\draw (311.41,157.12) node [anchor=north west][inner sep=0.75pt]    {$\beta _{3}$};

\end{tikzpicture}\]
    \caption{3S}
    \label{fig:3S}
\end{figure}

\item[(6AS)]
Considering all possible pants decompositions of a surface of type $(1,2)$ (Figure~\ref{fig: 6AS}) the equation
\begin{equation*}
    A_{a_3,e_3}S_{a_1,a_2}A_{e_3,e_2}A_{e_2,e_1}S_{a_2,a_3}A_{e_1,a_1}=1
\end{equation*} in $\bS(1,2)$.

\begin{figure}[h]
\[
\begin{tikzpicture}[x=0.75pt,y=0.75pt,yscale=-0.8,xscale=0.8]

\draw   (30.42,72.73) .. controls (27.16,72.68) and (24.64,62.11) .. (24.79,49.11) .. controls (24.94,36.11) and (27.7,25.61) .. (30.96,25.65) .. controls (34.22,25.69) and (36.74,36.27) .. (36.59,49.27) .. controls (36.44,62.26) and (33.68,72.77) .. (30.42,72.73) -- cycle ;
\draw   (163.32,72.73) .. controls (160.07,72.68) and (157.55,62.11) .. (157.69,49.11) .. controls (157.84,36.11) and (160.61,25.61) .. (163.87,25.65) .. controls (167.12,25.69) and (169.64,36.27) .. (169.49,49.27) .. controls (169.34,62.26) and (166.58,72.77) .. (163.32,72.73) -- cycle ;
\draw    (30.96,25.65) .. controls (65.51,22.34) and (66.21,7) .. (96.28,7) .. controls (126.36,7) and (127.76,21.61) .. (163.87,25.65) ;
\draw    (163.32,72.73) .. controls (127.06,77.85) and (127.06,91.24) .. (96.98,91) .. controls (66.91,90.75) and (64.81,77.85) .. (30.42,72.73) ;
\draw    (74.3,47.72) .. controls (90.82,56.14) and (100.11,57.24) .. (120.42,48.09) ;
\draw    (77.01,48.87) .. controls (90.09,41.91) and (102.67,40.38) .. (117.12,49.53) ;

\draw    (192,50.5) -- (239,50.5) ;
\draw [shift={(241,50.5)}, rotate = 180] [color={rgb, 255:red, 0; green, 0; blue, 0 }  ][line width=0.75]    (10.93,-3.29) .. controls (6.95,-1.4) and (3.31,-0.3) .. (0,0) .. controls (3.31,0.3) and (6.95,1.4) .. (10.93,3.29)   ;
\draw [shift={(190,50.5)}, rotate = 0] [color={rgb, 255:red, 0; green, 0; blue, 0 }  ][line width=0.75]    (10.93,-3.29) .. controls (6.95,-1.4) and (3.31,-0.3) .. (0,0) .. controls (3.31,0.3) and (6.95,1.4) .. (10.93,3.29)   ;
\draw   (266.42,72.73) .. controls (263.16,72.68) and (260.64,62.11) .. (260.79,49.11) .. controls (260.94,36.11) and (263.7,25.61) .. (266.96,25.65) .. controls (270.22,25.69) and (272.74,36.27) .. (272.59,49.27) .. controls (272.44,62.26) and (269.68,72.77) .. (266.42,72.73) -- cycle ;
\draw   (399.32,72.73) .. controls (396.07,72.68) and (393.55,62.11) .. (393.69,49.11) .. controls (393.84,36.11) and (396.61,25.61) .. (399.87,25.65) .. controls (403.12,25.69) and (405.64,36.27) .. (405.49,49.27) .. controls (405.34,62.26) and (402.58,72.77) .. (399.32,72.73) -- cycle ;
\draw    (266.96,25.65) .. controls (301.51,22.34) and (302.21,7) .. (332.28,7) .. controls (362.36,7) and (363.76,21.61) .. (399.87,25.65) ;
\draw    (399.32,72.73) .. controls (363.06,77.85) and (363.06,91.24) .. (332.98,91) .. controls (302.91,90.75) and (300.81,77.85) .. (266.42,72.73) ;
\draw    (310.3,47.72) .. controls (326.82,56.14) and (336.11,57.24) .. (356.42,48.09) ;
\draw    (313.01,48.87) .. controls (326.09,41.91) and (338.67,40.38) .. (353.12,49.53) ;

\draw    (428,50.5) -- (475,50.5) ;
\draw [shift={(477,50.5)}, rotate = 180] [color={rgb, 255:red, 0; green, 0; blue, 0 }  ][line width=0.75]    (10.93,-3.29) .. controls (6.95,-1.4) and (3.31,-0.3) .. (0,0) .. controls (3.31,0.3) and (6.95,1.4) .. (10.93,3.29)   ;
\draw [shift={(426,50.5)}, rotate = 0] [color={rgb, 255:red, 0; green, 0; blue, 0 }  ][line width=0.75]    (10.93,-3.29) .. controls (6.95,-1.4) and (3.31,-0.3) .. (0,0) .. controls (3.31,0.3) and (6.95,1.4) .. (10.93,3.29)   ;
\draw   (502.42,72.73) .. controls (499.16,72.68) and (496.64,62.11) .. (496.79,49.11) .. controls (496.94,36.11) and (499.7,25.61) .. (502.96,25.65) .. controls (506.22,25.69) and (508.74,36.27) .. (508.59,49.27) .. controls (508.44,62.26) and (505.68,72.77) .. (502.42,72.73) -- cycle ;
\draw   (635.32,72.73) .. controls (632.07,72.68) and (629.55,62.11) .. (629.69,49.11) .. controls (629.84,36.11) and (632.61,25.61) .. (635.87,25.65) .. controls (639.12,25.69) and (641.64,36.27) .. (641.49,49.27) .. controls (641.34,62.26) and (638.58,72.77) .. (635.32,72.73) -- cycle ;
\draw    (502.96,25.65) .. controls (537.51,22.34) and (538.21,7) .. (568.28,7) .. controls (598.36,7) and (599.76,21.61) .. (635.87,25.65) ;
\draw    (635.32,72.73) .. controls (599.06,77.85) and (599.06,91.24) .. (568.98,91) .. controls (538.91,90.75) and (536.81,77.85) .. (502.42,72.73) ;
\draw    (546.3,47.72) .. controls (562.82,56.14) and (572.11,57.24) .. (592.42,48.09) ;
\draw    (549.01,48.87) .. controls (562.09,41.91) and (574.67,40.38) .. (589.12,49.53) ;

\draw    (95,107.5) -- (95,154.5) ;
\draw [shift={(95,156.5)}, rotate = 270] [color={rgb, 255:red, 0; green, 0; blue, 0 }  ][line width=0.75]    (10.93,-3.29) .. controls (6.95,-1.4) and (3.31,-0.3) .. (0,0) .. controls (3.31,0.3) and (6.95,1.4) .. (10.93,3.29)   ;
\draw [shift={(95,105.5)}, rotate = 90] [color={rgb, 255:red, 0; green, 0; blue, 0 }  ][line width=0.75]    (10.93,-3.29) .. controls (6.95,-1.4) and (3.31,-0.3) .. (0,0) .. controls (3.31,0.3) and (6.95,1.4) .. (10.93,3.29)   ;
\draw    (567,107.5) -- (567,154.5) ;
\draw [shift={(567,156.5)}, rotate = 270] [color={rgb, 255:red, 0; green, 0; blue, 0 }  ][line width=0.75]    (10.93,-3.29) .. controls (6.95,-1.4) and (3.31,-0.3) .. (0,0) .. controls (3.31,0.3) and (6.95,1.4) .. (10.93,3.29)   ;
\draw [shift={(567,105.5)}, rotate = 90] [color={rgb, 255:red, 0; green, 0; blue, 0 }  ][line width=0.75]    (10.93,-3.29) .. controls (6.95,-1.4) and (3.31,-0.3) .. (0,0) .. controls (3.31,0.3) and (6.95,1.4) .. (10.93,3.29)   ;
\draw   (29.42,234.73) .. controls (26.16,234.68) and (23.64,224.11) .. (23.79,211.11) .. controls (23.94,198.11) and (26.7,187.61) .. (29.96,187.65) .. controls (33.22,187.69) and (35.74,198.27) .. (35.59,211.27) .. controls (35.44,224.26) and (32.68,234.77) .. (29.42,234.73) -- cycle ;
\draw   (162.32,234.73) .. controls (159.07,234.68) and (156.55,224.11) .. (156.69,211.11) .. controls (156.84,198.11) and (159.61,187.61) .. (162.87,187.65) .. controls (166.12,187.69) and (168.64,198.27) .. (168.49,211.27) .. controls (168.34,224.26) and (165.58,234.77) .. (162.32,234.73) -- cycle ;
\draw    (29.96,187.65) .. controls (64.51,184.34) and (65.21,169) .. (95.28,169) .. controls (125.36,169) and (126.76,183.61) .. (162.87,187.65) ;
\draw    (162.32,234.73) .. controls (126.06,239.85) and (126.06,253.24) .. (95.98,253) .. controls (65.91,252.75) and (63.81,239.85) .. (29.42,234.73) ;
\draw    (73.3,209.72) .. controls (89.82,218.14) and (99.11,219.24) .. (119.42,210.09) ;
\draw    (76.01,210.87) .. controls (89.09,203.91) and (101.67,202.38) .. (116.12,211.53) ;

\draw    (191,212.5) -- (238,212.5) ;
\draw [shift={(240,212.5)}, rotate = 180] [color={rgb, 255:red, 0; green, 0; blue, 0 }  ][line width=0.75]    (10.93,-3.29) .. controls (6.95,-1.4) and (3.31,-0.3) .. (0,0) .. controls (3.31,0.3) and (6.95,1.4) .. (10.93,3.29)   ;
\draw [shift={(189,212.5)}, rotate = 0] [color={rgb, 255:red, 0; green, 0; blue, 0 }  ][line width=0.75]    (10.93,-3.29) .. controls (6.95,-1.4) and (3.31,-0.3) .. (0,0) .. controls (3.31,0.3) and (6.95,1.4) .. (10.93,3.29)   ;
\draw   (265.42,234.73) .. controls (262.16,234.68) and (259.64,224.11) .. (259.79,211.11) .. controls (259.94,198.11) and (262.7,187.61) .. (265.96,187.65) .. controls (269.22,187.69) and (271.74,198.27) .. (271.59,211.27) .. controls (271.44,224.26) and (268.68,234.77) .. (265.42,234.73) -- cycle ;
\draw   (398.32,234.73) .. controls (395.07,234.68) and (392.55,224.11) .. (392.69,211.11) .. controls (392.84,198.11) and (395.61,187.61) .. (398.87,187.65) .. controls (402.12,187.69) and (404.64,198.27) .. (404.49,211.27) .. controls (404.34,224.26) and (401.58,234.77) .. (398.32,234.73) -- cycle ;
\draw    (265.96,187.65) .. controls (300.51,184.34) and (301.21,169) .. (331.28,169) .. controls (361.36,169) and (362.76,183.61) .. (398.87,187.65) ;
\draw    (398.32,234.73) .. controls (362.06,239.85) and (362.06,253.24) .. (331.98,253) .. controls (301.91,252.75) and (299.81,239.85) .. (265.42,234.73) ;
\draw    (309.3,209.72) .. controls (325.82,218.14) and (335.11,219.24) .. (355.42,210.09) ;
\draw    (312.01,210.87) .. controls (325.09,203.91) and (337.67,202.38) .. (352.12,211.53) ;

\draw    (427,212.5) -- (474,212.5) ;
\draw [shift={(476,212.5)}, rotate = 180] [color={rgb, 255:red, 0; green, 0; blue, 0 }  ][line width=0.75]    (10.93,-3.29) .. controls (6.95,-1.4) and (3.31,-0.3) .. (0,0) .. controls (3.31,0.3) and (6.95,1.4) .. (10.93,3.29)   ;
\draw [shift={(425,212.5)}, rotate = 0] [color={rgb, 255:red, 0; green, 0; blue, 0 }  ][line width=0.75]    (10.93,-3.29) .. controls (6.95,-1.4) and (3.31,-0.3) .. (0,0) .. controls (3.31,0.3) and (6.95,1.4) .. (10.93,3.29)   ;
\draw   (501.42,234.73) .. controls (498.16,234.68) and (495.64,224.11) .. (495.79,211.11) .. controls (495.94,198.11) and (498.7,187.61) .. (501.96,187.65) .. controls (505.22,187.69) and (507.74,198.27) .. (507.59,211.27) .. controls (507.44,224.26) and (504.68,234.77) .. (501.42,234.73) -- cycle ;
\draw   (634.32,234.73) .. controls (631.07,234.68) and (628.55,224.11) .. (628.69,211.11) .. controls (628.84,198.11) and (631.61,187.61) .. (634.87,187.65) .. controls (638.12,187.69) and (640.64,198.27) .. (640.49,211.27) .. controls (640.34,224.26) and (637.58,234.77) .. (634.32,234.73) -- cycle ;
\draw    (501.96,187.65) .. controls (536.51,184.34) and (537.21,169) .. (567.28,169) .. controls (597.36,169) and (598.76,183.61) .. (634.87,187.65) ;
\draw    (634.32,234.73) .. controls (598.06,239.85) and (598.06,253.24) .. (567.98,253) .. controls (537.91,252.75) and (535.81,239.85) .. (501.42,234.73) ;
\draw    (545.3,209.72) .. controls (561.82,218.14) and (571.11,219.24) .. (591.42,210.09) ;
\draw    (548.01,210.87) .. controls (561.09,203.91) and (573.67,202.38) .. (588.12,211.53) ;

\draw [color={rgb, 255:red, 4; green, 146; blue, 194 }  ,draw opacity=1 ] [dash pattern={on 4.5pt off 4.5pt}]  (96.37,91.38) .. controls (90.93,91.35) and (90.69,54.89) .. (95.99,54.99) ;
\draw [color={rgb, 255:red, 4; green, 146; blue, 194 }  ,draw opacity=1 ]   (96.37,91.38) .. controls (101.95,91.4) and (101.57,54.93) .. (95.99,54.99) ;

\draw [color={rgb, 255:red, 4; green, 146; blue, 194 }  ,draw opacity=1 ] [dash pattern={on 4.5pt off 4.5pt}]  (332.37,91.38) .. controls (326.93,91.35) and (326.69,54.89) .. (331.99,54.99) ;
\draw [color={rgb, 255:red, 4; green, 146; blue, 194 }  ,draw opacity=1 ]   (332.37,91.38) .. controls (337.95,91.4) and (337.57,54.93) .. (331.99,54.99) ;

\draw [color={rgb, 255:red, 4; green, 146; blue, 194 }  ,draw opacity=1 ] [dash pattern={on 4.5pt off 4.5pt}]  (332.37,43.38) .. controls (326.93,43.35) and (326.69,6.89) .. (331.99,6.99) ;
\draw [color={rgb, 255:red, 4; green, 146; blue, 194 }  ,draw opacity=1 ]   (332.37,43.38) .. controls (337.95,43.4) and (337.57,6.93) .. (331.99,6.99) ;

\draw [color={rgb, 255:red, 4; green, 146; blue, 194 }  ,draw opacity=1 ] [dash pattern={on 4.5pt off 4.5pt}]  (568.37,43.38) .. controls (562.93,43.35) and (562.69,6.89) .. (567.99,6.99) ;
\draw [color={rgb, 255:red, 4; green, 146; blue, 194 }  ,draw opacity=1 ]   (568.37,43.38) .. controls (573.95,43.4) and (573.57,6.93) .. (567.99,6.99) ;

\draw [color={rgb, 255:red, 4; green, 146; blue, 194 }  ,draw opacity=1 ]   (140.67,78.67) .. controls (145.2,52.2) and (142.67,22.67) .. (93.67,22.67) .. controls (44.67,22.67) and (60.2,60.2) .. (49.67,77.67) ;
\draw [color={rgb, 255:red, 4; green, 146; blue, 194 }  ,draw opacity=1 ] [dash pattern={on 4.5pt off 4.5pt}]  (140.67,78.67) .. controls (128.31,67.23) and (147.2,28.6) .. (98.2,28.4) .. controls (49.2,28.2) and (43.2,56.2) .. (49.67,77.67) ;
\draw [color={rgb, 255:red, 4; green, 146; blue, 194 }  ,draw opacity=1 ]   (139.67,240.67) .. controls (144.2,214.2) and (141.67,184.67) .. (92.67,184.67) .. controls (43.67,184.67) and (59.2,222.2) .. (48.67,239.67) ;
\draw [color={rgb, 255:red, 4; green, 146; blue, 194 }  ,draw opacity=1 ] [dash pattern={on 4.5pt off 4.5pt}]  (139.67,240.67) .. controls (127.31,229.23) and (146.2,190.6) .. (97.2,190.4) .. controls (48.2,190.2) and (42.2,218.2) .. (48.67,239.67) ;
\draw [color={rgb, 255:red, 4; green, 146; blue, 194 }  ,draw opacity=1 ]   (521,22.6) .. controls (516.39,49.05) and (520,76.52) .. (569,76.67) .. controls (618,76.82) and (601.41,37.04) .. (612,19.6) ;
\draw [color={rgb, 255:red, 4; green, 146; blue, 194 }  ,draw opacity=1 ] [dash pattern={on 4.5pt off 4.5pt}]  (521,22.6) .. controls (533.32,34.08) and (515.48,70.58) .. (564.48,70.92) .. controls (613.48,71.27) and (618.4,41.09) .. (612,19.6) ;
\draw [color={rgb, 255:red, 4; green, 146; blue, 194 }  ,draw opacity=1 ]   (520,184.6) .. controls (515.39,211.05) and (519,238.52) .. (568,238.67) .. controls (617,238.82) and (600.41,199.04) .. (611,181.6) ;
\draw [color={rgb, 255:red, 4; green, 146; blue, 194 }  ,draw opacity=1 ] [dash pattern={on 4.5pt off 4.5pt}]  (520,184.6) .. controls (532.32,196.08) and (514.48,232.58) .. (563.48,232.92) .. controls (612.48,233.27) and (617.4,203.09) .. (611,181.6) ;
\draw  [color={rgb, 255:red, 4; green, 146; blue, 194 }  ,draw opacity=1 ] (66.12,211.24) .. controls (66.09,204.65) and (79.07,199.24) .. (95.13,199.16) .. controls (111.18,199.07) and (124.23,204.34) .. (124.26,210.92) .. controls (124.3,217.51) and (111.31,222.92) .. (95.26,223.01) .. controls (79.2,223.09) and (66.16,217.82) .. (66.12,211.24) -- cycle ;
\draw  [color={rgb, 255:red, 4; green, 146; blue, 194 }  ,draw opacity=1 ] (302.12,211.24) .. controls (302.09,204.65) and (315.07,199.24) .. (331.13,199.16) .. controls (347.18,199.07) and (360.23,204.34) .. (360.26,210.92) .. controls (360.3,217.51) and (347.31,222.92) .. (331.26,223.01) .. controls (315.2,223.09) and (302.16,217.82) .. (302.12,211.24) -- cycle ;
\draw  [color={rgb, 255:red, 4; green, 146; blue, 194 }  ,draw opacity=1 ] (538.12,211.24) .. controls (538.09,204.65) and (551.07,199.24) .. (567.13,199.16) .. controls (583.18,199.07) and (596.23,204.34) .. (596.26,210.92) .. controls (596.3,217.51) and (583.31,222.92) .. (567.26,223.01) .. controls (551.2,223.09) and (538.16,217.82) .. (538.12,211.24) -- cycle ;
\draw [color={rgb, 255:red, 4; green, 146; blue, 194 }  ,draw opacity=1 ]   (289.6,240) .. controls (318.2,246.8) and (342.2,240.4) .. (360.2,231.4) .. controls (378.2,222.4) and (378.2,202.4) .. (360.2,190.4) .. controls (342.2,178.4) and (322.6,173) .. (289.6,182) ;
\draw [color={rgb, 255:red, 4; green, 146; blue, 194 }  ,draw opacity=1 ] [dash pattern={on 4.5pt off 4.5pt}]  (289.6,182) .. controls (279,185.6) and (278.6,198) .. (278.6,212) .. controls (278.6,227) and (282.6,233) .. (289.6,240) ;

\draw (104.41,61.2) node [anchor=north west][inner sep=0.75pt]    {$\alpha _{3}$};
\draw (342.41,60.2) node [anchor=north west][inner sep=0.75pt]    {$\alpha _{3}$};
\draw (343.41,15.2) node [anchor=north west][inner sep=0.75pt]    {$\alpha _{1}$};
\draw (579.41,14.4) node [anchor=north west][inner sep=0.75pt]    {$\alpha _{1}$};
\draw (86.41,223.2) node [anchor=north west][inner sep=0.75pt]    {$\alpha _{2}$};
\draw (291.6,185.4) node [anchor=north west][inner sep=0.75pt]    {$\alpha _{2}$};
\draw (560.41,178.2) node [anchor=north west][inner sep=0.75pt]    {$\alpha _{2}$};
\draw (141.67,76.07) node [anchor=north west][inner sep=0.75pt]    {$\varepsilon _{1}$};
\draw (137.67,240.07) node [anchor=north west][inner sep=0.75pt]    {$\varepsilon _{1}$};
\draw (276.67,242.07) node [anchor=north west][inner sep=0.75pt]    {$\varepsilon _{2}$};
\draw (508.67,162.07) node [anchor=north west][inner sep=0.75pt]    {$\varepsilon _{3}$};
\draw (510.67,-0.93) node [anchor=north west][inner sep=0.75pt]    {$\varepsilon _{3}$};
\draw (210,30) node [anchor=north west][inner sep=0.75pt]   [align=left] {A};
\draw (446,30) node [anchor=north west][inner sep=0.75pt]   [align=left] {A};
\draw (209,192) node [anchor=north west][inner sep=0.75pt]   [align=left] {A};
\draw (445,192) node [anchor=north west][inner sep=0.75pt]   [align=left] {A};
\draw (104,122) node [anchor=north west][inner sep=0.75pt]   [align=left] {S};
\draw (576,122) node [anchor=north west][inner sep=0.75pt]   [align=left] {S};
\end{tikzpicture}
\]
\caption{3AS}\label{fig: 6AS}
\end{figure}

\item[(C)] Any two moves $\alpha_1 \to \alpha_2$ and $\beta_1 \to \beta_2$ (either S-moves or A-moves) supported in disjoint subsurfaces commute. 
\end{itemize}

\subsection{The Grothendieck-Teichm\"uller and Nakumara-Schneps groups}\label{sec:gt}
Recall from our introduction that the goal of this series is to convince the reader that the modular operad of seamed surfaces, $B\bS$, gives a reasonable model for the Teichm\"uller tower (after completion). In order to argue this, we need to show that $(1)$ the absolute Galois group acts on our proposed model and $(2)$ that this action commutes with the modular operad structure.  In practice, we don't really know how to do this without passing through some intermediate profinte groups. Here, we introduce two such groups: Grothendieck-Teichm\"uller and Nakumara-Schneps groups. 

\medskip

Let $\widehat{F}_2$ denote the profinite completion of the free group on two letters $F_2=<x,y>$. If the reader still finds the abstract definition of profinite group confusing, it can be helpful to think of an element $f\in\widehat{F}_2$ as a (possibly infinite) word in $x$ and $y$. Any homomorphism of profinite groups \[\begin{tikzcd} \widehat{F}_2\arrow[r]& \widehat{G} \end{tikzcd}\] will necessarily be determined by where it sends the generators, $(x,y)\mapsto (a,b)$, and we will write $f(a,b)$ for the image of a word $f\in \widehat{F}_2$ in $\widehat{G}$. For example, we will write $f(y,x)$ for the image of any $f\in \widehat{F}_2$ under the map $ \widehat{F}_2\rightarrow \widehat{F}_2$ given by $(x,y)\mapsto (y,x)$.

\begin{definition}\label{def: gt} The Grothendieck-Teichm\"uller group $\gt$ is the group of pairs $$(\lambda, f) \in \widehat{\mathbb{Z}}^{*}\times \widehat{F}'_2$$ which satisfy the property that \[ x\mapsto x^{\lambda} \quad \text{and} \quad y \mapsto f^{-1}y^{\lambda} f\] induces an automorphism of $\widehat{F}_2$. Moreover, we require the pair $(\lambda, f)$ satisfy the following axioms: 
\begin{itemize} 
\item[(I)] $f(x,y)f(y,x)=1$,
\item[(II)] $f(x,y)x^mf(z,x)z^mf(y,z)y^m=1$ where $xyz=1$ and $m=(\lambda-1)/2$,
\item[(III)] $f(b_{3},b_{4})f(b_{5},b_{1})f(b_{2},b_{3})f(b_{4},b_{5})f(b_{1},b_{2})=1$ in $\widehat{\Gamma}_{0,5}$ where $b_{i}$ is a \emph{Dehn twist} along a loop $\beta_i$ depicted in Figure~\ref{fig:5A}. 
\end{itemize}
\end{definition}

The profinite \emph{Grothendieck-Teichm\"uller group}, $\gt$, is closely related to the absolute Galois group. In particular, a theorem of Ihara says:

\begin{theorem}\cite{ihara} There is an injection $Gal(\mathbb{Q})\hookrightarrow \gt$. \end{theorem} 

A related group, defined by Nakumara and Schneps in \cite{ns} is defined by adding a ``higher genus'' relation to $\gt$. 

\begin{definition}\label{def: ns} Let $\ns$ denote the group of pairs $$(\lambda, f) \in \widehat{\mathbb{Z}}^{*}\times \widehat{F}'_2$$ which satisfy the property that \[ x\mapsto x^{\lambda} \quad \text{and} \quad y \mapsto f^{-1}y^{\lambda} f\] induce an automorphism of $\widehat{F}_2$. Moreover, we require pairs $(\lambda, f)$ satisfy relations $(I)-(III)$ of $\gt$ and : 
\begin{itemize} 
\item[(IV)]$ f(e_1, a_1)a_3^{-8\rho_2}f(a_2^2,a_3^2)(a_3a_2a_3)^{2m}f(e_2,e_1)e_2^{2m} f(e_3,e_2)a_2^{-2m}(a_1a_2a_1)^{2m}f(a_1^{2},a_2^2)a_1^{8\rho_2}f(a_3,e_3)=1$ where $a_1,a_2,a_3, e_1,e_2$ are Dehn twists in $\Gamma_{1,2}$ corresponding to the curves in Figure~\ref{fig: 6AS}.
\end{itemize} 
\end{definition} 

Nakumara and Schneps show that $\ns$ is a subgroup of $\gt$ and that Ihara's injection $\Gal(\mathbb{Q}) \hookrightarrow \widehat{GT}$ also maps $\Gal(\mathbb{Q})$ into $\ns$ (Theorem 1.2 \cite{ns}). 

\begin{remark}
At this time, it is not known if $\ns$ is a \emph{proper} subgroup of $\gt$.
\end{remark}

The goal for the remainder of this lecture is to investigate the Galois actions on our proposed model for the Teichm\"uller tower by studying actions of $\ns$ and $\gt$. 

\subsubsection{The genus zero case}
The Grothendieck-Teichm\"uller group is closely related to the \emph{genus zero} component of our tower. Let $\bS_{0}=\{\bS(0,n)\}$ for the restriction of $\bS$ to genus $0$. In other words, $\bS_0$ is the \emph{underlying cyclic operad} of $\bS$ via the adjunction in \eqref{adjunction}. The underlying \emph{operad} $\calS=\{\calS(0,n+1)\}$ is obtained from the cyclic operad $\bS_{0}=\{\bS(0,n)\}$ by marking one boundary of each surface as the distinguished output of the surface. Operad composition is then defined by gluing the \emph{marked} boundary component of a surface in $\calS(0,m+1)$ to the $i$th free boundary component of $\calS(0,n+1)$:\[\begin{tikzcd} \calS(0,n+1)\times \calS(0,m+1) \arrow[r, "\circ_i"] & \calS(0, n+m+1).\end{tikzcd}\]

\begin{remark}
As we mentioned, the modular operad $\bS=\{\bS(g,n)\}$ is an extension of the surface operad of Tillmann~\cite{Till} and Wahl~\cite[Section 3.1]{W2}. Our genus zero operad $\calS=\{\calS(0,n+1)\}$ is precisely the genus zero part of Tillmann's operad (\cite[Definition 6.5]{bhr1}).  This operad is equivalent to the operad of parenthesized ribbon braids $\mathsf{PaRB} =\{\mathsf{PaRB}(n)\}$ and closely related to the operad of framed discs, in the sense that there are homotopy equivalences $$B\calS(0,n+1) \simeq B(\mathsf{PaRB}(n)) \simeq \mathsf{fD}(n).$$
\end{remark}

Each groupoid $\calS(0, n+1)$ in our operad $\calS$ has finitely many objects and thus, after applying the profinite completion functor entrywise, we obtain an operad in profinite groupoids $$\widehat{\calS}=\{\widehat{\calS}(0, n+1).\}$$ We write $\End_{0}$ for the set of operad endomorphisms which fix objects. Proposition 7.3 and Proposition 8.1 of \cite{bhr1} combine to show:  
\begin{prop}\label{prop:gt} There is an isomorphism
\[\gt \cong \textnormal{End}_{0}(\widehat{\calS}).\] 
\end{prop}

Putting this together with the nerve theorem from the first lecture (Theorem~\ref{thm: nerve}), we identify the group $\gt$ with the group of (path components of) self maps of the $\infty$-operad $N\widehat{B\calS}$ (\cite[Theorem 8.4]{bhr1}): 
\begin{theorem}
There is an isomorphism \[\gt  \cong \pi_0\mathbb{R}\End (N\widehat{B\calS}).\]
\end{theorem} 

\begin{remark}
In \cite{bhr1}, we actually show $\gt \cong \textnormal{End}_{0}(\widehat{\mathsf{PaRB}})$. We have been a bit loose with the translation, because one can show that the operads $\mathsf{PaRB}$ and  $\calS$ are equivalent. This presentation just translates a bit easier to the higher genus case. 
\end{remark}

\subsubsection{The $\gt$ action}\label{sec: gt action}
The proof of Proposition~\ref{prop:gt} is outside of the scope of these lectures. However, we can describe the arity-wise action of $\gt$ 
\[\begin{tikzcd} \gt \arrow[r] &  \textnormal{End}_{0}(\widehat{\bS}(0,n)) \end{tikzcd}\] on the profinte cyclic operad rather easily by translating the action of $\gt$ on the $\mathcal{SHT}$ complex from \cite{ns} to our groupoids. 

\medskip

Recall that an object $\Sigma$ in one of the groupoids $\bS(0,n)$ is equivalent to fixing a surface of type $(0,n)$ together with an ``atomic'' quilted pants decomposition. For each $(\lambda, f)\in\gt$, we define an $F_{(\lambda,f)}:\bS(0,n)\rightarrow\bS(0,n)$ which is the identity on objects and acts on elementary morphisms by:
\begin{equation}\label{gt action}
\begin{tikzcd} (\lambda, f) \arrow[r] & \begin{cases} a_{\frac{1}{2}} \mapsto a^{\lambda}_{\frac{1}{2}} \\ A_{\alpha,b}\mapsto A_{\alpha,b} \cdot f(a,b)a^{n(\lambda-1)/2}. \end{cases} \end{tikzcd} \end{equation} Here $a_{\frac{1}{2}}$ is a half Dehn twist around the boundary components or any curve $\{\alpha_i\}$ in the pants decomposition of the relevant object. The integer $n$ which arises in the action on an A-move $A_{\alpha,\beta}:\Sigma\rightarrow\Sigma'$ can be calculated based on the interaction of the A-move with the quilt on $\Sigma$. We wont discuss how this integer is computed in these notes, but a similar computation is done in Section 7 of \cite{bhr1}. Full details will appear in \cite{BR22}. 

\begin{remark}
Note that there are no S-moves in the genus zero groupoids $\bS(0,n)$. 
\end{remark}


\medskip
To check that the action of $\gt$ in \eqref{gt action} is well-defined, we need to check that the maps $F_{(\lambda,f)}$ commute with the defining relations of $\gt$. In general, this is a bit involved, but we can sketch how one shows that the map $F_{(\lambda,f)}$ commutes with relation (III) from Definition~\ref{def: gt}. We fix a surface $\Sigma$ of type $(0,5)$. Then the action of $\gt$ on the automorphism of $\Sigma$ given by the composite of A-moves $$A_{b_3,b_4} A_{b_1,b_2}A_{b_4,b_5}A_{b_2,b_3}A_{b_5,b_1}$$ becomes: 
\begin{equation}\label{eq 1}
  A_{b_3,b_4} A_{b_1,b_2}A_{b_4,b_5}A_{b_2,b_3}A_{b_5,b_1} \mapsto   A_{b_3,b_4} f(b_1,b_2)b_1^{n(\lambda-1)/2} \ldots A_{b_5,b_1} f(b_5,b_1)b_5^{n(\lambda-1)/2}.
\end{equation}

A quick computation shows that for this particular action that $n=0$ (A proof is similar to \cite[Proposition 8.3]{ns}. Full details in this setting will appear in \cite{BR22}). The action commutes with categorical composition, and so Equation \eqref{eq 1}, simplifies to  
\begin{equation}\label{eq 2}
 A_{b_3,b_4} f(b_3,b_4) \ldots A_{b_5,b_1} f(b_5,b_1) = (A_{b_3,b_4} \ldots A_{b_5, b_1})\cdot (f(b_3,b_4) \ldots f(b_5,b_1)) =1.
\end{equation}
 
But now relation $(5A)$ between morphisms in $\bS(0,5)$, reduces Equation \eqref{eq 2} to \[f(b_3,b_4) f(b_1,b_2)f(b_4,b_5)f(b_2,b_3)f(b_5,b_1) =1.\] This is precisely relation (III) in the definition of $\gt$. 

The other relations follow a similar pattern. The difficult part is showing that the $\gt$-action commutes with the cyclic operad structure maps. As we did in \cite{bhr1}, we can overcome this by showing: 

\begin{prop} There is an isomorphism
\[\gt \cong \textnormal{End}_{0}(\widehat{\bS}_{0}).\] 
\end{prop}

\subsubsection{The higher genus action}
The Nakamura-Schneps group $\ns$ acts on the full modular operad $\bS$ in such a way that the restriction to genus zero is precisely the action of $\gt$ on $\bS_0$ we have just described. In \cite{BR22} we show: 

\begin{prop}\cite{BR22}\label{ns is end}
The profinte group $\ns$ acts on $\textnormal{End}_{0}(\widehat{\bS}).$
\end{prop}

As in the genus zero case, we can described the arity wise action \[\begin{tikzcd} \ns \arrow[r] &  \textnormal{End}_{0}(\widehat{\bS}(g,n)). \end{tikzcd}\] 
 
Given a $(\lambda, f)\in\ns$ we wish to define a functor $F_{(\lambda, f)}:\bS(g,n)\rightarrow \bS(g,n)$ which fixes objects and acts on elementary morphisms via
\[\begin{tikzcd} (\lambda, f) \arrow[r,"F_{(\lambda, f)}"] & \begin{cases} a_{\frac{1}{2}} \mapsto a^{\lambda}_{\frac{1}{2}} \\ A_{\alpha,b}\mapsto A_{\alpha,b} \cdot f(a,b)a^{n(\lambda-1)/2}\\ S_{\alpha, b} \mapsto S_{\alpha, b}\cdot (aba)^{\lambda-1}b^{n(\lambda -1)/2-8\rho_2} f(a^2,b^2)a^{8\rho_2}. \end{cases} \end{tikzcd} \] 
As before the, if $S_{a,b}:\Sigma\rightarrow\Sigma'$ is our S-move, the integer $n$ is dependent on the quilt of $\Sigma$. The integer $\rho_2$ is the Kummer $1$-cocycle with respect to the roots of $2$ (See Section 5 of \cite{ns} for full details). 

To check that the action we have given is well-defined, it remains to verify that it is compatible with the defining relations of our group $\ns$. For example, if we consider the $\ns$ action on $\bS(1,2)$ then we know that all $A$ and $S$ moves necessarily satisfy the relation (6AS):
\begin{equation}
    A_{\alpha_3,\epsilon_3}S_{\alpha_1,\alpha_2}A_{\epsilon_3,\epsilon_2}A_{\epsilon_2,\epsilon_1}S_{\alpha_2,\alpha_3}A_{\epsilon_1,\alpha_1}=1. 
\end{equation}
Acting by $(\lambda, f)$ gives the equation
\begin{equation}
A_{\alpha_3,\epsilon_3} f(a_3,e_3)a_3^{n(\lambda-1)/2} S_{\alpha_1,\alpha_2} (a_1a_2a_1)^{\lambda-1}a_2^{n(\lambda -1)/2-8\rho_2} f(a_1^2,a_2^2)a_1^{8\rho_2}\ldots A_{\epsilon_1,\alpha_1}f(e_1,a_1)e_1^{n(\lambda-1)/2}.
\end{equation}

The integer $n$ is computed based on the quilting (see \cite[Proposition 8.3]{ns}) and the action, coming from a group homomorphism, commutes with the categorical compositions. This yields the relation (IV) from Definition~\ref{def: ns}: 

\begin{equation*} 
f(e_1, a_1)a_3^{-8\rho_2}f(a_2^2,a_3^2)(a_3a_2a_3)^{2m}f(e_2,e_1)e_2^{2m} f(e_3,e_2)a_2^{-2m}(a_1a_2a_1)^{2m}f(a_1^{2},a_2^2)a_1^{8\rho_2}f(a_3,e_3)=1
\end{equation*} 

\subsubsection{An operadic two level principle}
The groupoids $\bS(g,n)$ are only homotopy approximations of 
$\Gamma_{g}^{n}$ in the sense that \[B\bS(g,n)\simeq B\Gamma_{g}^{n}.\] 
Therefore, in order to complete our description of Teichm\"uller tower we 
want to see the action of $\ns$ on the profinite completion of the modular 
operad $B\bS$. Applying the profinite completion functor \eqref{completion 
adjunction space} entrywise results in a sequence of profinite spaces 
\[\widehat{B\bS} =\{\widehat{B\bS(g,n)},\}\] but, unfortunately, these 
profinite spaces do not form a modular $\infty$-operad. This is because we 
do not know if the mapping class groups $\Gamma_{g}^{n}$ are good groups 
(Definition~\ref{def: good}) for $g\geq 2$ and thus we cannot apply 
Proposition~\ref{prof spaces} to get a family of weak composition maps.  
One can show, however, that we have a modular dendroidal space: 
\[N\widehat{B\bS}: \bU^{op}\rightarrow\widehat{\textbf{sSet}}.\] 

The truncation of the modular operad $\bS$ at genus $1$ to defines a modular $B\bS_{1}$ with \[B\bS_{1}(g,n) =\begin{cases} B\bS(g,n) \ \text{if} \ g\leq 1\\
\emptyset \ \text{otherwise.}\end{cases}\]In this case, applying the profinte completion functor entrywise results in a modular $\infty$-operad, $N\widehat{B\bS}_{1}$.

The modular operad $\bS$ is generated by a single object (our standard pair of pants) and morphisms in genus zero and one. It follows that one can show: 
\begin{prop}
There is an isomorphism of profinite groups: 
\[\End_{0}(\widehat{\bS})\cong \End_{0}(\widehat{\bS}_{1}).\]
\end{prop}

The classifying space functor $B:\mathbf{Gpd}\rightarrow\mathbf{sSet}$ is homotopically fully faithful, meaning that for any two groupoids $\mathbf{C}$ and $\mathbf{D}$: \[\map(\mathbf{C},\mathbf{D})\cong \map(B\mathbf{C},B\mathbf{D}).\] Combining this with the observation that the truncation functor 
\[\begin{tikzcd} \ModOp_{1}(\widehat{\textbf{Gpd}})  \arrow[r, shift left = .15cm] & \ModOp(\widehat{\textbf{Gpd}})  \arrow[l, shift left=.15cm, "\tau_{1}"] \end{tikzcd}\]is part of a Quillen adjunction leads us to our final theorem of this lecture series:
\begin{theorem}\label{thm:truncated}\cite{BR22}
There is an action of the profinite group $\ns$ on the profinite modular $\infty$-operad $N\widehat{B\bS}$.  
\end{theorem}

\bibliographystyle{amsalpha}
\bibliography{Barcelona}

\end{document}